\documentclass[aos,MSNbibl,rotating,dvips]{arximspdf}

\doi{10.1214/12-AOS1056} 
\volume{40}
\issue{6}
\pubyear{2012}
\firstpage{2943}
\lastpage{2972}

\makeatletter
\newcommand{\rrVert}{\Vert}
\newcommand{\rrvert}{\vert}
\newcommand{\llVert}{\Vert}
\newcommand{\llvert}{\vert}

\newcommand{\cal}{\mathcal}
\newcommand{\calS}{\mathcal{S}}
\newcommand{\calD}{\mathcal{D}}

\newcommand{\statT}{\mathrm{T}}

\newcommand{\ryskus}{\mathbf}

\newcommand{\bmu}{\bolds\mu}

\newcommand{\btheta}{\bolds\theta}

\newcommand{\bSigma}{\bolds\Sigma}
\newcommand{\balpha}{\bolds\alpha}
\newcommand{\bvarepsilon}{\bolds\varepsilon}
\newcommand{\bbeta}{\bolds\beta}
\newcommand{\bldeta}{\bolds\eta}
\newcommand{\bgamma}{\bolds\gamma}

\newtheorem{theorem}{Theorem}
\newtheorem{lemma}{Lemma}

\newproclaim{remark}{Remark}
\newproclaim{example}{Example}

\makeatother

\begin{document}
\begin{frontmatter}

\title{Two-step spline estimating equations for generalized additive
partially linear models with~large cluster sizes}
\runtitle{Spline GEE for GAPLMS with large cluster sizes}

\begin{aug}
\author[A]{\fnms{Shujie} \snm{Ma}\corref{}\ead[label=e1]{shujie.ma@ucr.edu}\ead[label=u1,url]{http://faculty.ucr.edu/\textasciitilde shujiema/}}
\runauthor{S. Ma}
\affiliation{University of California, Riverside}
\address[A]{Department of Statistics\\
University of California, Riverside\\
Riverside, California 92521\\
USA\\
\printead{e1} \\
\printead{u1}} 
\end{aug}

\received{\smonth{1} \syear{2012}}
\revised{\smonth{6} \syear{2012}}

%
\begin{abstract}
We propose a two-step estimating procedure for generalized additive
partially linear models with clustered data using estimating equations.
Our proposed method applies to the case that the number of observations
per cluster is allowed to increase with the number of independent
subjects. We establish oracle properties for the two-step estimator of
each function component such that it performs as well as the univariate
function estimator by assuming that the parametric vector and all other
function components are known. Asymptotic distributions and consistency
properties of the estimators are obtained. Finite-sample experiments
with both simulated continuous and binary response variables confirm
the asymptotic results. We illustrate the methods with an application
to a U.S. unemployment data set.
\end{abstract}

%

\end{frontmatter}

\section{Introduction}\label{SECIntroduction}

The generalized estimating equations (GEE) approach has been widely applied
to the analysis of clustered data. Reference~\cite{LZ86} introduced the GEEs to
estimate the regression parameters of generalized linear models with
possible unknown correlations between responses. The GEE approach only
requires the first two marginal moments and a working correlation matrix
that accounts for the form of within-subject correlations of responses, and
it can yield consistent parameter estimators even when the covariance
structure is misspecified, as long as the mean function is correctly
specified.

Parametric GEEs enjoy simplicity by assuming a fully predetermined
parametric form for the mean function, but they have suffered from
inflexibility in modeling complicated relationships between the
response and covariates in clustered data studies. To allow for
flexibility,~\cite{WY96,HRWY98} and~\cite{LC00} proposed to model
covariate effects nonparametrically via GEE. The proposed nonparametric
GEE method enables us to capture the underlying structure that
otherwise can be missed. Reference~\cite{LC01} extended the kernel
estimating equations in~\cite{LC00} to generalized partially linear
models (GPLMs), which assume that the mean of the outcome variable
depends on a vector of covariates parametrically and a scalar predictor
nonparametrically to overcome the ``curse of dimensionality'' of
nonparametric models. As an extension,~\cite{HFZ05} and~\cite{HZZ06}
approximated the nonparametric function in GPLMs by regression splines.
It is pointed out in~\cite{WLC02} and~\cite{LWWC04} that splines
effectively account for the correlations of clustered data and are more
efficient in nonparametric models with longitudinal data than
conventional local-polynomials. Splines also provide optimal
convergence rates in partially linear models~\cite{H86,HS96}. To allow
the nonparametric part in partially linear models to include
multivariate covariates,~\cite{MSW12} extended the estimating equations
method to generalized additive partially linear models (GAPLMs) with an
identity link for continuous response cases, and obtained estimators
for the parametric vector and the nonparametric additive functions via
a one-step spline estimation.

To introduce GAPLMs for clustered data, denote $\{(Y_{ij},\ryskus
{X}_{ij},%
\ryskus{Z}_{ij}),1\leq i\leq n,1\leq j\leq m_{i}\}$ as the $j$th
repeated observation for the $i$th subject or experimental unit, where
$Y_{ij}$ is the response variable, $\ryskus{X}_{ij}= (
1,X_{ij1},\ldots,X_{ij(d_{1}-1)} )^{\statT}$ and
$\ryskus{Z}_{ij}= (
Z_{ij1},\ldots,Z_{ijd_{2}} )^{\statT}$ are $d_{1}$-dimensional and $%
d_{2}$-dimensional vectors of covariates, respectively. The marginal model
assumes that $Y_{ij}=\mu_{ij}+\varepsilon_{ij}$, and the marginal mean $
\mu_{ij}$ depends on $\ryskus{X}_{ij}$ and $\ryskus{Z}_{ij}$ through a
known monotonic and differentiable link function $\vartheta$, so that the
GAPLM is given as
%
\begin{equation}\label{model}
\eta_{ij}=\vartheta(\mu_{ij})=\ryskus{X}_{ij}^{\statT}
\bbeta%
+\sum_{l=1}^{d_{2}}
\theta_{l}(Z_{ijl}),\qquad j=1,\ldots,m_{i},
i=1,\ldots,n,
\end{equation}
where $\bbeta$ is a $d_{1}$-dimensional regression parameter, and $%
\theta_{l}$, $l=1,\ldots,d_{2}$, are unknown but smooth functions. We
assume $%
\underline{\bvarepsilon}_{i}= ( \varepsilon_{i1},\ldots,\varepsilon
_{im_{i}} )^{\statT}\sim(
\ryskus{0},\bSigma_{i} ) $. For identifiability, both the additive
and linear
components must be centered, that is, $E\theta_{l}(Z_{ijl})\equiv0$, $%
l=1,\ldots,d_{2}$, $EX_{ijk}=0$, $k=1,\ldots,d_{1}$. Model (\ref{model}) can
either become a generalized additive model~\cite{HT90} if the parameter
vector $\bbeta=\ryskus{0}$ or be a generalized linear model if $\theta
_{l}(\cdot)=0,1\leq l\leq d_{2}$. Model (\ref{model}) is more parsimonious
and easier to interpret than purely generalized additive models by allowing
a subset of predictors to be discrete and unbounded, modeled as some of the
variables $ ( X_{ijk} )_{k=0}^{d_{1}-1}$ and more flexible than
generalized linear models by allowing nonlinear relationships.

The GEE methods have been widely applied to analyze clustered data with
small cluster sizes and a large number of subjects $n$. However, data with
large cluster sizes have occurred frequently in various fields such as
machine learning, pattern recognition, image analysis, information
retrieval and bioinformatics. Reference~\cite{XY03} first studied the asymptotics for
parametric GEE estimators with large cluster sizes. As an extension, we
develop asymptotic properties of the spline GEE estimators in the GAPLMs
(\ref{model}) when the cluster sizes are allowed to increase with $n$,
that is,
the maximum cluster size $m_{ ( n ) }=\max_{1\leq i\leq
n}m_{i}$ is a function of $n$, such that $m_{ ( n ) }\rightarrow
\infty$ as $n\rightarrow\infty$.

The one-step spline estimation in~\cite{MSW12} for GAPLMs with identity link
is fast to compute but lacks limiting distribution. The traditional
backfitting approach has been widely used to estimate additive models for
independent and identically distributed (i.i.d.) and weekly-dependent
data~\cite{HT90,OR97,MLN99}. It, however, has computational
burden issues, due to its iterative nature. Moreover, it is pointed out in
\cite{HM05} that derivation of the asymptotic properties of a backfitting
estimator for a model with a link function is very complicated. As an
alternative,~\cite{L00,HM05,HKM06} and~\cite{HL05} proposed
two-stage kernel based estimators for i.i.d. data including one step
backfitting of the integration estimators in~\cite{L00} and one step
backfitting of the projection estimators in~\cite{HKM06}, one Newton step
from the nonlinear least squares estimators in~\cite{HM05}, and the
extension of the method in~\cite{HM05} to additive quantile regression
models. The two-stage estimator enjoys the oracle property which backfitting
estimators do not have, that is, it performs as well as the univariate
function estimator by assuming that other components are known.

In this paper, we propose a two-step spline GEE approach to approximate
$%
\theta_{l} ( \cdot) $ for $1\leq l\leq d_{2}$ in model (\ref{model})
with $m_{ ( n ) }$ going to infinity or bounded, and establish oracle
efficiency such that the two-step spline GEE estimator
of $%
\theta_{l} ( \cdot) $ achieves the same asymptotic distribution
of the oracle estimator obtained by assuming that $\bbeta$ and
other functions $\theta_{l^{\prime}} ( \cdot) $ for $1\leq
l^{\prime}\leq d_{2}$ and $l^{\prime}\neq l$ are known. In the first step,
the additive components $\theta_{l^{\prime}} ( \cdot) $ for $%
1\leq l^{\prime}\leq d_{2}$ and $l^{\prime}\neq l$ are pre-estimated
by their pilot estimators through an undersmoothed spline procedure. In
the second step, a more smoothed spline estimating procedure is applied
to the univariate data to estimate $\theta_{l} ( \cdot) $
with asymptotic distribution established. The proposed two-step
estimators achieve uniform oracle efficiency by ``reducing bias via
undersmoothing'' in the first step and ``averaging out the
variance'' in the second step. We establish asymptotic
consistency and normality of the one-step estimator for the parameter
vector and the two-step estimators of the nonparametric components. The
two-step spline GEE approach is inspired by the idea of
``spline-backfitted kernel/spline
smoothing'' of~\cite{WY07,SY10,LY10}
and~\cite{MY11} for additive models, additive coefficient models and
additive partially linear models with i.i.d or weekly-dependent data by
using least squares. The complex correlations within the clusters as
well as the non-Gaussian nature of discrete data make the estimation
and development of asymptotic properties in the framework studied in
this paper much more challenging.

\section{Two-step spline estimating equations}
\label{SECMethodology}

For simplicity, we denote vectors $\underline{\ryskus{Y}}_{i}= \{ (
Y_{i1},\ldots,Y_{im_{i}} )^{\statT} \}_{m_{i}\times1}$ and
$\underline{\bldeta}_{i}= \{ ( \eta_{i1},\ldots,\eta_{im_{i}} )^{\statT
} \}_{m_{i}\times1}$, $1\leq m_{i}\leq m_{ ( n ) }$, $1\leq i\leq
n$.\vadjust{\goodbreak}
Let $\varepsilon_{ij}=Y_{ij}-\mu_{ij}$,
and $\underline{\bvarepsilon}_{i}= ( \varepsilon_{i1},\ldots,\varepsilon_{im_{i}} )^{\statT}$. Similarly, let $%
\underline{\ryskus{X}}_{i}= \{ ( \ryskus{X}_{i1},\ldots,\ryskus{X}%
_{im_{i}} )^{\statT} \}_{m_{i}\times d_{1}}$ and
$\underline{\ryskus{Z}}_{i}= \{ ( \ryskus{Z}_{i1},\ldots,\ryskus{Z}%
_{im_{i}} )^{\statT} \}_{m_{i}\times d_{2}}$. Assume that $%
Z_{ijl}$ is distributed on a compact interval $ [ a_{l},b_{l} ],
1\leq l\leq d_{2}$, and, without loss of generality, we take all
intervals $%
[ a_{l},b_{l} ] = [ 0,1 ],1\leq l\leq d_{2}$. We further
let $\theta_{l}(\ryskus{Z}_{il})= \{ \{ \theta_{l}(Z_{i1l}),\ldots,\theta_{l}(Z_{im_{i}l}) \}^{\statT} \}_{m_{i}\times1}$, for
$l=1,\ldots,d_{2}$. The mean function in model
(\ref{model}) can be written in matrix notation as $\underline{\bldeta}%
_{i}=\underline{\ryskus{X}}_{i}\bbeta+\sum_{l=1}^{d_{2}}%
\theta_{l}(\ryskus{Z}_{il})$, which is the marginal model~\cite{LZ86}.
Let $\mu( \cdot) =\vartheta^{-1} ( \cdot) $ be the
inverse of the link function and $\mu( \underline{\bldeta}%
_{i} ) = [ \{ \mu( \eta_{i1} ),\ldots,\mu(
\eta_{im_{i}} ) \}^{\statT} ]_{m_{i}\times1}$.

As in~\cite{WCL05}, we allow $\underline{\ryskus{X}}_{i}$ and
$\underline{\ryskus{Z}}_{i}$ to be dependent. Let $\ryskus
{V}_{i}=\ryskus{V}%
_{i} ( \underline{\ryskus{X}}_{i},\underline{\ryskus{Z}}_{i} ) $ be
the assumed ``working'' covariance of $%
\underline{\ryskus{Y}}_{i}$, where $\ryskus{V}_{i}=\ryskus{A}_{i}^{1/2}%
\ryskus{R}_{i} ( \balpha) \ryskus{A}_{i}^{1/2}$, $\ryskus{%
A}_{i}=\ryskus{A}_{i} ( \underline{\ryskus{X}}_{i},\underline{\ryskus
{Z}}%
_{i} ) $ denotes an $m_{i}\times m_{i}$ diagonal matrix that contains
the marginal variances of $Y_{ij}$, and $\ryskus{R}_{i}$ is an invertible
working correlation matrix, which depends on a nuisance parameter
vector $%
\balpha$. Let $\bSigma_{i}=\bSigma_{i} (
\underline{\ryskus{X}}_{i},\underline{\ryskus{Z}}_{i} ) $ be the\vspace*{1pt} true
covariance of $\underline{\ryskus{Y}}_{i}$. If $\ryskus{R}_{i}$ is
equal to
the true correlation matrix $\overline{\ryskus{R}}_{i}$, then $\ryskus
{V}_{i}=%
\bSigma_{i}$.

Following~\cite{WY07}, we approximate the nonparametric functions
$\theta_{l}$'s by centered polynomial splines. Let $G_{n}$ be the space of
polynomial splines of degree $q\geq1$. We introduce a knot sequence
with $%
N_{n}$ interior knots
\[
t_{-q}=\cdots=t_{-1}=t_{0}=0<t_{1}<\cdots<t_{N}<1=t_{N+1}=\cdots=t_{N+q+1},
\]
where $N\equiv N_{n}$ increases when the number of subjects $n$ increases,
with order assumption given in condition (A4). Then $G_{n}$ consists of
functions $\varpi$ satisfying the following: (i) $\varpi$ is a
polynomial of degree $q$
on each of the subintervals $I_{s}= [ t_{s},t_{s+1} ) $, $%
s=0,\ldots,N_{n}-1$, $I_{Nn}= [ t_{N_{n}},1 ] $; (ii) for $q\geq1$, $%
\varpi$ is $q-1$ time continuously differentiable on $[0,1]$. Let $%
J_{n}=N_{n}+q+1$. Let $ \{ b_{s,l}\dvtx 1\leq l\leq d_{2},1\leq s\leq
J_{n}+1 \}^{\statT}$ be a basis system of the space $G_{n}$. We
adopt the centered B-spline space $G_{n}^{0}$ introduced in~\cite{XY06},
where $\ryskus{B} ( \ryskus{z} ) = \{ B_{s,l} (
z_{l} )\dvtx 1\leq l\leq d_{2},1\leq s\leq J_{n} \}^{\statT}$ is
a basis system of the space $G_{n}^{0}$ with $B_{s,l} ( z_{l} ) =%
\sqrt{N_{n}} [ b_{s+1,l} ( z_{l} ) - \{ E (
b_{s+1,l} ) /E ( b_{1,l} ) \} b_{1,l} ( z_{l} ) %
] $ and $\ryskus{z}= ( z_{l} )_{l=1}^{d_{2}}$.\vspace*{1pt}

Equally-spaced knots are used in this article for simplicity of proof. Other
regular knot sequences can also be used, with similar asymptotic
results.\vspace*{9pt}

\textit{Step} I. Pilot estimators of $\bbeta$ and $\theta_{l}(\cdot
)$. Suppose that $\theta_{l}$ can be approximated well by a
spline function in $G_{n}^{0}$, so that
%
\begin{equation}\label{EQspline}
\theta_{l} ( z_{l} ) \approx\widetilde{\theta
}_{l} ( z_{l} ) =\sum_{s=1}^{J_{n}}
\gamma_{sl}B_{s,l} ( z_{l} ).
\end{equation}
Let $\bgamma= ( \gamma_{sl}\dvtx s=1,\ldots,J_{n},l=1,\ldots,d_{2} )^{\statT}$ be
the collection of the coefficients in (\ref{EQspline}),
and denote $\ryskus{B}_{ijl}= [ \{ B_{s,l} (
Z_{ijl} )\dvtx
s=1,\ldots,J_{n} \}^{\statT} ]_{J_{n}\times1}$ and $\ryskus{B}%
_{ij}= \{ ( \ryskus{B}_{ij1}^{\statT},\ldots,\ryskus{B}_{ijd_{2}}^{%
\statT} )^{\statT} \}_{d_{2}J_{n}\times1}$, then we
have an approximation $\eta_{ij}\approx\widetilde{\eta}_{ij}=\ryskus
{X}%
_{ij}^{\statT}\bbeta+\ryskus{B}_{ij}^{\statT}\bgamma$. We
can also write the approximation in matrix notation as $\underline
{\bldeta}_{i}
\approx\underline{\widetilde{\bldeta}}_{i}=\underline{\ryskus
{X}}_{i}\bbeta+
\underline{\ryskus{B}}_{i}\bgamma$, where $%
\underline{\ryskus{B}}_{i}= \{ ( \ryskus{B}_{i1},\ldots,\ryskus{B}%
_{im_{i}} )^{\statT} \}_{m_{i}\times d_{2}J_{n}}$. Let $\mu
( \widetilde{\underline{\bldeta}}_{i} ) = [ \{ \mu
( \widetilde{\eta}_{i1} ),\ldots,\break\mu( \widetilde{\eta}%
_{im_{i}} ) \}^{\statT} ]_{m_{i}\times1}$. Let $%
\widehat{\bbeta}_{n}= ( \widehat{\beta}_{n,1},\ldots,\widehat{%
\beta}_{n,d_{1}} )^{\statT}$ and $\widehat{\bgamma}%
_{n}= \{ \widehat{\gamma}_{n,sl}\dvtx
s=1,\ldots,J_{n},l=1,\ldots,d_{2} \}^{\statT}$ be the minimizer
of
%
\begin{equation}\label{Qnbetagamma}
Q_{n} ( \bbeta,\bgamma) =\frac{1}{2}\sum
_{i=1}^{n} \bigl\{ \underline{\ryskus{Y}}_{i}-
\mu( \underline{\ryskus{X}}_{i}\bbeta+\underline{
\ryskus{B}}_{i}\bgamma) \bigr\}^{\statT}\ryskus{V}_{i}^{-1}
( \bbeta,\bgamma) \bigl\{ \underline{\ryskus{Y}}_{i}-\mu(
\underline{\ryskus{X}}_{i}\bbeta+%
\underline{
\ryskus{B}}_{i}\bgamma) \bigr\},\hspace*{-30pt}
\end{equation}
which is corresponding to the class of working covariance matrices $ \{
\ryskus{V}_{i},1\leq i
\leq n \} $. Then $\widehat{\bbeta}_{n}$ and $\widehat{%
\bgamma}_{n}$ solve the estimating equations
%
\begin{equation}\label{EQLS}
\ryskus{g}_{n} ( \bbeta,\bgamma) =\sum_{i=1}^{n}
\underline{\ryskus{D}}_{i}^{\statT}\Delta_{i} (
\bbeta,\bgamma%
) \ryskus{V}_{i}^{-1} ( \bbeta,
\bgamma) \bigl\{ \underline{\ryskus{Y}}_{i}-\mu( \underline{
\ryskus{X}}_{i} \bbeta+\underline{\ryskus{B}}_{i}\bgamma
) \bigr\} =\ryskus{0},
\end{equation}
where $\underline{\ryskus{D}}_{i}= ( \underline{\ryskus{X}}_{i},
\underline{\ryskus{B}}_{i} )_{m_{i}\times( d_{1}+d_{2}J_{n} ) }$,
and
\[
\Delta_{i} ( \bbeta,\bgamma) =\operatorname{diag} \bigl(
\Delta_{i1} ( \bbeta,\bgamma),\ldots,\Delta_{im_{i}} ( \bbeta,
\bgamma) \bigr)
\]
is a diagonal matrix with the diagonal elements being the first derivative
of $\mu( \cdot) $ evaluated at $\widetilde{\eta}_{ij}$, $%
j=1,\ldots,m_{i}$. Then we let $\widehat{\bbeta}_{n}$ be the
estimator of the parameter vector $\bbeta$. For each $1\leq l\leq
d_{2}$, the pilot estimator of the $l$th nonparametric function $\theta
_{l}(z_{l})$ is $\widehat{\theta}_{n,l}(z_{l})=\sum
_{s=1}^{J_{n}}\widehat{%
\gamma}_{n,sl}B_{s,l} ( z_{l} ) $. The one-step spline estimator
of each function component has consistency properties, but lacks limiting
distribution \mbox{\cite{WY07,MY11,MSW12}}.\vspace*{9pt}

\textit{Step} II. Two-step spline GEE estimator of $\theta_{l}(\cdot)$.
Next, we propose a two-step spline estimator of $\theta_{l}(\cdot)$ for
given $1\leq l\leq d_{2}$. The basic idea is that for every $1\leq l\leq
d_{2}$, we estimate the $l$th function $\theta_{l}(\cdot)$ in model
(\ref{model}) nonparametrically with the GEE method by assuming that the parameter
vector $\bbeta$ and other nonparametric components $\btheta%
_{-l}= \{ \theta_{l^{\prime}}(\cdot)\dvtx 1\leq l^{\prime}\leq
d_{2},l^{\prime}\neq l \} $ are known. The problem turns into a
univariate function estimation problem. Because the true parameter
vector $%
\bbeta$ and functions $\btheta_{-l}$ are not known in
reality, we replace them by their pilot estimators from step I to
obtain the
two-step estimator of $\theta_{l}(\cdot)$. Both kernel and spline based
methods can be employed in the second step to estimate $\theta
_{l}(\cdot)$.
Here we choose the spline method described in the beginning of this
section. We use the splines of the same degree $q$ as in step I. Denote
$%
\ryskus{B}_{ijl}^{\calS}= [ \{ B_{s,l}^{\calS} (
Z_{ijl} )\dvtx s=1,\ldots,J_{n}^{\calS} \}^{\statT} ]_{J_{n}^{\calS}\times1}$,
where $B_{s,l}^{\calS} (
z_{l} ) $ is the spline function defined in the same way as $%
B_{s,l} ( z_{l} ) $ in step I, but with $N^{\calS}\equiv
N_{n}^{\calS}$ the number of interior knots and let $J_{n}^{\calS%
}=N^{\calS}+q+1$. Denote $\ryskus{B}_{l}^{\calS}
(z_{l} ) =\{ B_{s,l}^{\calS} ( z_{l} )$,
$s=1,\ldots,J_{n}^{\calS}\}^{\statT}$, $\ryskus{B}_{i\cdot l}^{%
\calS}= \{ ( \ryskus{B}_{i1l},\ldots,\ryskus{B}%
_{im_{i}l} )^{\statT} \}_{m_{i}\times J_{n}^{\calS}}$,
and $\bgamma_{l}^{\calS}= ( \gamma_{sl}\dvtx s=1,\ldots,J_{n}^{%
\calS} )^{\statT}$. By assuming that $\bbeta$ and $%
\btheta_{-l}= \{ \theta_{l^{\prime}}(\cdot)\dvtx l^{\prime}\neq
l$, $1\leq l^{\prime}\leq d_{2} \} $ are known, $\theta_{l}(z_{l})$
is estimated by the oracle estimator
%
\begin{equation}\label{EQORACLE}
\widehat{\theta}_{n,l}^{\calS}(z_{l},\bbeta,
\btheta%
_{-l})=\sum_{s=1}^{J_{n}}
\widehat{\gamma}_{n,sl}^{\calS} ( \bbeta,\btheta_{-l}
) B_{s,l}^{\calS} ( z_{l} ) =\ryskus{B}_{l}^{\calS}
( z_{l} )^{\statT}\widehat{%
\bgamma
}_{n,l}^{\calS} ( \bbeta,\btheta%
_{-l} )
\end{equation}
with $\widehat{\bgamma}_{n,l}^{\calS} ( \bbeta,\btheta_{-l} ) = \{
\widehat{\gamma}_{n,sl}^{\calS} (
\bbeta,\btheta_{-l} ) \}_{s=1}^{J_{n}^{\calS}}$
solving the estimating equation%
%
\begin{eqnarray}
\label{EQgnlstar}
&&\ryskus{g}_{n,l}^{\calS} \bigl(
\bgamma_{l}^{\calS},%
\bbeta,\btheta_{-l}
\bigr) \nonumber\\[-2pt]
&&\qquad=\sum_{i=1}^{n} \bigl(
\ryskus{B}%
_{i\cdot l}^{\calS} \bigr)^{\statT}
\Delta_{i} \bigl( \bbeta,\btheta_{-l},\bgamma_{l}^{\calS}
\bigr) \ryskus{V}%
_{i}^{-1} \bigl( \bbeta,
\btheta_{-l},\bgamma_{l}^{\calS} \bigr)
\nonumber\\[-9pt]\\[-9pt]
&&\hspace*{12pt}\qquad\quad{}\times \Biggl\{ \underline{\ryskus{Y}}_{i}-\mu\Biggl( \underline{
\ryskus{X}}_{i}%
\bbeta+\sum_{l^{\prime}=1,l^{\prime}\neq
l}^{d_{2}}
\theta_{l^{\prime}}(\ryskus{Z}_{il^{\prime}})+ \bigl( \ryskus{B}%
_{i\cdot l}^{\calS}
\bigr)^{\statT}\bgamma_{l}^{\calS} \Biggr) \Biggr\} \nonumber\\[-2pt]
&&\qquad=
\ryskus{0},
\nonumber
\end{eqnarray}
where $\Delta_{i} ( \bbeta,\btheta_{-l},\bgamma_{l}^{%
\calS} )=\operatorname{diag} ( \Delta_{i1} ( \eta_{i1}^{%
\calS} ),\ldots,\Delta_{im_{i}} ( \eta_{im_{i}}^{\calS} ) ) $, and
$\Delta_{ij} ( \eta_{ij}^{\calS%
} ) $ is the first derivative of $\mu( \cdot) $ evaluated
at $\eta_{ij}^{\calS}=\ryskus{X}_{ij}^{\statT}\bbeta%
+\sum_{l^{\prime}=1,l^{\prime}\neq l}^{d_{2}}\theta_{l^{\prime
}}(Z_{ijl^{\prime}})+ ( \ryskus{B}_{ijl}^{\calS} )^{\statT}\bgamma
_{l}^{\calS}$, $j=1,\ldots,m_{i}$.
We replace the
true parameter vector $\bbeta$ and the true functions $\btheta_{-l}=
\{ \theta_{l^{\prime}}(\cdot),1\leq l^{\prime}\leq
d_{2},l^{\prime}\neq l \} $ with the pilot estimators $\widehat{%
\bbeta}_{n}$ and $\widehat{\btheta}_{n,-l}= \{
\widehat{\theta}_{n,l^{\prime}}(\cdot),1\leq l^{\prime}\leq
d_{2},l^{\prime}\neq l \} $, where $\widehat{\theta}_{n,l^{\prime
}}(z_{l^{\prime}})=\sum_{s=1}^{J_{n}}\widehat{\gamma}_{n,sl^{\prime
}}B_{s,l^{\prime}} ( z_{l^{\prime}} ) $, so that $\theta_{l}(z_{l})$
is estimated by the two-step spline estimator
%
\begin{equation}\label{EQTSEl}
\widehat{\theta}{}^{\calS}_{n,l}(z_{l},\widehat{ \bbeta}_{n},
\widehat{\btheta}_{n,-l})=\ryskus{B}_{l}^{\calS} ( z_{l}
)^{\statT}\widehat{\bgamma}_{n,l}^{\calS}
( \widehat{\bbeta}_{n},\widehat{\btheta}%
_{n,-l}).
\end{equation}
The Newton--Raphson algorithm of GEE is applied to obtain
$\widehat{\bgamma}_{n,l}^{\calS}$. Define
\begin{eqnarray*}
\calD_{n} ( \bbeta,\bgamma) & = & \bigl\{ -\partial
\ryskus{g}_{n} ( \bbeta,\bgamma) /\partial\bigl(
\bbeta^{\statT},\bgamma^{\statT} \bigr) \bigr\}_{ (
d_{1}+d_{2}J_{n} ) \times( d_{1}+d_{2}J_{n} ) },
\\
\Psi_{n} ( \bbeta,\bgamma) & = & \Biggl\{ \sum
_{i=1}^{n}\underline{\ryskus{D}}_{i}^{\statT}
\Delta_{i} ( \bbeta,\bgamma) \ryskus{V}_{i}^{-1}
( \bbeta,\bgamma) \Delta_{i} ( \bbeta,\bgamma) \underline{
\ryskus{D}}_{i} \Biggr\}_{ ( d_{1}+d_{2}J_{n} )
\times( d_{1}+d_{2}J_{n} ) }.
\end{eqnarray*}

\section{Asymptotic properties of the estimators}\label{SECAsymptotics}

For any $s\times s$ symmetric matrix~$\ryskus{A}$, denote by $\lambda
_{\min
} ( \ryskus{A} ) $ and $\lambda_{\max} ( \ryskus{A} ) $
its smallest and largest eigenvalues. For any vector $\balpha=
( \alpha_{1},\ldots,\alpha_{s} )^{\statT}$, let its
Euclidean norm be $\llVert\balpha\rrVert=\sqrt{\alpha
_{1}^{2}+\cdots+\alpha_{s}^{2}}$. Let $C^{0,1} ( \cal{X} %
_{w} ) $ be the space of Lipschitz continuous functions on $\cal{X} %
_{w}$, that is,%
\[
C^{0,1} ( \cal{X}_{w} ) = \biggl\{ \varphi\dvtx \llVert\varphi
\rrVert_{0,1}=\sup_{w\neq w^{\prime},w,w^{\prime}\in\cal{X} _{w}}%
\frac{\llvert\varphi(w)-\varphi( w^{\prime} ) \rrvert}{%
\llvert w-w^{\prime}\rrvert}<+\infty
\biggr\},
\]
in which $\llVert\varphi\rrVert_{0,1}$ is the $C^{0,1}$-norm of $%
\varphi$. Throughout the paper, we assume the following regularity
conditions:

\begin{longlist}[(C5)]
\item[(C1)] The random variables $Z_{ijl}$ are bounded, uniformly in
$1\leq
j\leq m_{i}$, $1\leq i\leq n$, $1\leq l\leq d_{2}$. The marginal
density $%
f_{ijl} ( \cdot) $ of $Z_{ijl}$ is bounded away from $0$ and $%
\infty$ on $ [ 0,1 ] $, uniformly in $1\leq j\leq m_{i}$, $1\leq i\leq
n$. The joint density $f_{ijlj^{\prime}l^{\prime}} ( \cdot ,\cdot) $ of
$ ( Z_{ijl},Z_{ij^{\prime}l^{\prime}} ) $ is\vspace*{1pt}
bounded away from $0$ and $\infty$ on $ [ 0,1 ] $, uniformly in $%
1\leq i\leq n$, $1\leq j,j^{\prime}\leq m_{i}$, and $1\leq l\neq
l^{\prime
}\leq d_{2}$.

\item[(C2)] The eigenvalues of the true correlation matrices $\overline{
\ryskus{R}}_{i}$ are bounded away from $0$, uniformly in $1\leq i\leq n$.\vadjust{\goodbreak}

\item[(C3)] The eigenvalues of the inverse of the working correlation
matrices $\ryskus{R}_{i} ( \alpha)^{-1}$ are bounded away from $%
0 $, uniformly in $1\leq i\leq n$.

\item[(C4)] Let $n_{\statT}=\sum_{i=1}^{n}m_{i}$. There
are constants $0<c<C<\infty$, such that
$cn_{\statT}\leq\lambda_{\mathrm{min}} (
\sum_{i=1}^{n}\underline{\ryskus{X}}_{i}^{\statT}
\underline{\ryskus{X}}_{i} ) \leq\lambda_{\mathrm{max}} (
\sum_{i=1}^{n}\underline{\ryskus{X}}_{i}^{\statT}
\underline{\ryskus{X}}_{i} ) \leq Cn_{\statT}$.

\item[(C5)] For $1\leq l\leq d_{2}$, $\theta_{l}^{{ ( p-1 )
}}(z_{l})\in C^{0,1} [ 0,1 ] $, for given integer $p\geq1$. The
spline degree satisfies $q+1\geq p$, and $\mu^{\prime} ( \eta)
\in C^{0,1} ( \cal{X}_{\eta} ) $. The number of interior knots
$N_{n}\rightarrow\infty$, as $n_{\statT}\rightarrow\infty$.
\end{longlist}

Conditions (C1)--(C4) are similar to conditions (A1)--(A4) in~\cite{MSW12},
and condition (C5) is weaker than the first part of condition (A5) in
\cite%
{MSW12}. Let $\bbeta_{0}$ be the true parameter vector and $\theta_{l0}
( \cdot) $ be the true $l$th additive function in model
(\ref{model}). According to the result on page 149 of~\cite{dB01}, for $
\theta_{l0} ( \cdot) $ satisfying condition (C5), there is a
function
%
\begin{equation}\label{DEFglzl}
\widetilde{\theta}_{l0} ( z_{l} ) =\sum
_{s=1}^{J_{n}}\gamma_{sl,0}B_{s,l} (
z_{l} ) \in G_{n}^{0},
\end{equation}
such that $\sup_{z_{l}\in[ 0,1 ] }\llvert\widetilde{%
\theta}_{l0} ( z_{l} ) -\theta_{l0} ( z_{l} )
\rrvert=O ( J_{n}^{-p} ) $. Thus, by letting $\bgamma%
_{0}= ( \gamma_{sl,0}\dvtx\break s=1,\ldots,J_{n},l=1,\ldots,d_{2} )^{\statT}$,
\[
\sup_{\ryskus{z}\in[ 0,1 ] ^{d_{2}}}\Biggl\llvert\ryskus{%
B} ( \ryskus{z} )^{\statT}\bgamma_{0}-\sum%
_{l=1}^{d_{2}}\theta_{l0} ( z_{l} ) \Biggr\rrvert\leq
\sum_{l=1}^{d_{2}}\sup_{z_{l}\in[ 0,1 ]
}\bigl\llvert\widetilde{\theta}_{l0} ( z_{l} ) -\theta_{l0} (
z_{l} ) \bigr\rrvert=O \bigl( d_{2}J_{n}^{-p} \bigr).
\]
In addition to the regularity conditions above, we need extra
conditions to ensure the existence and weak consistency of the
estimators in (\ref{EQLS}). Let $ 
\lambda_{n}^{\min}=\min_{1\leq i\leq n}\lambda_{\mathrm{min}}%
\{ \ryskus{R}_{i}^{-1} ( \alpha) \} $, $\lambda_{n}^{\max}=\max
_{1\leq i\leq n}\lambda_{\max} \{ \ryskus{R}%
_{i}^{-1} ( \alpha) \} $, $\tau_{n}^{\max}=\break\max_{1\leq
i\leq n} \{ \lambda_{\max} ( \ryskus{R}_{i}^{-1} ( \alpha
) \overline{\ryskus{R}}_{i} ) \} $ and $\tau_{n}^{\min
}=\min_{1\leq i\leq n} \{ \lambda_{\min} ( \ryskus{R}%
_{i}^{-1} ( \alpha) \overline{\ryskus{R}}_{i} ) \} $.
The additional conditions are as follows:

\begin{longlist}[(A2)]
\item[(A1)] $ ( \lambda_{n}^{\min}/\tau_{n}^{\max} )
n_{\statT}/J_{n}^{1/2}\rightarrow\infty$.

\item[(A2)] There is a constant $c_{0}>0$, for any $r>0$, such that $%
P \{ \calD_{n} ( \bbeta,\bgamma) \geq
c_{0}\Psi_{n} ( \bbeta_{0},\bgamma_{0} )$
and $\calD_{n} ( \bbeta,\bgamma)$ is\vspace*{2pt}
nonsingular, for all $( \bbeta^{\statT},\bgamma^{%
\statT} )^{\statT}\in\xi_{n} ( r ) \}
\rightarrow1$, where $\xi_{n} ( r ) = \{ (
\bbeta^{\statT},\bgamma^{\statT} )^{\statT}\dvtx \llVert\{ \Psi_{n} (
\bbeta_{0},\bgamma_{0} ) \}^{1/2} ( ( \bbeta-\bbeta_{0} )^{
\statT}$,
$( \bgamma-\bgamma_{0} )^{\statT} )^{\statT}\rrVert\leq\break( \tau
_{n}^{\max} )^{1/2}r \} $.
\end{longlist}

Conditions (A1) and (A2) are used to ensure the existence and weak
consistency of the solutions in (\ref{EQLS}). Condition (A2)
corresponds to
condition (L$_{\mathrm{w}}^{*}$) in~\cite{XY03} for generalized linear
models. Conditions (A1) and (C4) imply condition (I$_{\mathrm{w}}^{*}$)
in~\cite{XY03}, which will be proved in the
\hyperref[app]{Appendix}.\vspace*{1pt}
Condition (A2) relates
to the true and the working correlation structures $\overline{\ryskus
{R}}%
_{i} $ and $\ryskus{R}_{i} ( \alpha) $. Since the true
correlations $\overline{\ryskus{R}}_{i}$ are often not completely specified
and $\max_{1\leq i\leq n}\lambda_{\max} ( \overline{\ryskus{R}}%
_{i} ) \leq m_{ ( n ) }$, then condition (A1) is implied by

\begin{longlist}[(A1$^*$)]
\item[(A1$^*$)] $ ( \lambda_{n}^{\min}/\lambda_{n}^{\max} )
m_{ ( n ) }^{-1}n_{\statT}/J_{n}^{1/2}\rightarrow\infty$.\vadjust{\goodbreak}
\end{longlist}

Condition (A1$^{*}$) does not contain $\overline{\ryskus{R}}_{i}$.
Thus, the order requirements of $n$, $m_{ ( n ) }$ and $J_{n}$
depend on the choice of the working correlations $\ryskus{R}_{i} (
\alpha) $. For instance, if the working correlation structures are
independent or AR(1) within each subject, then there exist constants $%
0<c_{R}\leq C_{R}<\infty$, such that \mbox{$c_{R}\leq( \lambda_{n}^{\max
} )^{-1}\lambda_{n}^{\min}\leq C_{R}$}. Thus, condition (A1$^{*%
} $) is equivalent to $m_{ ( n ) }^{-1}n_{\statT
}/\break J_{n}^{1/2}\rightarrow\infty$. For exchangeable working correlation
structures, there exist constants $0<C_{R}^{\prime}<\infty$, such
that $%
\max_{1\leq i\leq n}\lambda_{\mathrm{max}} \{ \ryskus{R}%
_{i} ( \alpha) \} \leq C_{R}^{\prime}m_{ ( n ) }$,
then\break $ ( \lambda_{n}^{\max} )^{-1}\lambda_{n}^{\min}\geq
c_{R}^{{\prime}}m_{ ( n ) }^{-1}$, for some constant $%
0<c_{R}^{\prime}<\infty$. Condition (A1$^{*}$) is implied by $%
m_{ ( n ) }^{-2}n_{\statT}/J_{n}^{1/2}\rightarrow\infty
$.
%
\begin{theorem}
\label{THMweakconsistency}Under conditions \textup{(A1)} and \textup{(A2)} or
\textup{(A1$^{*}$)} and \textup{(A2)}, as $n_{\statT}\rightarrow\infty$, there exist
sequences of random variables $\widehat{\bbeta}_{n}$ and $\widehat{%
\bgamma}_{n}$, such that $P \{ \ryskus{g}_{n} ( \widehat{%
\bbeta}_{n},\widehat{\bgamma}_{n} ) =0 \}
\rightarrow1$, and $\llVert\widehat{\bbeta}_{n}-\bbeta%
_{0}\rrVert\rightarrow0$ and $\llVert\widehat{\gamma}_{n}-%
\bgamma_{0}\rrVert\rightarrow0$ in probability.
\end{theorem}

Next we derive the asymptotic properties of $\widehat{\bbeta}_{n}$.
Let $\cal{X} $ and $\cal{Z} $ be the collections of all $X_{ijk}$'s
and $Z_{ijl}$'s, respectively, that is, $\cal{X}_{n_{\statT}\times
d_{1}}= ( \underline{\ryskus{X}}_{1}^{\statT},\ldots,
\underline{\ryskus{X}}_{n}^{\statT} )^{\statT}$ and $\cal{Z}_{n_{%
\statT}\times d_{2}}= ( \underline{\ryskus{Z}}_{1}^{\statT},\ldots,
\underline{\ryskus{Z}}_{n}^{\statT} )^{\statT}$. Let $\Delta_{i}$ be
the diagonal matrix with the diagonal elements being the first
derivative\vspace*{-1pt} of $\mu( \cdot) $ evaluated at $\ryskus{X}_{ij}^{%
\statT}\bbeta_{0}+\sum_{l=1}^{d_{2}}\theta_{l0}(Z_{ijl})$, $j=1,\ldots,m_{i}$, and $\ryskus{V}_{i}=\ryskus{A}%
_{i}^{1/2}\ryskus{R}_{i} ( \balpha) \ryskus{A}_{i}^{1/2}$
with $\ryskus{A}_{i}$ being the marginal variance of $\underline{\ryskus
{Y}}%
_{i}$ evaluated at the true parameters and additive functions. To make $
\bbeta$ estimable, we need a condition to ensure $\cal{X} $ and $%
\cal{Z} $ not functionally related, which is similar to the condition
given in~\cite{MSW12}. Define the Hilbert space $\cal{H} = \{ \psi
( \ryskus{z} ) =\sum_{l=1}^{d_{2}}\psi_{l} ( z_{l} )
,E\psi_{l} ( z_{l} ) =0,\llVert\psi_{l}\rrVert_{2}<\infty\} $ of
theoretically centered $L_{2}$
additive functions on $ [ 0,1 ]^{d_{2}}$, where $\llVert\psi
_{l}\rrVert_{2}= \{ \int_{0}^{1}\psi_{l}^{2} (
z_{l} ) \,dz_{l} \}^{1/2}$. Let $\psi_{k}^{\ast}$ be the function
$\psi\in\cal{H} $ that minimizes\break $\sum_{i=1}^{n}E [
\{ \underline{\ryskus{X}}_{i}^{ ( k ) }-\psi(
\underline{\ryskus{Z}}_{i} ) \}^{\statT}\Delta_{i}\ryskus{V}%
_{i}^{-1}
\Delta_{i} \{ \underline{\ryskus{X}}_{i}^{ ( k )
}-\psi( \underline{\ryskus{Z}}_{i} ) \} ] $, where $%
\underline{\ryskus{X}}_{i}^{ ( k ) }= ( X_{i1k},\ldots,\break
X_{im_{i}k} )^{\statT},1\leq k\leq d_{1}$. Some other assumptions
needed are given as follows.

\begin{longlist}[(A3)]
\item[(A3)] Given $1\leq k\leq d_{1}$, $\psi_{k}^{\ast^{ ( p-1 )
}}\in C^{0,1} [ 0,1 ] $, for $1\leq p\leq q+1$.
\end{longlist}

The order requirements of the number of interior knots $N$ and
$N^{\calS}$ in steps I and~II are given in the following assumption:

\begin{longlist}[(A4)]
\item[(A4)] (i) $\sqrt{ ( \log n_{\statT} ) N^{\calS}/n_{%
\statT}} ( \tau_{n}^{\max}/\lambda_{n}^{\min} )^{1/2}=o ( 1 ) $,
$ ( N^{\calS} )^{-p-1/2}n_{\statT}^{1/2} ( \lambda_{n}^{\max}/\break\lambda
_{n}^{\min} ) ( \lambda_{n}^{\max}/\tau_{n}^{\min} )^{1/2}=O ( 1 ) $,
and (ii)\vspace*{2pt} $ ( \lambda_{n}^{\max}/\tau_{n}^{\min} )^{1/2} (
\lambda_{n}^{\max}/\lambda_{n}^{\min} )^{2} ( n_{\statT%
}/\break N^{\calS} )^{1/2}\*N^{-p}=o ( 1 ) $,
$ ( \lambda_{n}^{\max}/\tau_{n}^{\min} )^{1/2} ( \lambda_{n}^{\max
}/\lambda_{n}^{\min} )^{2} ( \log n_{\statT}/N^{\calS} )^{1/2}=o ( 1)
$,\break
$ ( N_{n}^{\calS}N_{n}\times\log n_{\statT} )^{1/2}n_{%
\statT}^{-1}=o ( 1 ) $.
\end{longlist}

Since $\lambda_{n}^{\min}\leq\tau_{n}^{\min}\leq\tau_{n}^{\max}\leq
m_{ ( n ) }\lambda_{n}^{\max}$, condition (A4) is implied by a
stronger condition as below:\vfill\eject

\begin{longlist}[(A4$^{*}$)]
\item[(A4$^{*}$)] (i) $\sqrt{ ( \log n_{\statT} ) N^{%
\calS}/n_{\statT}}m_{ ( n ) }^{1/2} ( \lambda_{n}^{\max}/\lambda
_{n}^{\min} )^{1/2}=o ( 1 ) $,
$ ( N^{\calS} )^{-p-1/2}\times\break n_{\statT}^{1/2} ( \lambda_{n}^{\max}/\lambda
_{n}^{\min} )^{3/2}=O ( 1 ) $,
and (ii) $ ( \lambda_{n}^{\max}/\lambda_{n}^{\min} )^{5/2} ( n_{\statT
}/N^{\calS} )^{1/2}N^{-p}=o (
1 ) $, $ ( \lambda_{n}^{\max}/\lambda_{n}^{\min} )^{5/2} ( \log
n_{\statT}/N^{\calS} )^{1/2}=o ( 1 ) $, $%
( N_{n}^{\calS}N_{n}\log n_{\statT} )^{1/2}n_{%
\statT}^{-1}=o ( 1 ) $.
\end{longlist}

Condition (A3) is weaker than the second part of condition (A5) in \cite
{MSW12}. Condition (A4$^{*}$) does not depend on the true correlation
matrices $\overline{\ryskus{R}}_{i}$, which are not specified. It is clear
that the first conditions in (A4) and (A4$^{*}$) ensure conditions
(A1) and (A1$^{*}$), respectively.
%
\begin{remark}\label{remark1} (A4)(i) lists the order requirements for
$N^{\calS}$ to obtain the asymptotic results of the oracle
estimator in Theorem~\ref{THMthetahatnormality}. (A4)(ii) ensures the
uniform oracle efficiency of the two-step spline estimator. It will be
shown in Theorem~\ref{THMthetahatconvergence} that the difference
between the two-step spline and the oracle estimators is of uniform
order $O_{P} \{ ( \lambda_{n}^{\max}/\lambda_{n}^{\min} )^{2} (
J_{n}^{-p}+\sqrt{\log n_{%
\statT}/n_{\statT}} ) \} $ with $%
O_{P} \{ ( \lambda_{n}^{\max}/\lambda_{n}^{\min} )^{2}%
\sqrt{\log n_{\statT}/n_{\statT}} \} $ and $%
O_{P} \{ ( \lambda_{n}^{\max}/\lambda_{n}^{\min} )^{2}J_{n}^{-p} \} $
caused by the noise and bias terms, respectively, in
the first step spline estimation. The inverse of the asymptotic standard
deviation of the oracle estimator is of order $O \{ \sqrt{n_{\statT%
}/J_{n}^{\calS}} ( \lambda_{n}^{\max}/\tau_{n}^{\min} )^{1/2} \} $.
The first two conditions of (A4)(ii) ensure that the
difference is asymptotically uniformly negligible. If we let $N$ have the
order $n_{\statT}^{1/ ( 2p ) }$, then the difference is
of uniform order $O_{P} \{ ( \lambda_{n}^{\max}/\break\lambda_{n}^{\min}
)^{2}\sqrt{\log n_{\statT}/n_{\statT}%
} \} $. Therefore, an undersmoothing procedure is applied in the first
step to reduce the bias. When $\lambda_{n}^{\min}$, $\lambda_{n}^{\max}$,
$\tau_{n}^{\min}$ and $\tau_{n}^{\max}$ are finite numbers, (A4)(i)
becomes $\sqrt{ ( \log n_{\statT} ) N^{\calS}/n_{%
\statT}}=o ( 1 ) $ and $ ( N^{\calS} )^{-p-1/2}n_{\statT}^{1/2}=O ( 1
) $. The optimal order of $N^{%
\calS}$ is $n_{\statT}^{1/ ( 2p+1 ) }$. Define
\[
\widetilde{\ryskus{X}}_{ik}=\underline{
\ryskus{X}}_{i}^{ (
k ) }-\psi_{k}^{\ast} (
\underline{\ryskus{Z}}_{i} ),\qquad 1\leq k\leq d_{1},\qquad \widetilde{
\underline{\ryskus{X}}}_{i}= ( \widetilde{%
\ryskus{X}}_{i1},\ldots,\widetilde{\ryskus{X}}_{id_{1}}
)_{m_{i}\times d_{1}}.
\]
Define $\widetilde{\Psi}_{n}=\sum_{i=1}^{n}\widetilde{\underline
{\ryskus{X}}%
}_{i}^{\statT}\Delta_{i}\ryskus{V}_{i}^{-1}\Delta_{i}\widetilde{%
\underline{\ryskus{X}}}_{i}$, $\widetilde{\Phi}_{n}=\sum_{i=1}^{n}%
\widetilde{\underline{\ryskus{X}}}_{i}^{\statT}\Delta_{i}\ryskus{V}%
_{i}^{-1}\Sigma_{i}\ryskus{V}_{i}^{-1}\Delta_{i}\widetilde{\underline{%
\ryskus{X}}}$, and
%
\begin{equation}\label{DEFgammatilda}
\widetilde{\Xi}_{n}= \bigl\{ E ( \widetilde{\Psi}_{n} )
\bigr\}^{-1}E ( \widetilde{\Phi}_{n} ) \bigl\{ E (
\widetilde{\Psi}%
_{n} ) \bigr\}^{-1}.
\end{equation}
\end{remark}

The following result gives the asymptotic distribution and consistency rate
of $\widehat{\bbeta}_{n}$ for general working covariance matrices.
%
\begin{theorem}
\label{THMbetahatnormality}Under conditions \textup{(A2)--(A4)}, as $n_{\statT%
}\rightarrow\infty$, $\widetilde{\Xi}_{n}^{-1/2} ( \widehat{
\bbeta}_{n}-\bbeta_{0} )
\rightarrow\operatorname{Normal} ( 0,\ryskus{I}_{d_{1}} ) $, and $%
\llVert\widehat{\bbeta}_{n}-\bbeta_{0}\rrVert
=O_{p} \{ n_{\statT}^{-1/2} ( \tau_{n}^{\max} )^{1/2} ( \lambda
_{n}^{\min} )^{-1/2} \} $. If condition
\textup{(A4)} is replaced by \textup{(A4$^{*}$)}, then
\[
\llVert\widehat{\bbeta}_{n}-\bbeta_{0}\rrVert
=O_{p} \bigl\{ n_{\statT}^{-1/2}m_{ ( n ) }^{1/2}
\bigl( \lambda_{n}^{\max} \bigr)^{1/2} \bigl(
\lambda_{n}^{\min} \bigr)^{-1/2} \bigr\}.
\]
\end{theorem}
%
\begin{remark}\label{remark2}
It is easy to show that the covariance $\widetilde{\Xi}%
_{n}$ in (\ref{DEFgammatilda}) is minimized when the working covariance
matrices are equal to the true covariance matrices such that $\ryskus{V}
_{i}=\Sigma_{i}$ for all $1\leq i\leq n$, and in this case equal to $%
\{ E ( \widetilde{\Psi}_{n} ) \}^{-1}$. To construct\vspace*{1pt}
the confidence sets for $\bbeta$, $\widetilde{\Xi}_{n}$ is
consistently estimated by $\widehat{\Xi}_{n}=\widehat{\Psi}_{n}^{-1}%
\widehat{\Phi}_{n}\widehat{\Psi}_{n}^{-1}$, where $\widehat{\Psi}%
_{n}=\sum_{i=1}^{n}\underline{\widehat{\ryskus{X}}}_{i}^{\statT}\Delta
_{i}\ryskus{V}_{i}^{-1}\Delta_{i}\underline{\widehat{\ryskus{X}}}_{i}$,
$%
\widehat{\Phi}_{n}=\sum_{i=1}^{n}\underline{\widehat{\ryskus{X}}}_{i}^{
\statT}\Delta_{i}\ryskus{V}_{i}^{-1}\widehat{\bSigma}_{i}%
\ryskus{V}_{i}^{-1}\Delta_{i}\underline{\widehat{\ryskus{X}}}_{i}$,\vspace*{1pt} and
$%
\widehat{\underline{\ryskus{X}}}_{i}=\underline{\ryskus{X}}_{i}-\operatorname{Proj}
_{G_{n}^{\ast}}\underline{\ryskus{X}}_{i}$, $i=1,\ldots,n$, in which Proj$
_{G_{n}^{\ast}}$ is the\vspace*{1pt} projection onto the empirically centered spline
inner product space and $\widehat{\bSigma}_{i}$ is a consistent
estimator of $\bSigma_{i}$.

For $1\leq l\leq d_{2}$, let $\bgamma_{l,0}^{\calS}=(
\gamma_{sl,0} )_{s=1}^{J_{n}^{\calS}}$, with $\gamma_{sl,0}$ defined in
the same fashion as given in (\ref{DEFglzl}), and
$\btheta_{-l0}=\{ \theta_{l^{\prime}0} ( \cdot)$,
$1\leq l^{\prime}\leq d_{2},l^{\prime}\neq l\} $. Define
\begin{eqnarray*}
\calD_{n,l}^{\ast} \bigl( \bgamma_{l}^{\calS}
\bigr) & = & \bigl\{ -\partial\ryskus{g}_{n,l}^{\ast} \bigl(
\bgamma_{l}^{%
\calS} \bigr) /\partial\bigl(
\bgamma_{l}^{\calS%
} \bigr)^{\statT} \bigr
\}_{J_{n}^{\calS}\times J_{n}^{\calS%
}},
\\
\Psi_{n,l}^{\ast} \bigl( \bgamma_{l,0}^{\calS}
\bigr) & = & \Biggl\{ \sum_{i=1}^{n} \bigl(
\ryskus{B}_{i\cdot l}^{\calS%
} \bigr)^{\statT}
\Delta_{i} \bigl( \bbeta_{0},\btheta_{-l0},
\bgamma_{l,0}^{\calS} \bigr) \ryskus{V}_{i}^{-1}
\bigl( \bbeta_{0},\btheta_{-l0},\bgamma_{l,0}^{\calS}
\bigr)
\\
&&\hspace*{109pt}{} \times\Delta_{i} \bigl( \bbeta_{0},
\btheta_{-l0},%
\bgamma_{l,0}^{\calS} \bigr)
\ryskus{B}_{i\cdot l}^{\calS}\Biggr\}_{J_{n}^{\calS}\times J_{n}^{\calS}}.
\end{eqnarray*}
In order to ensure the existence and uniformly weak convergence of the
oracle estimator $\widehat{\theta}_{n,l}^{\calS}(z_{l},\bbeta%
_{0},\btheta_{-l0})$, we need the following conditions:

\begin{longlist}[(A5)]
\item[(A5)] For $1\leq l\leq d_{2}$, there is a constant $c_{l}>0$, for
any $%
r>0$, such that $P \{ \calD_{n,l}^{\ast} ( \bgamma%
_{l}^{\calS} ) \geq c_{l}\Psi_{n,l}^{\ast} ( \bgamma_{l,0}^{\calS} )$
and $\calD_{n,l}^{\ast} (
\bgamma_{l}^{\calS} )$ is nonsingular, for all
$\bgamma_{l}^{\calS}
\in\xi_{n} ( r ) \} \rightarrow1$, where $\xi_{n} ( r ) = \{ \bgamma
_{l}^{\calS}\dvtx \llVert
\{ \Psi_{n,l}^{\ast} ( \bgamma_{l,0}^{\calS%
} ) \}^{1/2} ( \bgamma_{l}^{\calS}-\bgamma_{l,0}^{\calS} ) \rrVert
\leq
( \tau_{n}^{\max} )^{1/2}r \} $.
\end{longlist}

For $1\leq l\leq d_{2}$, define $\Xi_{n,l}^{\ast}= \{ E ( \Psi
_{n,l}^{\ast} ) \}^{-1}E ( \Phi_{n,l}^{\ast} )
\{ E ( \Psi_{n,l}^{\ast} ) \}^{-1}$, where
\[
\Phi_{n,l}^{\ast}=\sum
_{i=1}^{n} \bigl( \ryskus{B}_{i\cdot l}^{%
\calS}
\bigr)^{\statT}\Delta_{i}\ryskus{V}_{i}^{-1}
\Sigma_{i}%
\ryskus{V}_{i}^{-1}
\Delta_{i}\ryskus{B}_{i\cdot l}^{\calS},\qquad
\Psi_{n,l}^{\ast}=\sum_{i=1}^{n}
\bigl( \ryskus{B}_{i\cdot l}^{\calS%
} \bigr)^{\statT}
\Delta_{i}\ryskus{V}_{i}^{-1}
\Delta_{i}\ryskus{B}%
_{i\cdot l}^{\calS}.
\]
\end{remark}
%
\begin{theorem}
\label{THMthetahatnormality}Let $\theta_{l0}^{\ast}(z_{l})=E \{
\widehat{\theta}_{n,l}^{\calS}(z_{l},\bbeta_{0},\btheta_{-l0})\vert
\cal{X},\cal{Z} \} $. Under conditions \textup{(A3), (A4)(i)} and
\textup{(A5)}, for $1\leq l\leq d_{2}$ and
$z_{l}\in%
[ 0,1 ] $, as $n_{\statT}\rightarrow\infty$,
\begin{eqnarray*}\qquad
\bigl( \ryskus{B}_{l}^{\calS} ( z_{l}
)^{\statT}\Xi_{n,l}^{\ast}\ryskus{B}_{l}^{\calS}
( z_{l} ) \bigr)^{-1/2} \bigl\{ \widehat{\theta
}_{n,l}^{\calS}(z_{l},\bbeta%
_{0},
\btheta_{-l0})-\theta_{l0}^{\ast}(z_{l})
\bigr\} &\longrightarrow& N ( 0,1 ),\\[-20pt]
\end{eqnarray*}
%
\begin{eqnarray}\label{EQbias}
\sup_{z_{l}\in[ 0,1 ] }\bigl\llvert\widehat{\theta}_{n,l}^{%
\calS}(z_{l},
\bbeta_{0},\btheta_{-l0})-\theta_{l0}^{\ast}(z_{l})
\bigr\rrvert&=&O_{P} \Bigl\{ \sqrt{ ( \log n_{\statT%
} )
J_{n}^{\calS}/n_{\statT}} \bigl( \tau_{n}^{\max}/
\lambda_{n}^{\min} \bigr)^{1/2} \Bigr\},
\nonumber\hspace*{-35pt}\\[-4pt]\\[-12pt]
\sup_{z_{l}\in[ 0,1 ] }\bigl\llvert\theta_{l0}^{\ast
}(z_{l})-
\theta_{l0}(z_{l})\bigr\rrvert&=&O_{P} \bigl\{
\bigl( \lambda_{n}^{\max}/\lambda_{n}^{\min}
\bigr) \bigl( J_{n}^{\calS} \bigr)^{-p} \bigr\},\nonumber\hspace*{-35pt}
\end{eqnarray}
and there are constants $0<c_{l,\Xi}\leq C_{l,\Xi}<\infty$, such
that for
all $z_{l}\in[ 0,1 ] $,
%
\begin{eqnarray}\label{EQBGB}
\bigl\{ \ryskus{B}_{l}^{\calS} ( z_{l}
)^{\statT}\Xi_{n,l}^{\ast}\ryskus{B}_{l}^{\calS}
( z_{l} ) \bigr\}^{1/2} &\geq&c_{l,\Xi}\sqrt
{J_{n}^{\calS}/n_{\statT}} \bigl(
\tau_{n}^{\min}/\lambda_{n}^{\max}
\bigr)^{1/2},
\nonumber\\[-8pt]\\[-8pt]
\bigl\{ \ryskus{B}_{l}^{\calS} ( z_{l}
)^{\statT}\Xi_{n,l}^{\ast}\ryskus{B}_{l}^{\calS}
( z_{l} ) \bigr\}^{1/2} &\leq&C_{l,\Xi}\sqrt
{J_{n}^{\calS}/n_{\statT}} \bigl(
\tau_{n}^{\max}/\lambda_{n}^{\min}
\bigr)^{1/2}.\nonumber
\end{eqnarray}
Replacing \textup{(A4)(i)} by \textup{(A4$^{*}$)(i)}, one has $\sup_{z_{l}\in[ 0,1%
] }\llvert\widehat{\theta}_{n,l}^{\calS}(z_{l},%
\bbeta_{0},\btheta_{-l0})-\theta_{l0}^{\ast}(z_{l})\rrvert
=O_{P} \{ \sqrt{ ( \log n_{\statT} ) J_{n}^{\calS%
}m_{ ( n ) }/n_{\statT}} ( \lambda_{n}^{\max}/\lambda_{n}^{\min}
)^{1/2} \} $.
\end{theorem}
%
\begin{remark}\label{remark3}
Pointwise confidence intervals for $\theta_{l0}(z_{l})$
can be constructed based on the results in Theorem~\ref{THMthetahatnormality}.
By (\ref{EQbias}) and (\ref{EQBGB}), the bias
term in (\ref{EQbias}) is~asymptotically uniformly negligible through
undersmoothing if\break
$ ( N^{\calS} )^{-p-1/2}n_{\statT}^{1/2} ( \lambda_{n}^{\max}/\allowbreak\lambda
_{n}^{\min} ) ( \lambda_{n}^{\max}/\tau_{n}^{\min} )^{1/2}=o ( 1 )
$. Thus, $N^{\calS}$ is of
the form $ [ ( \lambda_{n}^{\max}/\lambda_{n}^{\min} )^{2} ( \lambda
_{n}^{\max}/\tau_{n}^{\min} )\* n_{\statT}%
]^{1/ ( 2p+1 ) }N^{\ast}$, where the sequence $N^{\ast}$\vspace*{1pt}
satisfies\break $N^{\ast}\rightarrow\infty$ and $n_{\statT}^{-\tau
}N^{\ast}\rightarrow0$ for any $\tau>0$. Under (A4$^*$)(i), $N^{\calS%
} $ is of the form $ [ ( \lambda_{n}^{\max}/\lambda_{n}^{\min
} )^{3}n_{\statT} ]^{1/ ( 2p+1 ) }N^{\ast}$.

Theorem~\ref{THMthetahatnormality} presents asymptotic normality and
uniform convergence rate of the oracle estimator $\widehat{\theta
}_{n,l}^{%
\calS}(z_{l},\bbeta_{0},\btheta_{-l0})$. The oracle
estimator achieves the convergence rate of univariate spline regression
function estimation. References~\cite{ZSW98} and~\cite{H03} studied asymptotic
normality of spline estimators for nonparametric regression functions with
i.i.d. data. Reference~\cite{HZZ06} established the asymptotic distribution for the
univariate spline estimator in partially linear models for clustered data
with $m_{ ( n ) }<\infty$. Reference~\cite{H03} discussed the difficulty of
obtaining asymptotic normality of spline estimators for additive models.
Reference~\cite{MSW12} studied convergence rate of the one-step additive spline
estimator for clustered data with $m_{ ( n ) }<\infty$, but it
lacks the limiting distribution. The next\vspace*{1pt} theorem will present the uniform
convergence rate of the two-step spline estimator $\widehat{\theta
}_{n,l}^{%
\calS}(z_{l},\widehat{\bbeta},\widehat{\btheta}%
_{n,-l})$ to the oracle estimator $\widehat{\theta}_{n,l}^{\calS%
}(z_{l},\bbeta_{0},\btheta_{-l0})$, and establish the
asymptotic normality of $\widehat{\theta}_{n,l}^{\calS}(z_{l},%
\widehat{\bbeta},\widehat{\btheta}_{n,-l})$.
\end{remark}
%
\begin{theorem}
\label{THMthetahatconvergence}Under conditions \textup{(A2)--(A5)}, for $1\leq
l\leq
d_{2}$,
%
\begin{eqnarray}
\label{EQuniformrate}
&&
\sup_{z_{l}\in[ 0,1 ] }\bigl\llvert\widehat
{\theta}_{n,l}^{\calS}(z_{l},\widehat{\bbeta},
\widehat{%
\btheta}_{n,-l})-\widehat{\theta
}_{n,l}^{\calS}(z_{l},%
\bbeta_{0},\btheta_{-l0})\bigr\rrvert
\nonumber\\
&&\qquad=O_{p} \bigl\{ \bigl( \lambda_{n}^{\max}/
\lambda_{n}^{\min} \bigr)^{2} \bigl( \sqrt{\log
n_{\statT}/n_{\statT}}+J_{n}^{-p} \bigr) \bigr
\} \\
&&\qquad=o_{p} \bigl\{ \bigl( J_{n}^{\calS}/n_{\statT}
\bigr)^{1/2} \bigl( \tau_{n}^{\min}/
\lambda_{n}^{\max} \bigr)^{1/2} \bigr\}
\nonumber
\end{eqnarray}
and replacing \textup{(A4)} by \textup{(A4$^{*}$)},
\begin{eqnarray*}
&&
\sup_{z_{l}\in[ 0,1 ] }\bigl\llvert\widehat{\theta}_{n,l}^{%
\calS}(z_{l},
\widehat{\bbeta},\widehat{%
\btheta}_{n,-l})-\widehat{\theta
}_{n,l}^{\calS}(z_{l},\bbeta%
_{0},
\btheta_{-l0})\bigr\rrvert\\
&&\qquad=o_{p} \bigl\{ \bigl(
J_{n}^{\calS}/n_{\statT} \bigr)^{1/2} \bigl(
\lambda_{n}^{\min}/\lambda_{n}^{\max}
\bigr)^{1/2} \bigr\}.
\end{eqnarray*}
Hence, for $1\leq l\leq d_{2}$ and $z_{l}\in[ 0,1 ] $, as $n_{%
\statT}\rightarrow\infty$,
\[
\bigl( \ryskus{B}_{l}^{\calS} ( z_{l}
)^{\statT}\Xi_{n,l}^{\ast}\ryskus{B}_{l}^{\calS}
( z_{l} ) \bigr)^{-1/2} \bigl\{ \widehat{\theta
}_{n,l}^{\calS}(z_{l},\widehat{%
\bbeta},
\widehat{\btheta}_{n,-l})-\theta_{l0}^{\ast
}(z_{l})
\bigr\} \longrightarrow N ( 0,1 ).
\]
\end{theorem}
%
\begin{remark}\label{remark4}
Similarly as $\widetilde{\Xi}_{n}$ in (\ref{DEFgammatilda}%
), $\Xi_{n,l}^{\ast}$ is minimized when $\ryskus{V}_{i}=\Sigma_{i}$ for
all $1\leq i\leq n$, and in this case is equal to $ \{ E ( \Psi
_{n,l}^{\ast} ) \}^{-1}$. To construct a pointwise confidence
interval for $\theta_{l0}(z_{l})$ at $z_{l}\in[ 0,1 ] $, $\Xi
_{n,l}^{\ast}$ is consistently estimated by $\widehat{\Xi
}_{n,l}^{\ast}=%
\widehat{\Psi}_{n,l}^{\ast-1}\widehat{\Phi}_{n,l}^{\ast}\widehat
{\Psi}%
_{n,l}^{\ast-1}$, where\vspace*{1pt} $\widehat{\Psi}_{n,l}^{\ast}=\sum_{i=1}^{n} (
\ryskus{B}_{i\cdot l}^{\calS} )^{\statT}\Delta_{i}\ryskus{V}%
_{i}^{-1}\Delta_{i}\ryskus{B}_{i\cdot l}^{\calS}$ and $\widehat{\Phi}%
_{n,l}^{\ast}=\sum_{i=1}^{n} ( \ryskus{B}_{i\cdot l}^{\calS%
} )^{\statT}\Delta_{i}\ryskus{V}_{i}^{-1}\widehat{\bSigma%
}_{i}\ryskus{V}_{i}^{-1}\Delta_{i}\ryskus{B}_{i\cdot l}^{\calS}$.
Then under the assumption given in Remark~\ref{remark3}, for any $\alpha\in(
0,1 ) $, an asymptotic $100 ( 1-\alpha) \%$ pointwise
confidence interval for $\theta_{l0}(z_{l})$ is
%
\begin{equation}\label{EQCIthetahat}
\widehat{\theta}_{n,l}^{\calS}(z_{l},\widehat{
\bbeta}, \widehat{\btheta}_{n,-l})\pm z_{\alpha/2} \bigl(
\ryskus{B}_{l}^{%
\calS} ( z_{l} )^{\statT}
\Xi_{n,l}^{\ast}\ryskus{B}%
_{l}^{\calS}
( z_{l} ) \bigr)^{1/2}.
\end{equation}
\end{remark}
%
\begin{remark}\label{remark5}
By letting $N$ have order
$n_{\statT}^{1/ ( 2p ) }$, the difference in
(\ref{EQuniformrate}) is of uniform order $O_{P} \{ ( \lambda
_{n}^{\max}/\lambda_{n}^{\min} )^{2}\sqrt{\log
n_{\statT}/n_{\statT}} \} $. So undersmoothing is applied
to reduce the approximation error caused by the bias in the first step.
\end{remark}

\section{Simulation}
\label{SECSimulation}

In this section we conduct simulations to illustrate the finite-sample
behavior of the proposed GEE estimators for both normal and binary
responses. For each procedure, we consider three different working
correlation structures: independence (IND), exchangeable (EX) and first
order auto-correlation (AR(1)). For notation simplicity,
denote the two-step
spline estimator $\widehat{\theta}_{n,l}^{\calS}(z_{l},\widehat{%
\bbeta}_{n},\widehat{\btheta}_{n,-l})$ defined in
(\ref{EQTSEl}) as $\widehat{\theta}_{n,l}^{\calS\calS
}(z_{l})=\ryskus{B}%
_{l}^{\calS} ( z_{l} )^{\statT}\widehat{\bgamma}%
_{n,l}^{\calS\calS}$, and the oracle estimator $\widehat{\theta
}_{n,l}^{%
\calS}(z_{l},\bbeta,\btheta_{-l})$ in (\ref{EQORACLE}) as $%
\widehat{\theta}_{n,l}^{\mathrm{OR}}(z_{l})=\ryskus{B}_{l}^{\calS%
} ( z_{l} )^{\statT}\widehat{\bgamma}_{n,l}^{\mathrm{OR}%
}$. In the first step, the pilot estimators are obtained by an undersmoothed
spline procedure to reduce bias. By the order requirements of the
number of
interior knots, we select a relatively large $N$ by letting $N= [ 2n_{%
\statT}^{1/ ( 2p ) } ] $, where $ [ a ] $
denotes the nearest integer to $a$. In the second step, $N^{\calS}$ is
selected from the interval $I_{N^{\calS}}= [ [ a_{n} ],%
[ 5a_{n} ] ] $, $a_{n}= ( n_{\statT}\log n_{%
\statT} )^{1/ ( 2p+1 ) }$, minimizing the BIC
criterion%
%
\begin{equation}\label{EQBIC}
\operatorname{BIC} \bigl( N^{\calS} \bigr) =\log\bigl\{
2Q_{n,l}^{\ast} \bigl( \widehat{\bgamma}_{n,l}^{\calS}
\bigr) /n \bigr\} +J_{n}^{%
\calS}\log( n ) /n,
\end{equation}
where $Q_{n,l}^{\ast} ( \widehat{\bgamma}_{n,l}^{\calS%
} ) =2^{-1}\sum_{i=1}^{n} ( \underline{\ryskus{Y}}_{i}-
\underline{\widehat{\bmu}}_{i} )^{\statT}\ryskus{V}_{i}^{-1} (
\widehat{\bbeta}_{n},\widehat{\btheta}_{n,-l},%
\widehat{\bgamma}_{n,l}^{\calS} ) (
\underline{\ryskus{Y}}_{i}-\underline{\widehat{\bmu}}_{i} ) $ with $%
\underline{\widehat{\bmu}}_{i}=\mu( \underline{\ryskus{X}}_{i}%
\widehat{\bbeta}_{n}+\sum_{l^{\prime}=1,l^{\prime}\neq l}^{d_{2}}%
\widehat{\theta}_{n,l^{\prime}}(\ryskus{Z}_{il^{\prime}})+ ( \ryskus
{B%
}_{i\cdot l}^{\calS} )^{\statT}\widehat{\bgamma}%
_{n,l}^{\calS} ) $. The optimal number of interior knots $N^{%
\calS}$ is chosen as $\widehat{N}{}^{\calS}=\arg\min_{N^{\calS}\in
I_{N^{\calS}}}$BIC$ ( N^{\calS} ) $. We use
cubic B-splines ($q=3$) to estimate the additive nonparametric
functions. We
generate $\mathrm{nsim}=500$ replications for each simulation study.

Given $1\leq l\leq d_{2}$, to compare the performance of the two-step
estimator $\widehat{\theta}_{n,l}^{\calS\calS}(z_{l})$ with the pilot
spline estimator $\widehat{\theta}_{n,l}(z_{l})$ and the oracle
estimator $%
\widehat{\theta}_{n,l}^{\mathrm{OR}}(z_{l})$, we define the mean integrated
squared error (MISE) for\break $\widehat{\theta}_{n,l}^{\calS\calS
}(z_{l})$ as
MISE$ ( \widehat{\theta}_{n,l}^{\calS\calS} ) =\frac{1}{\mathrm{%
nsim}}\sum_{\alpha=1}^{\mathrm{nsim}}\operatorname{ISE} ( \widehat{\theta}%
_{n,l,\alpha}^{\calS\calS} ) $, where\vspace*{1pt} ISE$ ( \widehat{\theta}%
_{n,l,\alpha}^{\calS\calS} ) =\break n_{\statT}^{-1}\sum%
_{i=1}^{n}\sum_{j=1}^{m_{i}}
( \widehat{\theta}_{n,l,\alpha}^{\calS\calS}(Z_{ijl,\alpha
})-\theta_{l}(Z_{ijl,\alpha}) )^{2}$, and $\widehat{\theta}%
_{n,l,\alpha}^{\calS\calS}$ is the estimator of $\theta_{l}$ and $%
Z_{ijl,\alpha}$ is the observation\vadjust{\goodbreak} of $Z_{ijl}$ in the $\alpha$th sample.
The MISEs for $\widehat{\theta}_{n,l}(z_{l})$ and $\widehat{\theta
}_{n,l}^{%
\mathrm{OR}}(z_{l})$ denoted as MISE$ ( \widehat{\theta}_{n,l} ) $
and MISE$ ( \widehat{\theta}_{n,l}^{\mathrm{OR}} ) $ are defined in
the same way. The empirical relative efficiency for the two-step estimator
in the $\alpha$th sample is defined as $\mathrm{eff}_{l,\alpha}= \{
\operatorname{ISE} ( \widehat{\theta}_{n,l,\alpha}^{\calS\calS} ) /%
\operatorname{ISE} ( \widehat{\theta}_{n,l,\alpha}^{\mathrm{OR}} ) \}
^{1/2}$. To construct confidence intervals for coefficient parameters $%
( \beta_{0,0},\ldots,\beta_{0, ( d_{1}-1 ) } ) $ by
using the first result in Theorem~\ref{THMbetahatnormality} and to construct
pointwise confidence intervals for the $l$th nonparametric function
$\theta_{l0}(z_{l})$ given in (\ref{EQCIthetahat}), the true
correlation matrix $%
\overline{\ryskus{R}}$ is consistently estimated by
\begin{eqnarray*}
\widehat{\ryskus{R}} &=&n^{-1}\sum_{i=1}^{n}
\ryskus{A}%
_{i}^{-1/2} ( \widehat{\bbeta
}_{n},\widehat{%
\bgamma}_{n} ) \bigl[
\underline{\ryskus{Y}}_{i}-\mu\bigl\{ \underline{\widetilde{
\bldeta}}_{i} ( \widehat{\bbeta}_{n}%
,\widehat{
\bgamma}_{n} ) \bigr\} \bigr]
\\
&&\hspace*{31pt}{}\times\bigl[ \underline{\ryskus{Y}}_{i}-\mu\bigl\{
\underline{\widetilde{ \bldeta}}_{i} ( \widehat{\bbeta
}_{n},\widehat{%
\bgamma}_{n} ) \bigr\}
\bigr]^{\statT}\ryskus{A}%
_{i}^{-1/2} (
\widehat{\bbeta}_{n},\widehat{%
\bgamma}_{n} ).
\end{eqnarray*}
And the covariance matrix $\Sigma_{i}$ is estimated by $\widehat{\Sigma
}%
_{i}=\ryskus{A}_{i}^{1/2}\widehat{\ryskus{R}}\ryskus{A}_{i}^{1/2}$. Let
$%
\bbeta_{0}= ( \beta_{0,k} )_{k=0}^{ ( d_{1}-1 )
}$and $\widehat{\bbeta}_{n}= ( \widehat{\beta}_{n,k} )_{k=0}^{ (
d_{1}-1 ) }$. For evaluating estimation accuracy of each
coefficient parameter,\vspace*{1pt} we report the root mean squared error (RMSE) defined
as $ \{ \sum_{\alpha=1}^{\mathrm{nsim}} ( \widehat{\beta}%
_{n,k}^{\alpha}-\beta_{0,k} )^{2}/\mathrm{nsim} \}^{1/2}$, for $%
0\leq k\leq d_{1}-1$, where $\widehat{\beta}_{n,k}^{\alpha}$ is the
estimate of $\beta_{0,k}$ obtained from the $\alpha$th
sample.
%
\begin{example}[(Continuous response)]\label{example1}
The correlated normal responses\vspace*{1pt}
are generated from the model $Y_{ij}=\ryskus{X}_{ij}^{\statT}
\bbeta+\theta_{1}(Z_{ij1})+\theta_{2}(Z_{ij2})+\theta
_{3}(Z_{ij3})+\varepsilon_{ij}$, where $\bbeta= (
1,-1,0.5 ) $, $\ryskus{X}_{ij}= ( X_{ij,1},X_{ij,2},X_{ij,3} )^{\statT
}$, $\theta_{l}(Z_{l})=\sin( 2\pi Z_{l} ) $, $1\leq
l\leq3$.\vspace*{1pt} For the covariates, let $Z_{ijl}=\Phi( Z_{ijl}^{\ast
} ) $, $1\leq l\leq3$, with $\ryskus{Z}_{ij}^{\ast}= (
Z_{ij1}^{\ast},Z_{ij2}^{\ast},Z_{ij3}^{\ast} )^{\statT}$\vspace*{1pt}
generated from the multivariate normal distribution with mean $0$ and an
AR(1) covariance with marginal variance $1$ and autocorrelation
coefficient $%
0.5$, $X_{ij,1}=\pm1/2$ with probability $1/2$, and $ (
X_{ij,2},X_{ij,3} )^{\statT}\sim\mathrm{N} [ (
0,0 )^{\statT},\break\operatorname{diag} ( a ( Z_{ij1} ),a (
Z_{ij2} ) ) ] $ with $a ( z ) =\frac{5-0.5\sin
( 2\pi z ) }{5+0.5\sin( 2\pi z ) }$. The error term\vspace*{1pt} $%
\underline{\bvarepsilon}_{i}= ( \varepsilon_{i1},\ldots,\break
\varepsilon_{im_{i}} )^{\statT}$ is generated from the
multivariate normal distribution with mean $0$, marginal variance $1$
and an
exchangeable correlation matrix with parameter $\rho=0.5$. We let $n=250$
and cluster size $m_{i}=m=20,50,100$, respectively. For computational
simplicity, we choose the same cluster size for each subject. The
computational algorithm can be easily extended to the case with varying
cluster sizes. Table~\ref{TABcontinuous} lists the empirical coverage rates
of the $95\%$ confidence intervals of the estimators $ ( \widehat{\beta
}%
_{n,k} )_{k=1}^{3}$ for coefficients $ ( \beta_{0,k} )_{k=1}^{3}$, the
RMSE and the absolute value of the empirical bias denoted
as Bias for IND, EX and AR(1) and $m=20,50,100$.

\begin{table}
\tabcolsep=0pt
\caption{The empirical coverage rates of the $95\%$ confidence
intervals for $( \protect\beta_{0,k} )_{k=1}^{3}$, the RMSE and Bias
for the IND, EX and AR(1) working correlation structures with
$m=20,50,100$} \label{TABcontinuous}
\begin{tabular*}{\tablewidth}{@{\extracolsep{\fill}}lcccccccccc@{}}
\hline
& & \multicolumn{3}{c}{\textbf{Coverage frequency}}
& \multicolumn{3}{c}{\textbf{RMSE}} & \multicolumn{3}{c@{}}{\textbf{Bias}} \\[-4pt]
& & \multicolumn{3}{c}{\hrulefill}
& \multicolumn{3}{c}{\hrulefill} & \multicolumn{3}{c@{}}{\hrulefill} \\
$\bolds{m}$ & & \multicolumn{1}{c}{$\bolds{\beta_{0,1}}$}
& \multicolumn{1}{c}{$\bolds{\beta_{0,2}}$}
& \multicolumn{1}{c}{$\bolds{\beta_{0,3}}$} & \multicolumn{1}{c}{$\bolds{\beta_{0,1}}$} &
\multicolumn{1}{c}{$\bolds{\beta_{0,2}}$} & \multicolumn{1}{c}{$\bolds{\beta_{0,3}}$}
& \multicolumn{1}{c}{$\bolds{\beta_{0,1}}$} & \multicolumn{1}{c}{$\bolds{\beta_{0,2}}$} &
\multicolumn{1}{c@{}}{$\bolds{\beta_{0,3}}$} \\
\hline
\hphantom{0}$20$ & IND & $0.948$ & $0.956$ & $0.950$ & $0.0279$ & $0.0137$ & $0.0137$ & $
0.0050$ & $0.0002$ & $0.0008$ \\
& EX & $0.954$ & $0.950$ & $0.948$ & $0.0196$ & $0.0098$ &
$0.0108$ & $%
0.0018$ & $0.0000$ & $0.0006$ \\
& AR(1) & $0.936$ & $0.954$ & $0.956$ & $0.0260$ & $0.0123$ & $0.0121$
& $%
0.0026$ & $0.0003$ & $0.0011$ \\
[6pt]
\hphantom{0}$50$ & IND & $0.948$ & $0.952$ & $0.948$ & $0.0177$ & $0.0092$ & $0.0091$ & $
0.0006$ & $0.0001$ & $0.0009$ \\
& EX & $0.946$ & $0.950$ & $0.948$ & $0.0126$ & $0.0063$ &
$0.0066$ & $%
0.0002$ & $0.0001$ & $0.0002$ \\
& AR(1) & $0.944$ & $0.956$ & $0.948$ & $0.0157$ & $0.0079$ & $0.0081$
& $%
0.0003$ & $0.0002$ & $0.0003$ \\
[6pt]
$100$ & IND & $0.948$ & $0.956$ & $0.958$ & $0.0126$ & $0.0063$ & $0.0064$ & $
0.0001$ & $0.0003$ & $0.0002$ \\
& EX & $0.950$ & $0.954$ & $0.948$ & $0.0084$ & $0.0044$ &
$0.0045$ & $%
0.0001$ & $0.0002$ & $0.0001$ \\
& AR(1) & $0.946$ & $0.954$ & $0.956$ & $0.0111$ & $0.0056$ & $0.0055$
& $%
0.0001$ & $0.0004$ & $0.0001$ \\ \hline
\end{tabular*}
\end{table}
%
%
\begin{table}[b]
\caption{The MISE$ ( \times10^{-3} ) $ for $\protect\widehat{%
\protect\theta}_{n,l}^{\calS\calS}(\cdot)$, $\protect\widehat
{\protect%
\theta}_{n,l}(\cdot)$ and $\protect\widehat{\protect\theta
}_{n,l}^{\mathrm{OR}}(\cdot)$, $l=1,2,3$, for the IND, EX and AR(1)
working correlation structures with $m=20,50,100$}
\label{TABcontinuousMISE}
\begin{tabular*}{\tablewidth}{@{\extracolsep{\fill}}lcccccccccc@{}}
\hline
$\bolds{m}$ & & $\bolds{\widehat{\theta}_{n,1}^{\calS\calS}}$
& $\bolds{\widehat{\theta}_{n,1}}$ & $\bolds{%
\widehat{\theta}_{n,1}^{\mathrm{OR}}}$ & $\bolds{\widehat{\theta}_{n,2}^{\calS\calS}}$
& $\bolds{%
\widehat{\theta}_{n,2}}$ & $\bolds{\widehat{\theta}_{n,2}^{\mathrm{OR}}}$ &
$\bolds{\widehat{%
\theta}_{n,3}^{\calS\calS}}$ & $\bolds{\widehat{\theta}_{n,3}}$ & $\bolds{\widehat{\theta
}_{n,3}^{%
\mathrm{OR}}}$ \\ \hline
\hphantom{0}$20$ & IND & $1.678$ & $2.231$ & $1.588$ & $1.659$ & $2.278$ & $1.517$ & $1.516$
& $2.118$ & $1.448$ \\
& EX & $0.883$ & $1.228$ & $0.836$ & $0.943$ & $1.232$ & $0.848$ &
$%
0.849$ & $1.167$ & $0.811$ \\
& AR(1) & $1.249$ & $1.710$ & $1.186$ & $1.324$ & $1.790$ & $1.205$ & $1.252$
& $1.713$ & $1.182$ \\
[6pt]
\hphantom{0}$50$ & IND & $0.633$ & $0.862$ & $0.601$ & $0.677$ & $0.927$ & $0.608$ & $0.631$
& $0.881$ & $0.601$ \\
& EX & $0.342$ & $0.463$ & $0.328$ & $0.348$ & $0.475$ & $0.321$ &
$%
0.353$ & $0.465$ & $0.335$ \\
& AR(1) & $0.473$ & $0.664$ & $0.459$ & $0.513$ & $0.690$ & $0.478$ & $0.486$
& $0.679$ & $0.464$ \\
[6pt]
$100$ & IND & $0.319$ & $0.440$ & $0.306$ & $0.346$ & $0.461$ & $0.317$ & $0.315$
& $0.436$ & $0.299$ \\
& EX & $0.173$ & $0.234$ & $0.166$ & $0.176$ & $0.237$ & $0.162$
& $%
0.172$ & $0.227$ & $0.164$ \\
& AR(1) & $0.247$ & $0.333$ & $0.235$ & $0.252$ & $0.348$ & $0.230$ & $0.244$
& $0.338$ & $0.232$ \\
\hline
\end{tabular*}
\end{table}

The empirical coverage rates are close to the nominal coverage probabilities
$95\%$ for all cases. The results are confirmative to Theorem~\ref{THMbetahatnormality}. EX has the smallest RMSE, since it is the true
correlation structure, which leads to the most efficient estimators (Remark
\ref{remark2}). The RMSEs decrease as cluster size increases for all three working
correlation structures. The last three columns show that the empirical
biases are close to zero for all cases.

Table~\ref{TABcontinuousMISE} shows the MISE$( \mbox{$\times$}10^{-3} ) $
for the two-step spline estimator $\widehat{\theta}_{n,l}^{\calS\calS
}(\cdot)$, the pilot estimator $\widehat{\theta}_{n,l}(\cdot)$ and the
oracle estimator $\widehat{\theta}_{n,l}^{\mathrm{OR}}(\cdot)$, $l=1,2,3$,
for IND, EX and AR(1) structures and cluster size $m=20,50,100$.
$\widehat{%
\theta}_{n,l}^{\calS\calS}(\cdot)$ and $\widehat{\theta
}_{n,l}^{\mathrm{OR%
}}(\cdot)$ have similar MISE values, while $\widehat{\theta
}_{n,l}(\cdot)$
has the largest MISE value. The EX structure has the smallest MISEs,
and the
MISEs decrease as the cluster size increases.

\begin{figure}

\includegraphics{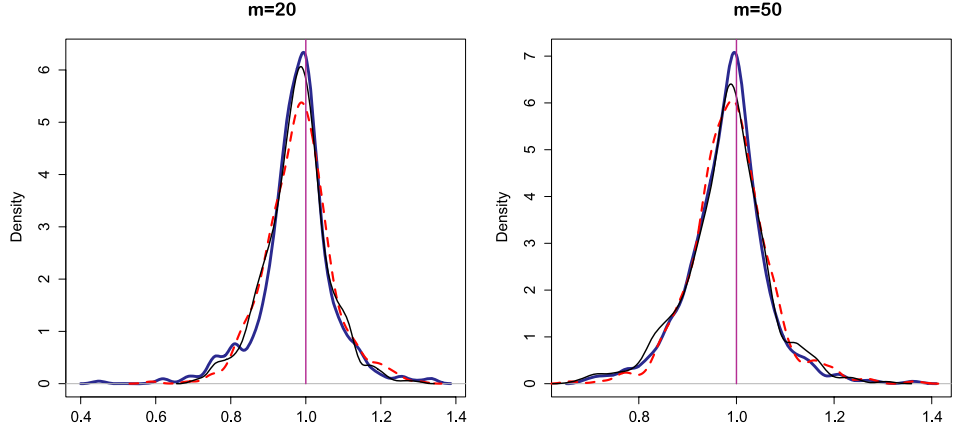}

\caption{Kernel density plots of the $500$ empirical efficiencies of the
two-step estimator to the oracle estimator of the first function
$\protect%
\theta_{1}(\cdot)$ for the IND (dashed lines), EX (thick lines) and AR(1)
(thin lines) working correlation structures with $m=20,50$.}
\label{FIGdensitycon}
\end{figure}

We plotted the kernel density estimates in Figure~\ref{FIGdensitycon}
of $%
500$ empirical efficiencies $\mathrm{eff}_{l,\alpha}$ for the
estimators of
the first function $\theta_{1}(\cdot)$ for IND (dashed lines), EX (thick
lines) and AR(1) (thin lines) structures with $m=20,50$ and $n=250$. The
vertical line at $\mbox{efficiency}=1$ is the standard line for the
comparison of
the two-step estimator (\ref{EQTSEl}) and the oracle estimator (\ref{EQORACLE}). The centers of density distributions are close to $1$ for all
working correlation structures, and EX has the narrowest distribution.
\end{example}
%
\begin{example}[(Binary response)]\label{example2}
The correlated binary responses $\{ Y_{ij} \} $ are
generated from a marginal logit model%
\[
\operatorname{logit}P ( Y_{ij}=1\vert\ryskus{X}_{ij},
\ryskus{Z}%
_{ij} ) =\ryskus{X}_{ij}^{\statT}
\bbeta+\theta_{1}(Z_{ij1})+\theta_{2}(Z_{ij2}),
\]
where $\bbeta= ( 0.5,-0.3,0.3 ) $, $\ryskus{X}_{ij}= (
1,X_{ij,1},X_{ij,2} )^{\statT}$, $\theta_{1}(Z_{1})=0.5\times
\sin( 2\pi Z_{1} ) $, and $\theta_{2}(Z_{2})=-0.5\times\{
Z_{2}-0.5+\sin( 2\pi Z_{2} ) \} $. For the covariates, we
generate $X_{ijk}$ and $Z_{ijl}$ independently from standard normal and
uniform distributions, respectively, such that $X_{ijk}\sim\mathrm{N} (
0,1 ) $ and $Z_{ijl}\sim\operatorname{Uniform} [ 0,1 ] $. We use
the R package ``mvtBinaryEP'' to generate the correlated binary responses with
exchangeable correlation structure with a correlation parameter of $0.1$
within each cluster. We let the number of clusters be $n=100,200,500$,
respectively, and let the cluster size be equal and increase with $n$, such
that $m_{ ( n ) }=m_{i}= \lfloor2n^{1/2} \rfloor$, for $%
1\leq i\leq n$, where $ \lfloor a \rfloor$ denotes the largest
integer no greater than $a$. So $m=20,28,44$ for $n=100,200,500$,
respectively. Table~\ref{TABbinary} shows the empirical coverage rates of
the $95\%$ confidence intervals of the estimators $ ( \widehat{\beta}%
_{n,k} )_{k=0}^{2}$ for the coefficients $ ( \beta_{0,k} )_{k=0}^{2}$
and the RMSEs for IND, EX and AR(1) and $n=100,200,500$. Table
\ref{TABbinaryMISE} shows that the empirical coverage rates are close to
the nominal coverage probabilities $95\%$ for all cases. EX has the smallest
RMSE values, and the RMSEs decrease as $n$ increases.

\begin{table}
\caption{The empirical coverage rates of the $95\%$ confidence
intervals for
$ ( \protect\beta_{0,k} )_{k=0}^{2}$ and the estimated MSE for
the IND, EX and AR(1) working correlation structures with $n=100,200,500$}
\label{TABbinary}
\begin{tabular*}{\tablewidth}{@{\extracolsep{\fill}}lccccccc@{}}
\hline
& & \multicolumn{3}{c}{\textbf{Coverage frequency}} & \multicolumn
{3}{c@{}}{\textbf{RMSE}} \\[-4pt]
& & \multicolumn{3}{c}{\hrulefill} & \multicolumn
{3}{c@{}}{\hrulefill} \\
& & $\bolds{\beta_{0,0}}$ & $\bolds{\beta_{0,1}}$ & $\bolds{\beta_{0,2}}$
& $\bolds{\beta_{0,0}}$ & $\bolds{%
\beta_{0,1}}$ & $\bolds{\beta_{0,2}}$ \\ \hline
$n=100,m=20$ & IND & $0.960$ & $0.946$ & $0.940$ & $0.0821$ & $0.0549$ & $0.0506$ \\
& EX & $0.940$ & $0.946$ & $0.946$ & $0.0763$ & $0.0469$ &
$%
0.0454$ \\
& AR(1) & $0.966$ & $0.930$ & $0.940$ & $0.0773$ & $0.0540$ & $0.0488$
\\
[6pt]
$n=200,m=28$ & IND & $0.944$ & $0.946$ & $0.940$ & $0.0559$ & $0.0299$ & $0.0328$ \\
& EX & $0.948$ & $0.952$ & $0.942$ & $0.0554$ & $0.0289$ &
$%
0.0310$ \\
& AR(1) & $0.940$ & $0.950$ & $0.940$ & $0.0556$ & $0.0291$ & $0.0325$
\\
[6pt]
$n=500,m=44$ & IND & $0.952$ & $0.946$ & $0.942$ & $0.0340$ & $0.0157$ & $0.0154$ \\
& EX & $0.948$ & $0.952$ & $0.946$ & $0.0336$ & $0.0136$ &
$%
0.0142$ \\
& AR(1) & $0.952$ & $0.952$ & $0.942$ & $0.0340$ & $0.0153$ & $0.0153$\\
\hline
\end{tabular*}
\end{table}

%
\begin{table}[b]
\caption{The MISE for $\protect\widehat{\protect\theta}_{n,l}^{\calS
\calS%
}(\cdot)$, $\protect\widehat{\protect\theta}_{n,l}(\cdot)$ and
$\protect%
\widehat{\protect\theta}_{n,l}^{\mathrm{OR}}(\cdot)$, $l=1,2$, for
the IND,
EX and AR(1) working correlation structures with $n=100,200,500$}
\label{TABbinaryMISE}
\begin{tabular*}{\tablewidth}{@{\extracolsep{\fill}}lccccccc@{}}
\hline
$\bolds{n}$ & & $\bolds{\widehat{\theta}_{n,1}^{\calS\calS}}$
& $\bolds{\widehat{\theta}_{n,1}}$ & $\bolds{%
\widehat{\theta}_{n,1}^{\mathrm{OR}}}$ & $\bolds{\widehat{\theta}_{n,2}^{\calS\calS}}$
& $\bolds{%
\widehat{\theta}_{n,2}}$ & $\bolds{\widehat{\theta}_{n,2}^{\mathrm{OR}}}$ \\
\hline
$100$ & IND & $0.0172$ & $0.0243$ & $0.0174$ & $0.0158$ & $0.0222$ & $0.0159$
\\
& EX & $0.0148$ & $0.0223$ & $0.0148$ & $0.0139$ & $0.0204$ & $0.0137$
\\
& AR(1) & $0.0178$ & $0.0265$ & $0.0176$ & $0.0161$ & $0.0234$ &
$0.0163$ \\
[6pt]
$200$ & IND & $0.0059$ & $0.0086$ & $0.0059$ & $0.0056$ & $0.0082$ & $0.0056$
\\
& EX & $0.0048$ & $0.0069$ & $0.0048$ & $0.0054$ & $0.0075$ & $0.0053$
\\
& AR(1) & $0.0058$ & $0.0085$ & $0.0058$ & $0.0056$ & $0.0081$ &
$0.0056$ \\
[6pt]
$500$ & IND & $0.0015$ & $0.0022$ & $0.0015$ & $0.0015$ & $0.0021$ & $0.0015$
\\
& EX & $0.0013$ & $0.0019$ & $0.0013$ & $0.0014$ & $0.0019$ & $0.0013$
\\
& AR(1) & $0.0015$ & $0.0022$ & $0.0015$ & $0.0015$ & $0.0020$ &
$0.0014$ \\
\hline
\end{tabular*}
\end{table}

Table~\ref{TABbinaryMISE} shows the MISE for the two-step spline
estimator $%
\widehat{\theta}_{n,l}^{\calS\calS}(\cdot)$, the pilot estimator $%
\widehat{\theta}_{n,l}(\cdot)$ and the oracle estimator $\widehat
{\theta}%
_{n,l}^{\mathrm{OR}}(\cdot)$, $l=1,2$, for the IND, EX and AR(1) structures
and $n=100,200,500$. The MISE values for $\widehat{\theta
}_{n,l}^{\calS\calS}(\cdot)$ and $\widehat{\theta}_{n,l}^{\mathrm
{OR}}(\cdot)$ are close
and $\widehat{\theta}_{n,l}(\cdot)$ has the largest MISE values. EX has
the smallest MISEs among the three working correlation structures, and the
MISEs decrease as $n$ increases.

%
\begin{figure}

\includegraphics{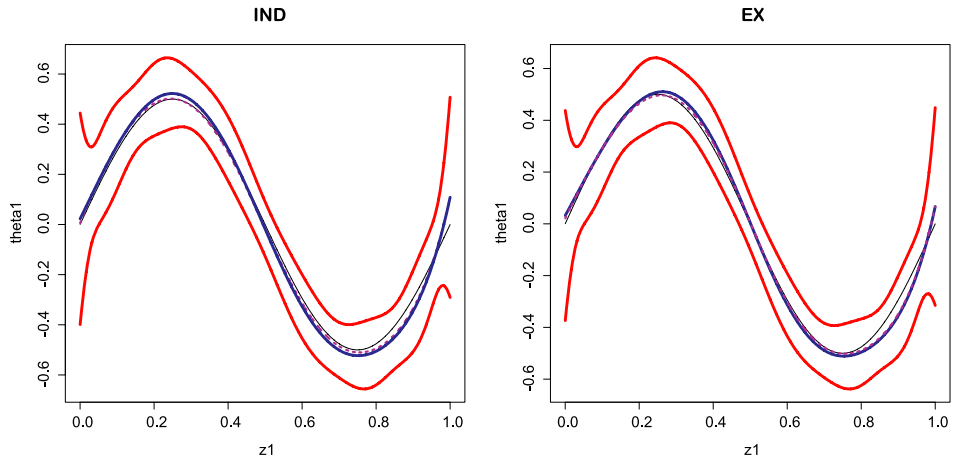}

\caption{Plots of oracle estimator (dashed curve), the two-step estimator
(thick curve) and the $95\%$ pointwise confidence intervals (upper and
lower curves) of $\protect\theta_{1}(\cdot)$ (thin curve) for $n=200$.}
\label{FIGestplot}
\end{figure}

For visualization of the actual function estimates, in Figure
\ref{FIGestplot} we plotted the oracle estimator given in
(\ref{EQORACLE}) (dashed curve), the two-step estimator given in
(\ref{EQTSEl}) (thick curve) and the $95\%$ pointwise confidence
intervals constructed in (\ref{EQCIthetahat}) (upper and lower curves)
of $\theta_{1} ( \cdot ) $ (thin curve) for $n=200$ based on one
simulated sample. The proposed two-step estimator seems satisfactory.
\end{example}

\section{Application}
\label{SECApplication}

In this section we apply the proposed estimation procedure to analyze
unemployment-economic growth and employment relationship at the U.S. state
level for the 1970--1986 period. Reference~\cite{BC05} has first studied the
effect of
economic growth on unemployment rate by establishing a parametric
unemployment-growth model. They concluded that relatively high economic
growth is more likely to show reduced unemployment rates when compared to
states in which the economy is growing more slowly by obtaining a negative
coefficient for growth. Reference~\cite{W99} demonstrated a strong negative
correlation between the change of unemployment rate and employment. We
restudy their relationship by considering possible nonlinear relations of
the unemployment rate with economic growth and time. The economic growth
rate is calculated from the logarithm difference of the gross state product
(GSP). The data for the unemployment rate, gross state product and
employment are available for the U.S. 48 contiguous states over the period
1970--1986. Details on this data set can be found in~\cite{M90}. The number
of time periods for each state in estimation is $m=16$, since the year 1970
is taken as the initial observation. We consider the following GAPLM:%
\[
U_{ij}=\beta_{0}+\beta_{1}E_{ij}+
\theta_{1}(T_{ij})+\theta_{2}(G_{ij})+
\varepsilon_{ij},\qquad j=2,\ldots,17,i=1,\ldots,48,
\]
where $U_{ij}$ is the change in the unemployment rate for the $j$th
year in
the $i$th state, $E_{ij}$ is the empirically centered value of the relative
change in employment, $G_{ij}$~is the GSP growth, and $T_{ij}$ is time.
$%
\theta_{1} ( \cdot) $ and $\theta_{2}(\cdot)$ are
nonparametric functions of time and GSP growth, respectively.

To test whether $\theta_{l} ( \cdot) $, $l=1,2$, has a specific
parametric form, we construct simultaneous confidence bands according to
Theorem 2 of~\cite{WY09}. For any $\alpha\in( 0,1 ) $, an
asymptotic $100 ( 1-\alpha) \%$ conservative confidence band for $%
\theta_{l0}(z_{l})$ over the domain of $z_{l}$ is given as
\[
\widehat{\theta}_{n,l}^{\calS}(z_{l},\widehat{
\bbeta},\widehat{\btheta}_{n,-l})\pm\bigl\{ 2\log\bigl(
N^{s}+1 \bigr) -2\log\alpha\bigr\}^{1/2} \bigl(
\ryskus{B}_{l}^{\calS} ( z_{l} )^{\statT}
\Xi_{n,l}^{\ast}\ryskus{B}_{l}^{\calS%
} (z_{l} ) \bigr)^{1/2}
\]
with $\widehat{\theta}_{n,l}^{\calS}$ obtained by linear splines with
degree $q=1$. We use linear splines in both steps of estimation.

We use three working correlation structures to analyze this data set,
including the working independence $\ryskus{R}_{i} ( \alpha) =%
\ryskus{I}_{m}$, where $\ryskus{I}_{m}$ is an $m\times m$ identity matrix,
the exchangeable $\ryskus{R}_{i} ( \alpha) =\alpha\times
1_{m}1_{m}^{\statT}+ ( 1-\alpha) \ryskus{I}_{m}$, where $1_{m}$
is the $m$-dimensional vector with $1$'s, and the AR(1) $\ryskus{R}%
_{i} ( \alpha) = ( R_{ijj^{\prime}} )_{j,j^{\prime
}=1}^{m}$ with $R_{ijj^{\prime}}=\alpha^{\llvert j-j^{\prime
}\rrvert}$.\vspace*{1pt} The parameter $\alpha$ is estimated by the R package
geepack from the first spline estimation step. We obtain the estimated
values for $\alpha$ which are $\widehat{\alpha}=0.088$ for the EX
structure and $\widehat{\alpha}=-0.199$ for the AR(1) structure,
respectively. Table~\ref{TABreal} shows the estimated
%
\begin{table}
\tablewidth=230pt
\caption{The estimated\vspace*{1pt} values $\protect\widehat{\protect\beta}_{0}$
and $%
\protect\widehat{\protect\beta}_{1}$ of $\protect\beta_{0}$ and
$\protect%
\beta_{1}$ and the standard errors SE$ ( \protect\widehat{\protect\beta
}_{0} ) $ and SE$ ( \protect\widehat{\protect\beta}_{1} ) $
for the IND, EX and AR(1) working correlation structures}
\label{TABreal}
\begin{tabular*}{\tablewidth}{@{\extracolsep{\fill}}lcccc@{}}
\hline
& $\bolds{\widehat{\beta}_{0}}$ & $\bolds{\operatorname{SE} ( \widehat{\beta}_{0} )}$
& $\bolds{\widehat{\beta}_{1}}$ & $\bolds{\operatorname{SE}( \widehat{\beta}_{1} ) }$ \\
\hline
IND & 0.127 & 0.0417 & $-$0.219 & 0.0230 \\
EX & 0.127 & 0.0494 & $-$0.249 & 0.0220 \\
AR(1) & 0.127 & 0.0484 & $-$0.250 & 0.0216 \\
\hline
\end{tabular*}
\end{table}
values $\widehat{
\beta}_{0}$ and $\widehat{\beta}_{1}$ of $\beta_{0}$ and $\beta_{1}$ and
the corresponding standard errors SE$ ( \widehat{\beta}_{0} ) $
and SE$ ( \widehat{\beta}_{1} ) $ for the three working
correlation structures. The estimation results are very similar for the
three structures. The negative values of $\widehat{\beta}_{1}$ imply a
negative relationship between $U_{ij}$ and $E_{ij}$, confirmative to the
result in~\cite{W99}. Both of the estimators are significant with $p$-values
close to $0$ for the three different working correlation structures. The
correlation coefficient $r=0.785$, $0.822$ and $0.762$ for the IND, EX and
AR(1) structures, respectively.

\begin{sidewaysfigure}
\centering
\includegraphics{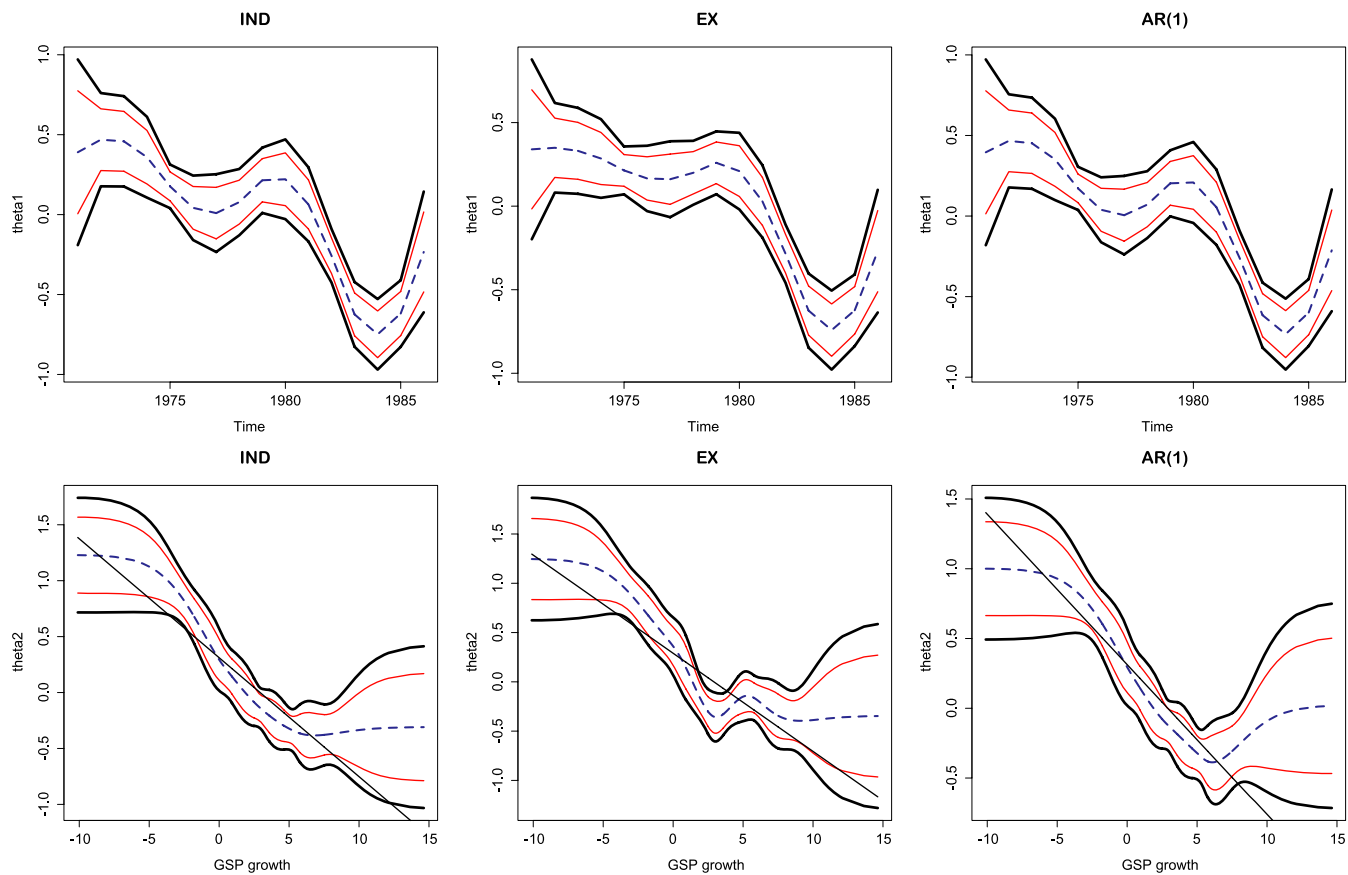}

\caption{Plots of the two-step spline estimated functions (dashed
line), the
$95\%$ pointwise confidence intervals (thin lines) and the $95\%$ confidence
bands (thick lines) for $\protect\theta_{1} ( \cdot) $ (upper
panel) and $\protect\theta_{2} ( \cdot) $ (lower panel), and the
GEE estimator of $\protect\theta_{2} ( \cdot) $ by assuming
linearity (straight solid line).}\vspace*{-5pt}
\label{FIGtheta1}
\end{sidewaysfigure}

Figure~\ref{FIGtheta1} displays the two-step spline estimators
$\widehat{%
\theta}_{n,1}^{\calS\calS} ( \cdot) $ (dashed lines) and $%
\widehat{\theta}_{n,2}^{\calS\calS} ( \cdot) $ (dashed lines)
of $\theta_{1} ( \cdot) $ and $\theta_{2} ( \cdot) $
and the corresponding $95\%$ pointwise confidence intervals (thin
lines) and
simultaneous confidence bands (thick lines) for the three working
structures. Figure~\ref{FIGtheta1} shows that the change patterns of $%
U_{ij} $ with $T_{ij}$ and $G_{ij}$ are very similar for the three working
structures. In the upper\vspace*{1pt} panel of Figure~\ref{FIGtheta1}, we can
observe a
declining trend for $\widehat{\theta}_{n,1}^{\calS\calS} ( \cdot
) $ in general. The values of $\widehat{\theta}_{n,1}^{\calS\calS%
} ( \cdot) $ were all positive before the year $1976$, which
means that the unemployment rate was increasing with time during that
period. The increasing unemployment rate was caused by a severe economic
recession that happened in the years 1973--1975. A local peak of
$\widehat{%
\theta}_{n,1}^{\calS\calS} ( \cdot) $ is observed around $1980$,
when another recession happened.

In order to test the linearity of the nonparametric function
$\theta_{2}$, we plotted straight solid lines in the lower panel of
Figure~\ref{FIGtheta1}, which are the regression lines obtained by
solving the GEE in (\ref{EQgnlstar}) by assuming that $\theta_{1} (
\cdot) $ is a linear function of GSP growth. All the three plots in the
lower panel of Figure~\ref{FIGtheta1} show that the confidence bands
with $95\%$ confidence level do not totally cover the straight
regression lines, that is, the linearity of the component function for
GSP growth is rejected at the significance level $0.05$. The lower
panel of Figure~\ref{FIGtheta1} indicates a general negative relation
between the GSP growth and the change in unemployment rate.

\section{Discussion}
\label{SECdiscussion}

In this paper we propose a two-step spline estimating
equations procedure for generalized additive partially linear models with
large cluster sizes. We develop asymptotic distributions and consistency
properties for the two-step estimators of the additive functions and the
one-step estimator of the parametric vector. We establish the oracle
properties of the two-step estimators. Because the two-step estimator
is a
mixture of two different spline bases, and an infinite number of observations
within clusters are correlated in complex ways, we encountered challenging
tasks when developing the theories. We demonstrate our proposed method by
two simulated examples and one real data example. Our proposed method
can be
extended to generalized additive models and generalized additive coefficient
models, and it provides a useful tool for studying clustered data. The
theoretical development in this paper helps us further investigate
semi-parametric models with clustered data. In the real data example, we
constructed confidence bands to test the linearity of the nonparametric
function. To establish confidence bands with rigorous theoretical proofs
will be our future work.

In this paper we focus on the two-step spline estimation procedure, which
is computationally expedient and theoretically reliable. As mentioned in
Section~\ref{SECMethodology}, that kernel smoothing method can be
applied to
the second step. Let $K_{h} ( \cdot) $ be a kernel weight
function, where $K_{h} ( z ) =h^{-1}K ( z/h ) $ with
bandwidth $h$. Let $G_{1} ( z_{l} ) = ( 1,z_{l} )^{\statT}$.
If we use local linear kernel estimation, then by assuming\vspace*{2pt} that $%
\bbeta$ and $\btheta_{-l}$ are known, $\theta_{l}(\cdot)$
is estimated by the oracle estimator $\widehat{\theta}_{l}^{\mathrm{OR}}
(Z_{l})=G_{1} ( Z_{l}-z_{l} )^{\statT}\widehat{\bgamma}%
_{l}^{\mathrm{OR}}$ at any given point $z_{l}$, where $\widehat{\bgamma
}_{l}^{\mathrm{OR}}= ( \widehat{\gamma
}_{l0}^{\mathrm{OR}},\widehat{\gamma
}_{l1}^{\mathrm{OR}} )^{\statT}$ with $\widehat{\bgamma}_{l}^{%
\mathrm{OR}}$ solving the kernel estimating equations%
\begin{eqnarray*}
&&\sum_{i=1}^{n}G_{i1} (
z_{l} )^{\statT}\Delta_{i} \bigl( \bbeta,
\btheta_{-l},\widehat{\bgamma}_{l}^{%
\mathrm{OR}} \bigr)
\ryskus{V}_{i}^{-1} \bigl( \bbeta,\btheta_{-l},%
\widehat{\bgamma}_{l}^{\mathrm{OR}} \bigr) \ryskus{K}_{ih}
( z_{l} )
\\
&&\qquad{}\times\Biggl\{ \underline{\ryskus{Y}}_{i}-\mu\Biggl( \underline{
\ryskus{X}}_{i}%
\bbeta+\sum_{l^{\prime}=1,l^{\prime}\neq
l}^{d_{2}}
\theta_{l^{\prime}}(\ryskus{Z}_{il^{\prime}})+G_{i1} (
z_{l} ) \widehat{\bgamma}_{l}^{\mathrm{OR}} \Biggr)
\Biggr\} =0,
\end{eqnarray*}
where $\ryskus{K}_{ih} ( z_{l} ) =\operatorname{diag} \{ K_{h} (
Z_{ijl}-z_{l} ) \} $ and $G_{i1} ( z_{l} ) = \{
G_{1} ( Z_{i1l}-z_{l} ),\ldots,\break
G_{1} ( Z_{im_{i}l}-z_{l} ) \}^{\statT}$. Then $%
\theta_{l}(z_{l})$ is estimated by $\widehat{\theta}_{l}^{\mathrm{OR}%
}(z_{l})=\widehat{\gamma}_{l0}^{\mathrm{OR}}$. The two-step spline backfitted
kernel (SBK) estimator $\widehat{\theta}_{l}^{\mathrm{SBK}}(z_{l})$ is
obtained by replacing $\bbeta$ and $\btheta_{-l}$ with the
pilot estimators $\widehat{\bbeta}_{n}$ and $\widehat{%
\btheta}_{n,-l}$ from step I. The asymptotic normality of the oracle
estimator $\widehat{\theta}_{l}^{\mathrm{OR}}(z_{l})$ which is a pure local
linear kernel estimator of $\theta_{l}(z_{l})$ by GEE can be obtained
following the same idea in the proofs for Theorem~\ref{THMthetahatnormality}
and the results in~\cite{LC00} for kernel estimators using GEE. The uniform
oracle efficiency of the SBK estimator $\widehat{\theta}_{l}^{\mathrm
{SBK}%
}(z_{l})$ is achievable by following the same procedure as the proofs for
Theorem~\ref{THMthetahatconvergence} and by studying the properties of
spline-kernel combination. See~\cite{WY07,LY10} and~\cite{MY11} for
the oracle properties of the SBK estimators in additive models, additive
coefficient models and additive partially linear models with
weekly-dependent data and a continuous response variable. The asymptotic
distributions and the oracle properties of the SBK estimators for GAPLMs
with large cluster sizes still need us to explore as future work.

\begin{appendix}\label{app}
\section*{Appendix}

We denote by the same letters $c,C$, any positive constants without
distinction. For any $s\times s^{\prime}$ matrix $\ryskus{M}$, let $%
\llVert\ryskus{M}\rrVert_{\infty}={\max_{1\leq i\leq
s}\sum_{j_{=1}}^{s^{\prime}}}\llvert M_{ij}\rrvert$. For
any vector $\balpha= ( \alpha_{1},\ldots,\alpha_{s} )^{%
\statT}$, denote$\llVert\balpha\rrVert_{\infty
}=\max_{1\leq i\leq s}\llvert\alpha_{i}\rrvert$ as the maximum
norm. Let $\ryskus{I}_{s}$ be the $s\times s$ identity matrix. Let
$\widehat{%
\Pi}_{n}$, $\Pi_{n}$ denote, respectively, the projection onto
$G_{n}^{0}$ relative to the empirical and the theoretical inner
products. For any\vspace*{1pt} function $\phi$, define the empirical
norm as $\llVert\phi\rrVert
_{n_{\statT}}^{2}=n_{\statT}^{-1}\sum_{i=1}^{n}%
\sum_{j=1}^{m_{i}}\phi( X_{ij},Z_{ij} )^{2}$. For positive
numbers $a_{n}$ and $b_{n}$, let $a_{n}\asymp b_{n}$ denote that $%
\lim_{n\rightarrow\infty}a_{n}/b_{n}=c$, where $c$ is some nonzero
constant.

\subsection{\texorpdfstring{Proof of Theorem \protect\ref{THMweakconsistency}}{Proof of Theorem 1}}

It can be proved following the similar reasoning as in~\cite{MSW12} that
under condition (A1) with $n_{\statT}\rightarrow\infty$, $%
J_{n}\rightarrow\infty$, and $J_{n}n^{-1}=o ( 1 ) $, there exist
constants $0<c^{\prime}<C^{\prime}<\infty$, such that with
probability $1$,
for $n_{\statT}$ sufficiently large,
\[
c^{\prime}n_{\statT}\leq\lambda_{\mathrm{min}} \Biggl( \sum
_{i=1}^{n}\underline{\ryskus{B}}_{i}^{\statT}
\underline{\ryskus{B}}_{i} \Biggr) \leq\lambda_{\mathrm{max}}
\Biggl( \sum_{i=1}^{n}\underline{
\ryskus{B}}_{i}^{\statT} \underline{\ryskus{B}}_{i}
\Biggr) \leq C^{\prime}n_{\statT}
\]
and $\llVert{\sum_{i=1}^{n}\underline{\ryskus{X}}_{i}^{\statT
}\underline{\ryskus{B}}_{i}}\rrVert_{\infty}=O_{\mathrm{a.s.}} \{
( n_{\statT}\log n_{\statT} )^{1/2} \}
$. By these results together with condition (C4), one has with
probability $%
1$,
%
\begin{equation}\label{EQDD}
c^{\prime\prime}n_{\statT}\leq\lambda_{\mathrm{min}} \Biggl( \sum
_{i=1}^{n}\underline{\ryskus{D}}_{i}^{\statT}%
\underline{\ryskus{D}}_{i} \Biggr) \leq\lambda_{\mathrm{max}}
\Biggl( \sum_{i=1}^{n}\underline{
\ryskus{D}}_{i}^{\statT} \underline{\ryskus{D}}_{i}
\Biggr) \leq C^{\prime\prime}n_{\statT}
\end{equation}
for some constants $0<c^{\prime\prime}<C^{\prime\prime}<\infty$. Then
by condition (A2),
\[
\bigl( \tau_{n}^{\max} \bigr)^{-1}
\lambda_{\mathrm{min}} \bigl\{ \Psi_{n} ( \bbeta_{0},
\bgamma_{0} ) \bigr\} \geq cc^{\prime\prime} \bigl(
\tau_{n}^{\max} \bigr)^{-1}\lambda_{n}^{\min}n_{\statT}
\rightarrow\infty.
\]
Results in Theorem~\ref{THMweakconsistency} can be proved similarly as
Theorems 1 and 2 in~\cite{XY03} with $r=\sqrt{2 (
d_{1}+d_{2}J_{n} ) /c_{0}\varepsilon}$ for any given $\varepsilon>0$.

\subsection{\texorpdfstring{Proof of Theorem \protect\ref{THMbetahatnormality}}{Proof of Theorem 2}}

By Taylor's expansion, one has
%
\begin{equation}\label{EQbetagammahat}
\ryskus{g}_{n} ( \widehat{\bbeta}_{n},
\widehat{%
\bgamma}_{n} ) -\ryskus{g}_{n} (
\bbeta_{0} ,\bgamma_{0} ) =-\calD_{n} \bigl(
\bbeta_{n}^{\ast}%
,\bgamma_{n}^{\ast}
\bigr) \pmatrix{ \widehat{\bbeta}_{n}-\bbeta_{0}
\cr
\widehat{\bgamma}_{n}-\bgamma_{0}},
\end{equation}
where $\bbeta_{n}^{\ast}=t_{1}\widehat{\bbeta}_{n}+ (
1-t_{1} ) \bbeta_{0}$, and $\bgamma_{n}^{\ast}=t_{2}%
\widehat{\bgamma}_{n}+ ( 1-t_{2} ) \bgamma_{0}$
for some $t_{1},t_{2}\in( 0,1 ) $. Let $\Pi_{i} (
\bbeta,\bgamma) =\Delta_{i} ( \bbeta,\bgamma%
) \ryskus{V}_{i}^{-1} ( \bbeta,\bgamma) $, for $%
1\leq i\leq n$. Then
\[
\calD_{n} \bigl( \bbeta_{n}^{\ast},
\bgamma_{n}^{\ast
} \bigr) =\Psi_{n} (
\bbeta_{0},\bgamma_{0} ) +\Pi_{n,1} \bigl(
\bbeta_{n}^{\ast},\bgamma_{n}^{\ast
} \bigr)
+\Pi_{n,2} \bigl( \bbeta_{n}^{\ast},\bgamma
_{n}^{\ast} \bigr) +\Pi_{n,3}+O \bigl(
n_{\statT}J_{n}^{-p} \bigr),
\]
where $\Pi_{n,1} ( \bbeta_{n}^{\ast},\bgamma%
_{n}^{\ast} ) =-\sum_{i=1}^{n}\underline{\ryskus{D}}_{i}^{\statT}%
\dot{\Pi}_{i} ( \bbeta_{n}^{\ast},\bgamma_{n}^{\ast
} ) \underline{\bvarepsilon}_{i}$,%
\[
\Pi_{n,2} \bigl( \bbeta_{n}^{\ast},
\bgamma_{n}^{\ast
} \bigr) =\sum_{i=1}^{n}
\underline{\ryskus{D}}_{i}^{\statT}\dot{%
\Pi}_{i} \bigl( \bbeta_{n}^{\ast},
\bgamma_{n}^{\ast
} \bigr) \Delta_{i} \bigl(
\bbeta_{n}^{\ast\ast},\bgamma%
_{n}^{\ast\ast}
\bigr) \underline{\ryskus{D}}_{i}\pmatrix{ \bbeta_{n}^{\ast}-
\bbeta_{0}
\cr
\bgamma_{n}^{\ast}-
\bgamma_{0}},
\]
$\Pi_{n,3}=\Psi_{n} ( \bbeta_{0},\bgamma_{0} )
-\Psi_{n} ( \bbeta_{n}^{\ast},\bgamma_{n}^{\ast
} ) $, where $\dot{\Pi}_{i} ( \bbeta_{n}^{\ast}
,\bgamma_{n}^{\ast} ) $ is the first order derivative of $\Pi_{i} (
\bbeta,\bgamma) $ evaluated at $ ( %
\bbeta_{n}^{\ast\statT},\bgamma_{n}^{\ast
\statT} )^{\statT}$, which is a $m_{i}\times
m_{i}\times(d_{1}+d_{2}J_{n})$-dimensional array, $\bbeta%
_{n}^{\ast\ast}$ is between $\bbeta_{n}^{\ast}$ and $%
\bbeta_{0}$, and $\bgamma_{n}^{\ast\ast}$ is between $%
\bgamma_{n}^{\ast}$ and $\bgamma_{0}$. By conditions (C3) and
(C4) and (\ref{EQDD}), for any given vector $\alpha_{n}\in
R^{(d_{1}+d_{2}J_{n})}$ with $\llVert\alpha_{n}\rrVert=1$, there
exists a constant $0<c<\infty$, such that with probability approaching $1$,
$\alpha_{n}^{\statT}\Psi_{n} ( \bbeta_{0}
,\bgamma_{0} ) \alpha_{n}\geq cn_{\statT}\lambda_{n}^{
\mathrm{min}}$. By Theorem~\ref{THMweakconsistency} and (\ref{EQDD}),
$%
\alpha_{n}^{\statT}\Pi_{n,2} ( \bbeta_{n}^{\ast}%
,\bgamma_{n}^{\ast} ) \alpha_{n}=o_{p} ( \lambda_{n}^
{\mathrm{max}} ) $. Since $E \{ \Pi_{n,1} ( %
\bbeta_{n}^{\ast},\bgamma_{n}^{\ast} ) \vert\cal{X%
},\cal{Z} \} =0$, it can be proved by Bernstein's
inequality of~\cite{B98} $\alpha_{n}^{\statT}\Pi_{n,1} (
\bbeta_{n}^{\ast},\bgamma_{n}^{\ast} ) \alpha_{n}=O_{P} \{ (
n_{\statT}\log n_{\statT} )^{1/2} \} $. By condition (C1), $\lambda
_{n}^{\mathrm{max}}
=O ( \tau_{n}^{\max} ) =o(n_{\statT}\lambda_{n}^{\min
}J_{n}^{-1/2})$. Therefore, $\Psi_{n} ( \bbeta_{0}%
,\bgamma_{0} ) $ dominates $\Pi_{n,1} ( \bbeta%
_{n}^{\ast},\bgamma_{n}^{\ast} ) $ and $\Pi_{n,2} (
\bbeta_{n}^{\ast},\bgamma_{n}^{\ast} ) $, and by
Theorem~\ref{THMweakconsistency}, $\Psi_{n} ( \bbeta_{0}%
,\bgamma_{0} ) $ dominates $\Pi_{n,3} ( \bbeta%
_{n}^{\ast},\bgamma_{n}^{\ast} ) $. Thus, from (\ref{EQbetagammahat}), one has
%
\begin{equation}\label{EQbetagamma}
\pmatrix{ \widehat{\bbeta}_{n}-\bbeta_{0}
\cr
\widehat{
\bgamma}_{n}-\bgamma_{0}} =\Psi_{n} (
\bbeta_{0},\bgamma_{0} )^{-1}\ryskus{g}_{n}
( \bbeta_{0},\bgamma_{0} ) \bigl\{ 1+o_{p} ( 1 )
\bigr\}.
\end{equation}
Let $\Delta_{i0}=\Delta_{i} ( \bbeta_{0},\bgamma%
_{0} ) $ and $\ryskus{V}_{i0}=\ryskus{V}_{i} ( \bbeta_{0}%
,\bgamma_{0} ) $. To obtain\vspace*{1pt} the closed-form expression of $%
\widehat{\bbeta}_{n}-\bbeta_{0}$, we need the following
block form of the inverse of $\sum_{i=1}^{n}\underline{\ryskus{D}}%
_{i}^{\statT}\Delta_{i0}\ryskus{V}_{i0}^{-1}\Delta_{i0}
\underline{\ryskus{D}}_{i}$:
%
\begin{eqnarray}\label{EQDVD}
&&\pmatrix{\displaystyle  \sum_{i=1}^{n}\underline{
\ryskus{X}}_{i}^{\statT}\Delta_{i0}%
\ryskus{V}_{i0}^{-1}\Delta_{i0}\underline{
\ryskus{X}}_{i} & \displaystyle \sum_{i=1}^{n}
\underline{\ryskus{X}}_{i}^{\statT}\Delta_{i0}%
\ryskus{V}_{i0}^{-1}\Delta_{i0}\underline{
\ryskus{B}}_{i}
\vspace*{2pt}\cr
\displaystyle \sum_{i=1}^{n}
\underline{\ryskus{B}}_{i}^{\statT}\Delta_{i0}%
\ryskus{V}_{i0}^{-1}\Delta_{i0}\underline{
\ryskus{X}}_{i} & \displaystyle \sum_{i=1}^{n}
\underline{\ryskus{B}}_{i}^{\statT}\Delta_{i0}%
\ryskus{V}_{i0}^{-1}\Delta_{i0}\underline{
\ryskus{B}}_{i}%
}^{-1}
\nonumber\\[-8pt]\\[-8pt]
&&\qquad=\pmatrix{ \ryskus{H}_{\mathbf{XX}} & \ryskus{H}_{\mathbf{XB}}
\cr
\ryskus{H}_{\mathbf{BX}} & \ryskus{H}_{\mathbf{BB}}%
}^{-1}=\pmatrix{ \ryskus{H}^{11} & \ryskus{H}^{12}
\cr
\ryskus{H}^{21} & \ryskus{H}^{22}},\nonumber
\end{eqnarray}
where $\ryskus{H}^{11}= ( \ryskus{H}_{\mathbf{XX}}-\ryskus{H}_{\mathbf
{XB%
}}\ryskus{H}_{\mathbf{BB}}^{-1}\ryskus{H}_{\mathbf{BX}} )^{-1}$, $%
\ryskus{H}^{22}= ( \ryskus{H}_{\mathbf{BB}}-\ryskus{H}_{\mathbf{BX}}%
\ryskus{H}_{\mathbf{XX}}^{-1}\ryskus{H}_{\mathbf{XB}} )^{-1}$, $%
\ryskus{H}^{12}=-\ryskus{H}^{11}\ryskus{H}_{\mathbf{XB}}\ryskus
{H}_{\mathbf{%
BB}}^{-1}$, and $\ryskus{H}^{21}=-\ryskus{H}^{22}\ryskus{H}_{\mathbf
{BX}}%
\ryskus{H}_{\mathbf{XX}}^{-1}$. Consequently, $\widehat{\bbeta}_{n}-%
\bbeta_{0}= ( \widetilde{\bbeta}_{n,e}+\widetilde{%
\bbeta}_{n,\mu} ) \{ 1+o_{p} ( 1 ) \} $,
in which%
\begin{eqnarray*}
\widetilde{\bbeta}_{n,e}&=&\ryskus{H}^{11} \Biggl\{ \sum
_{i=1}^{n}%
\underline{
\ryskus{X}}_{i}^{\statT}\Delta_{i0}
\ryskus{V}_{i0}^{-1}%
\underline{\bvarepsilon
}_{i}-\ryskus{H}_{\mathbf{XB}}\ryskus{H}_{%
\mathbf{BB}}^{-1}
\sum_{i=1}^{n}\underline{
\ryskus{B}}_{i}^{\statT}\Delta_{i0}
\ryskus{V}_{i0}^{-1}\underline{\bvarepsilon}_{i}
\Biggr\},
\\
\widetilde{\bbeta}_{n,\mu}&=&\ryskus{H}^{11}\Biggl[ \sum
_{i=1}^{n}%
\underline{
\ryskus{X}}_{i}^{\statT}\Delta_{i0}
\ryskus{V}%
_{i0}^{-1} \Biggl\{ \mu\Biggl(
\underline{\ryskus{X}}_{i}\bbeta_{0}%
+\sum
_{l=1}^{d_{2}}\theta_{l0}(
\ryskus{Z}_{il}) \Biggr) -\mu( \underline{\ryskus{X}}_{i}
\bbeta_{0}+\underline{%
\ryskus{B}}_{i}
\bgamma_{0} ) \Biggr\}
\\
&&\hspace*{21.5pt}{} -\ryskus{H}_{\mathbf{XB}}\ryskus{H}_{\mathbf{BB}}^{-1}\sum
_{i=1}^{n}%
\underline{
\ryskus{B}}_{i}^{\statT}\Delta_{i0}
\ryskus{V}%
_{i0}^{-1} \Biggl\{ \mu\Biggl(
\underline{\ryskus{X}}_{i}\bbeta_{0}%
+\sum
_{l=1}^{d_{2}}\theta_{l0}(
\ryskus{Z}_{il}) \Biggr)\\
&&\hspace*{203.5pt}{} -\mu( \underline{\ryskus{X}}_{i}
\bbeta_{0}+\underline{ \ryskus{B}}_{i}\bgamma_{0}
) \Biggr\} \Biggr].
\end{eqnarray*}

\begin{lemma}
\label{LEMH11}Under condition \textup{(A4)}, there are constants $%
0<c_{H_{1}}<C_{H_{1}}<\infty$, such that with probability approaching $1$,
for $n_{\statT}$ sufficiently large,\break $c_{H_{1}} ( \lambda_{n}^{\max
}\*n_{\statT} )^{-1}\ryskus{I}_{d_{1}}\leq\ryskus{H}%
^{11}\leq C_{H_{1}} ( \lambda_{n}^{\max}n_{\statT} )^{-1}
\ryskus{I}_{d_{1}}$ with $\ryskus{H}^{11}$in (\ref{EQDVD}).
\end{lemma}
\begin{pf}
The proof of Lemma~\ref{LEMH11} follows the same fashion as the proof of
Lemma A.4 in~\cite{MSW12}, and is hence omitted.
\end{pf}
%
\begin{lemma}
\label{LEMbetahatmu}Under conditions \textup{(A2)} and \textup{(A4)},
$\llVert
\widetilde{%
\bbeta}_{n,\mu}\rrVert=O_{P} \{ ( \lambda_{n}^{\max
}/\lambda_{n}^{\min} )\* J_{n}^{-2p} \} $.
\end{lemma}
\begin{pf}
Let $\Delta\mu(\underline{\bldeta}_{i})=\mu(
\underline{\ryskus{X}}_{i}\bbeta_{0}+\sum_{l=1}^{d_{2}}%
\theta_{l0}(\ryskus{Z}_{il}) ) -\mu( \underline{\ryskus{X}}_{i}%
\bbeta_{0}+\underline{\ryskus{B}}_{i}\bgamma_{0} )
= \break \{ \Delta\mu(\eta_{ij}) \}_{j=1}^{m_{i}}$, then
\begin{eqnarray*}
\widetilde{\bbeta}_{n,\mu}& = &\ryskus{H}^{11} \Biggl[ \sum
_{i=1}^{n}%
\underline{
\ryskus{X}}_{i}^{\statT}\Delta_{i0}
\ryskus{V}%
_{i0}^{-1} \bigl\{ \Delta\mu(
\underline{\bldeta}_{i}) \bigr\} -%
\ryskus{H}_{\mathbf{XB}}
\ryskus{H}_{\mathbf{BB}}^{-1}\sum_{i=1}^{n}
\underline{\ryskus{B}}_{i}^{\statT}\Delta_{i0}
\ryskus{V}_{i0}^{-1} \bigl\{ \Delta\mu(\underline{\bldeta
}_{i}) \bigr\} \Biggr]
\\
& = &\ryskus{H}^{11}\sum_{i=1}^{n}
\underline{\ryskus{X}}_{i}^{\mathrm{%
T}}\Delta_{i0}
\ryskus{V}_{i0}^{-1} \bigl[ \bigl\{ \Delta\mu( \underline{
\bldeta}_{i}) \bigr\} -\widehat{\Pi}_{n} \bigl\{ \Delta\mu(
\underline{\bldeta}_{i}) \bigr\} \bigr] =n_{\statT}
\ryskus{H}^{11}%
\ryskus{W},
\end{eqnarray*}
where $\ryskus{W}= ( W_{1},\ldots,W_{d_{1}} ) $, with
\begin{eqnarray*}
\llvert W_{k}\rrvert& = &n_{\statT}^{-1}\Biggl\llvert
\sum_{i=1}^{n} \bigl( \underline{
\ryskus{X}}_{i}^{ ( k )
} \bigr)^{\statT}
\Delta_{i0}\ryskus{V}_{i0}^{-1} \bigl[ \bigl\{
\Delta\mu(\underline{\bldeta}_{i}) \bigr\} -\widehat{
\Pi}_{n} \bigl\{ \Delta\mu(\underline{\bldeta}_{i})
\bigr\} \bigr] \Biggr\rrvert
\\
& \leq &C\lambda_{n}^{\max}n_{\statT}^{-1}
\sum%
_{i=1}^{n}\sum_{j=1}^{m_{i}}
\bigl\llvert X_{ijk} \bigl\{ \Delta\mu(\eta_{ij}) \bigr\} -
\widehat{\Pi}_{n} \bigl\{ \Delta\mu(\eta_{ij}) \bigr\} \bigr
\rrvert.
\end{eqnarray*}
Following similar reasoning as in the proof of Lemma A.5 in~\cite{MSW12},
it can be proved that $n_{\statT}^{-1}\sum_{i=1}^{n}%
\sum_{j=1}^{m_{i}}\llvert X_{ijk} \{ \Delta\mu(\eta_{ij}) \} -%
\widehat{\Pi}_{n} \{ \Delta\mu(\eta_{ij}) \} \rrvert
=O_{P} ( J_{n}^{-2p} ) $. Therefore, $\llvert W_{k}\rrvert
=O_{P} ( \lambda_{n}^{\max}J_{n}^{-2p} ) $. By the above result
and Lemma~\ref{LEMH11}, one has $\llVert\widetilde{\bbeta}%
_{n,\mu}\rrVert=O_{P} \{ ( \lambda_{n}^{\min} )^{-1}\lambda
_{n}^{\max}J_{n}^{-2p} \} $.
\end{pf}
%
\begin{lemma}
\label{LEMbetahateasymp}Under conditions \textup{(A2)--(A4)}, as $n_{\statT}%
\rightarrow\infty$, $\widetilde{\Xi}_{n}^{-1/2} ( \widetilde{%
\bbeta}_{n,e} ) \longrightarrow N ( 0,\ryskus{I}_{d_{1}} ) $,
where $\widetilde{\Xi}_{n}$ is defined in (\ref{DEFgammatilda}).
\end{lemma}
\begin{pf}
Lemma~\ref{LEMbetahateasymp} can be proved by using the
Linderberg--Feller CLT
and similar techniques for the proofs of Lemmas A.6 and A.7 in~\cite{MSW12}.
\end{pf}
%
\begin{lemma}
\label{LEMbetahatcon}Under conditions \textup{(A2)} and \textup{(A4)},
there exist
constants $%
0<c_{\Xi}\leq C_{\Xi}<\infty$, such that
\[
c_{\Xi}n_{\statT}^{-1} \bigl( \lambda_{n}^{\max}
\bigr)^{-1}\tau_{n}^{\min}\ryskus{I}_{d_{1}}
\leq\widetilde{\Xi}_{n}\leq C_{\Xi}n_{\statT%
}^{-1}
\tau_{n}^{\max} \bigl( \lambda_{n}^{\min}
\bigr)^{-1}\ryskus{I}%
_{d_{1}}
\]
and $\llVert\widetilde{\bbeta}_{n,e}\rrVert=O_{p} \{ n_{%
\statT}^{-1/2} ( \tau_{n}^{\max} )^{1/2} ( \lambda_{n}^{\min} )^{-1/2}
\} $.
\end{lemma}
\begin{pf}
For any vector $a\in R^{d_{1}}$ with $\llVert a\rrVert=1$, one has
\begin{eqnarray*}
a^{\statT}\widetilde{\Xi}_{n}a&\leq&\tau_{n}^{\max}a^{\statT}
\Biggl\{ E \Biggl( \sum_{i=1}^{n}
\widetilde{\underline{\ryskus{X}}}_{i}^{\statT%
}
\Delta_{i0}\ryskus{V}_{i0}^{-1}
\Delta_{i0}\widetilde{\underline{\ryskus{X}}}%
_{i}
\Biggr) \Biggr\}^{-1}a\leq C_{\Xi}n_{\statT}^{-1}
\tau_{n}^{\max
} \bigl( \lambda_{n}^{\min}
\bigr)^{-1},
\\
a^{\statT}\widetilde{\Xi}_{n}a&\geq&\Biggl\{ E \Biggl( \sum
_{i=1}^{n}%
\widetilde{
\underline{\ryskus{X}}}_{i}^{\statT}\Delta_{i0}
\ryskus{V}%
_{i0}^{-1}\Delta_{i0}
\widetilde{\underline{\ryskus{X}}}_{i} \Biggr) \Biggr\}^{-1}
\tau_{n}^{\min}\geq c_{\Xi}n_{\statT}^{-1}
\bigl( \lambda_{n}^{\max} \bigr)^{-1}
\tau_{n}^{\min},
\end{eqnarray*}
and the second result in Lemma~\ref{LEMbetahatcon} follows from Chebyshev's
inequality.
\end{pf}
\begin{pf*}{Proof of Theorem \protect\ref{THMbetahatnormality}}
By Lemmas~\ref{LEMbetahatmu} and~\ref{LEMbetahatcon}, for any
vector $a\in R^{d_{1}}$ with $\llVert a\rrVert=1$, one has%
\begin{eqnarray*}
a^{\statT}\widetilde{\Xi}_{n}^{-1/2}\widetilde{
\bbeta}_{n,\mu
}a&\leq& c_{\Xi}^{-1/2}n_{\statT}^{1/2}
\bigl( \lambda_{n}^{\max} \bigr)^{1/2} \bigl(
\tau_{n}^{\min} \bigr)^{-1/2}O_{P} \bigl\{
\bigl( \lambda_{n}^{\min} \bigr)^{-1}
\lambda_{n}^{\max}J_{n}^{-2p} \bigr\}
\\
&=&O_{P} \bigl\{ n_{\statT}^{1/2}J_{n}^{-2p}
\bigl( \lambda_{n}^{\max
} \bigr)^{3/2} \bigl(
\lambda_{n}^{\min} \bigr)^{-1} \bigl(
\tau_{n}^{\min} \bigr)^{-1/2} \bigr\}
=o_{p} ( 1 ).
\end{eqnarray*}
Therefore, Theorem~\ref{THMbetahatnormality} follows from Lemma
\ref{LEMbetahateasymp}, the above result and Slutsky's theorem.
\end{pf*}

\subsection{\texorpdfstring{Proof of Theorem \protect\ref{THMthetahatnormality}}{Proof of Theorem 3}}

Following the same reasoning as deriving (\ref{EQbetagamma}), it can be
proved that
%
\begin{eqnarray}\label{EQgammatildas}
\widehat{\bgamma}_{n,l}^{\calS} ( \bbeta_{0}%
,\btheta_{-l0} ) -\bgamma_{l,0}^{\calS} &=&
\Psi_{n,l}^{\ast} \bigl( \bgamma_{l,0}^{\calS}
\bigr)^{-1}%
\ryskus{g}_{n,l}^{\ast} (
\bgamma_{l,0} ) \bigl( 1+o_{p} ( 1 ) \bigr)
\nonumber\\[-8pt]\\[-8pt]
&=& \bigl( \widetilde{\bgamma}_{n,e,l}^{\calS}+
\widetilde{%
\bgamma}_{n,\mu,l}^{\calS} \bigr) \bigl(
1+o_{p} ( 1 ) \bigr),\nonumber
\end{eqnarray}
where
\begin{eqnarray*}
\widetilde{\bgamma}_{n,e,l}^{\calS}&=&\widetilde{%
\bgamma}_{n,e,l}^{\calS} ( \bbeta_{0},\btheta_{-l0} )\\
&=&\Psi_{n,l}^{\ast} \bigl( \bgamma_{l,0}^{\calS}
\bigr)^{-1}\sum_{i=1}^{n} \bigl(
\ryskus{B}_{i\cdot l}^{\calS%
} \bigr)^{\statT}
\Delta_{i} \bigl( \bbeta_{0},\btheta_{-l0},
\bgamma_{l,0}^{\calS} \bigr) \ryskus{V}_{i}^{-1}
\bigl( \bbeta_{0},\btheta_{-l0},\bgamma_{l,0}^{\calS}
\bigr) \underline{\bvarepsilon}_{i},
\\
\widetilde{\bgamma}_{n,\mu,l}^{\calS}&=&\widetilde{%
\bgamma}_{n,\mu,l}^{\calS} ( \bbeta_{0},\btheta
_{-l0} ) = \bigl( \widetilde{\gamma}_{n,\mu,sl}^{\calS}
\bigr)_{s=1}^{J_{n}^{\calS}}
\\
&=&\Psi_{n,l}^{\ast} \bigl( \bgamma_{l,0}^{\calS}
\bigr)^{-1}\sum_{i=1}^{n} \bigl(
\ryskus{B}_{i\cdot l}^{\calS%
} \bigr)^{\statT}
\Delta_{i} \bigl( \bbeta_{0},\btheta%
_{-l0},
\bgamma_{l,0}^{\calS} \bigr) \ryskus{V}_{i}^{-1}
\bigl( \bbeta_{0},\btheta_{-l0},\bgamma_{l,0}^{\calS}
\bigr)
\\
&&\hspace*{70.2pt}{}\times\biggl\{ \mu\biggl( \underline{\ryskus{X}}_{i}\bbeta_{0}+
\sum_{l^{\prime}\neq l}\theta_{l^{\prime}0}(
\ryskus{Z}_{il^{\prime
}})+\theta_{l0}(\ryskus{Z}_{il})
\biggr) \\
&&\hspace*{87.4pt}{}-\mu\biggl( \underline{\ryskus{X}}%
_{i}
\bbeta_{0}+\sum_{l^{\prime}\neq l}\theta_{l^{\prime}0}(%
\ryskus{Z}_{il^{\prime}})+\ryskus{B}_{i\cdot l}^{\calS}
\bgamma_{l,0}^{\calS} \biggr) \biggr\}.
\end{eqnarray*}
By the decomposition in (\ref{EQgammatildas}),%
\begin{eqnarray*}
\widehat{\theta}_{n,l}^{\calS}(z_{l},
\bbeta_{0} ,\btheta_{-l0})-\theta_{l0}^{\ast}(z_{l})&=&
\ryskus{B}_{l}^{\calS%
} ( z_{l} )^{\statT}
\widetilde{\bgamma}_{n,e,l}^{%
\calS} \bigl( 1+o_{p} (
1 ) \bigr),
\\
\theta_{l0}^{\ast}(z_{l})-\theta_{l0}(z_{l})
&=& \bigl\{ \ryskus{B}_{l}^{%
\calS} ( z_{l}
)^{\statT}\widetilde{\bgamma}%
_{n,\mu,l}^{\calS}+
\ryskus{B}_{l}^{\calS} ( z_{l} )^{%
\statT}
\bgamma_{l,0}^{\calS}-\theta_{l0}(z_{l})
\bigr\}\\
&&{}\times \bigl( 1+o_{p} ( 1 ) \bigr).
\end{eqnarray*}
It can be proved by the Linderberg--Feller CLT that as
$n_{\statT}\rightarrow\infty$,
\[
\bigl( \ryskus{B}_{l}^{\calS} ( z_{l}
)^{\statT}\Xi_{n,l}^{\ast}\ryskus{B}_{l}^{\calS}
( z_{l} ) \bigr)^{-1/2} \bigl( \ryskus{B}_{l}^{\calS}
( z_{l} )^{\statT}%
\widetilde{\bgamma
}_{n,e,l}^{\calS} \bigr) \longrightarrow N ( 0,1 ).
\]
Following similar reasoning as in the proofs in Lemma \ref
{LEMgammahatnm}, it
can be proved
\[
\sup_{1\leq s\leq J_{n}^{\calS}}\bigl\llvert\widetilde{\gamma}%
_{n,\mu,sl}^{\calS}
\bigr\rrvert=O_{P} \bigl\{ \bigl( \lambda_{n}^{\min}
\bigr)^{-1}\lambda_{n}^{\max} \bigl(
J_{n}^{\calS%
} \bigr)^{-p-1/2} \bigr\}
\]
and
\[
\bigl\llVert\widetilde{\bgamma}_{n,\varepsilon,l}^{%
\calS}
\bigr\rrVert_{\infty}=O_{P} \bigl\{ ( \log n_{\statT%
}/n_{\statT}
)^{1/2} \bigl( \tau_{n}^{\max} \bigr)^{1/2}
\bigl( \lambda_{n}^{\min} \bigr)^{-1/2} \bigr\}.
\]
By B-spline properties, $\sup_{z_{l}\in[ 0,1 ] }\llvert
\ryskus{B}_{l}^{\calS} ( z_{l} )^{\statT}\widetilde{%
\bgamma}_{n,\mu,l}^{\calS}\rrvert=O_{P} \{ (
\lambda_{n}^{\max}/\lambda_{n}^{\min} ) ( J_{n}^{\calS%
} )^{-p} \} $, and $\sup_{z_{l}\in[ 0,1 ] }\llvert
\ryskus{B}_{l}^{\calS} ( z_{l} )^{\statT}\widetilde{%
\bgamma}_{n,\varepsilon,l}^{\calS}\rrvert=O_{P} \{
\sqrt{ ( \log n_{\statT} ) J_{n}^{\calS}/n_{\statT}}%
( \tau_{n}^{\max}/\lambda_{n}^{\min} )^{1/2} \} $, so%
\begin{eqnarray*}
\sup_{z_{l}\in[ 0,1 ] }\bigl\llvert\theta_{l0}^{\ast
}(z_{l})-
\theta_{l0}(z_{l})\bigr\rrvert&\leq&\sup_{z_{l}\in[
0,1 ] }
\bigl\llvert\ryskus{B}_{l}^{\calS} ( z_{l}
)^{%
\statT}\widetilde{\bgamma}_{n,\mu,l}^{\calS}\bigr\rrvert
\\
&&{}+\sup_{z_{l}\in[ 0,1 ] }\bigl\llvert\ryskus{B}_{l}^{%
\calS} (
z_{l} )^{\statT}\bgamma_{l,0}^{%
\calS}-
\theta_{l0}(z_{l})\bigr\rrvert\\
&=&O_{P} \bigl\{
\bigl( \lambda_{n}^{\min} \bigr)^{-1}
\lambda_{n}^{\max} \bigl( J_{n}^{\calS%
}
\bigr)^{-p} \bigr\},
\end{eqnarray*}
$\sup_{z_{l}\in[ 0,1 ] }\llvert\widehat{\theta}_{n,l}^{%
\calS}(z_{l},\bbeta_{0},\btheta_{-l0})-\theta_{l0}^{\ast
}(z_{l})\rrvert=O_{P} \{ \sqrt{ ( \log n_{\statT%
} ) J_{n}^{\calS}/n_{\statT}} ( \tau_{n}^{\max}/\lambda_{n}^{\min}
)^{1/2} \} $.

\subsection{\texorpdfstring{Proof of Theorem \protect\ref{THMthetahatconvergence}}{Proof of Theorem 4}}

\begin{lemma}
\label{LEMgammahatnm}Under conditions \textup{(A2)--(A4)},
\begin{eqnarray*}
\llVert\widehat{\bgamma}_{n}-\bgamma_{0}\rrVert
&=&O_{P} \bigl\{ J_{n}^{1/2}n_{\statT}^{-1/2}
\bigl( \tau_{n}^{\max
}/\lambda_{n}^{\min}
\bigr)^{1/2}+ \bigl( \lambda_{n}^{\max}/
\lambda_{n}^{\min} \bigr) J_{n}^{-p} \bigr
\},
\\[-2pt]
\llVert\widehat{\bgamma}_{n}-\bgamma_{0}\rrVert_{\infty}&=&O_{P} \bigl\{
( \log n_{\statT}/n_{\statT} )^{1/2} \bigl( \tau_{n}^{\max}/\lambda
_{n}^{\min} \bigr)^{1/2}+ \bigl(
\lambda_{n}^{\max}/\lambda_{n}^{\min} \bigr) J_{n}^{-p-1/2} \bigr\}.
\end{eqnarray*}
\end{lemma}
\begin{pf}
From (\ref{EQbetagamma}) and (\ref{EQDVD}), one obtains $\widehat{%
\bgamma}_{n}-\bgamma_{0}= ( \widetilde{\bgamma}_{n,e}+%
\widetilde{\bgamma}_{n,\mu} ) ( 1+o_{p} ( 1 )
) $, where%
\begin{eqnarray*}
\widetilde{\bgamma}_{n,e}&=&\ryskus{H}^{22} \Biggl\{ \sum
_{i=1}^{n}\underline{\ryskus{B}}_{i}^{\statT}
\Delta_{i0}%
\ryskus{V}_{i0}^{-1}
\underline{\bvarepsilon}_{i}-\ryskus{H}_{%
\mathbf{BX}}
\ryskus{H}_{\mathbf{XX}}^{-1}\sum_{i=1}^{n}
\underline{\ryskus{X}}_{i}^{\statT}\Delta_{i0}
\ryskus{V}_{i0}^{-1} \underline{\bvarepsilon}_{i}
\Biggr\},
\\[-2pt]
\widetilde{\bgamma}_{n,\mu}&=&\ryskus{H}^{22}\Biggl[ \sum
_{i=1}^{n}\underline{\ryskus{B}}_{i}^{\statT}
\Delta_{i0}%
\ryskus{V}_{i0}^{-1} \Biggl
\{ \mu\Biggl( \underline{\ryskus{X}}_{i}%
\bbeta_{0}+\sum_{l=1}^{d_{2}}
\theta_{l0}(\ryskus{Z}%
_{il}) \Biggr) -\mu(
\underline{\ryskus{X}}_{i}\bbeta_{0}+%
\underline{
\ryskus{B}}_{i}\bgamma_{0} ) \Biggr\}
\\[-2pt]
&&\hspace*{22pt}{} -\ryskus{H}_{\mathbf{BX}}\ryskus{H}_{\mathbf{XX}}^{-1}\sum
_{i=1}^{n}%
\underline{
\ryskus{X}}_{i}^{\statT}\Delta_{i0}
\ryskus{V}%
_{i0}^{-1} \Biggl\{ \mu\Biggl(
\underline{\ryskus{X}}_{i}\bbeta_{0}%
+\sum
_{l=1}^{d_{2}}\theta_{l0}(
\ryskus{Z}_{il}) \Biggr) \\[-2pt]
&&\hspace*{203pt}{}-\mu( \underline{\ryskus{X}}_{i}
\bbeta_{0}+\underline{ \ryskus{B}}_{i}\bgamma_{0}
) \Biggr\} \Biggr].
\end{eqnarray*}
It can be proved that there exist constants
$0<c_{H_{2}}<C_{H_{2}}<\infty$,
such that with probability approaching $1$, for $n_{\statT}$
sufficiently large,
\[
c_{H_{2}} \bigl( \lambda_{n}^{\max}
\bigr)^{-1}n_{\statT}^{-1}%
\ryskus{I}_{d_{1}}\leq\ryskus{H}^{22}\leq C_{H_{2}}
\bigl( \lambda_{n}^{\min} \bigr)^{-1}n_{\statT}^{-1}
\ryskus{I}_{d_{1}}.
\]
Letting $\widehat{\Pi}_{n,\ryskus{X}}$ be the projection on $ \{
\underline{\ryskus{X}}_{i} \}_{i=1}^{n}$ to the empirical inner
product,
\[
\widetilde{\bgamma}_{n,\mu}=\ryskus{H}^{22}\sum
_{i=1}^{n} \underline{\ryskus{B}}_{i}^{\statT}
\Delta_{i0}\ryskus{V}_{i0}^{-1} \bigl[ \bigl\{
\Delta\mu(\underline{\bldeta}_{i}) \bigr\} -\widehat{\Pi
}%
_{n,\ryskus{X}} \bigl\{ \Delta\mu(\underline{\bldeta
}_{i}) \bigr\} %
\bigr] =n_{\statT}
\ryskus{H}^{22}\ryskus{W},
\]
where $\ryskus{W}= ( W_{1},\ldots,W_{J_{n}d_{2}} ) $, with
\[
W_{s,l}=n_{\statT}^{-1}\sum
_{i=1}^{n} \bigl( \ryskus{B}%
_{i}^{ ( s,l ) }
\bigr)^{\statT}\Delta_{i0}\ryskus{V}%
_{i0}^{-1}
\bigl[ \bigl\{ \Delta\mu(\underline{\bldeta}_{i}) \bigr\} -
\widehat{\Pi}_{n,\ryskus{X}} \bigl\{ \Delta\mu(\underline{\bldeta
}%
_{i}) \bigr\} \bigr],
\]
$\ryskus{B}_{i}^{ ( s,l ) }= [ \{ B_{s,l} (
Z_{i1l} ),\ldots,B_{s,l} ( Z_{im_{i}l} ) \}^{\statT}%
] $. The Cauchy--Schwarz inequality implies
\begin{eqnarray*}
\llvert W_{s,l}\rrvert& \leq &C\lambda_{n}^{\max
}n_{\statT}^{-1}
\sum_{i=1}^{n}\sum
_{j=1}^{m_{i}}\bigl\llvert B_{s,l} (
Z_{ijl} ) \bigl\{ \Delta\mu(\eta_{ij}) \bigr\} -\widehat{
\Pi}_{n,%
\ryskus{X}} \bigl\{ \Delta\mu(\eta_{ij}) \bigr\} \bigr
\rrvert
\\[-2pt]
& \leq &C\lambda_{n}^{\max}\llVert B_{s,l}
\rrVert_{n_{\statT}}\bigl\llVert\Delta\mu-\widehat{\Pi}_{n,\ryskus
{X}} (
\Delta\mu) \bigr\rrVert_{n_{\statT}}=O_{P} \bigl(
\lambda_{n}^{\max}J_{n}^{-p-1/2} \bigr),
\end{eqnarray*}
thus, $\llVert\widetilde{\bgamma}_{n,\mu}\rrVert
=O_{P} \{ ( \lambda_{n}^{\max}/\lambda_{n}^{\min} )
J_{n}^{-p} \} $, $\llVert\widetilde{\bgamma}_{n,\mu
}\rrVert_{\infty}=O_{P} \{ ( \lambda_{n}^{\max}/\lambda_{n}^{\min} )
J_{n}^{-p-1/2} \} $. For any $\omega\in\cal{R} %
^{J_{n}d_{2}}$ with $\llVert\omega\rrVert=1$, it can be proved\vadjust{\goodbreak}
that $\operatorname{Var}( \omega^{\statT}\widetilde{\bgamma}%
_{n,e}\vert\cal{X},\cal{Z} ) \leq O_{P} \{ n_{%
\statT}^{-1} ( \tau_{n}^{\max}/\lambda_{n}^{\min} )
\} $, thus, $\omega^{\statT}\widetilde{\bgamma}%
_{n,e}=O_{P} \{ n_{\statT}^{-1/2} ( \tau_{n}^{\max}/\lambda_{n}^{\min
} )^{1/2} \} $. Therefore, $\llVert\widetilde{%
\bgamma}_{n,e}\rrVert\leq J_{n}^{1/2}\llvert\omega^{%
\statT}\widetilde{\bgamma}_{n,e}\rrvert=O_{P} \{
J_{n}^{1/2}n_{\statT}^{-1/2} ( \tau_{n}^{\max}/\lambda_{n}^{\min
} )^{1/2} \} $, and by Bernstein's inequality of~\cite{B98} that $%
\llVert\widetilde{\bgamma}_{n,e}\rrVert_{\infty
}=O_{P} \{ ( \log n_{\statT}/n_{\statT} )^{1/2} (
\tau_{n}^{\max}/\lambda_{n}^{\min} )^{1/2} \}$.
\end{pf}
%
\begin{lemma}
\label{LEMgammahats}Under conditions \textup{(A2)--(A4)},
\[
\bigl\llVert\widehat{\bgamma}_{n,l}^{\calS\calS}-
\widehat{%
\bgamma}_{n,l}^{\mathrm{OR}}\bigr
\rrVert_{\infty}=O_{p} \Bigl\{ \bigl( \lambda_{n}^{\max}/
\lambda_{n}^{\min} \bigr)^{2} \Bigl( \sqrt{
\log n_{%
\statT}/ \bigl( J_{n}^{\calS}n_{\statT}
\bigr) }%
+ \bigl( J_{n}^{\calS} \bigr)^{-1/2}J_{n}^{-p}
\Bigr) \Bigr\}.
\]
\end{lemma}
\begin{pf}
Let $\widetilde{\btheta}_{-l0}= \{ \widetilde{\theta}%
_{l^{\prime}0} ( \cdot),l^{\prime}\neq l \} $, where $%
\widetilde{\theta}_{l^{\prime}0} ( \cdot) $ is defined in (\ref{DEFglzl}). Let $\widehat{\bgamma}_{n,-l}= ( \widehat{\gamma}%
_{n,sl^{\prime}}\dvtx 1\leq s\leq J_{n},l^{\prime}\neq l )^{\statT}$
and $\bgamma_{-l0}= ( \gamma_{sl^{\prime},0}\dvtx 1\leq s\leq
J_{n},l^{\prime}\neq l )^{\statT}$. By the Taylor expansion, $\ryskus{%
g}_{n,l}^{\calS} ( \widehat{\bgamma}_{n,l}^{\mathrm{OR}},%
\widehat{\bbeta}_{n},\widehat{\btheta}%
_{n,-l} ) -\ryskus{g}_{n,l}^{\calS}
( \widehat{\bgamma}_{n,l}^{\mathrm{OR}},\widehat{\bbeta
}_{n},\widetilde{%
\btheta}_{-l0} ) =
\{ \partial\ryskus{g}_{n,l}^{\calS} ( \widehat{%
\bgamma}_{n,l}^{\mathrm{OR}},\break\widehat{\bbeta}_{n},%
\widetilde{\btheta}_{-l} ) /\partial\widetilde{\bgamma}_{-l}^{\statT
} \} ( \widehat{\bgamma}_{n,-l}-%
\bgamma_{-l0} ) $, where $\widetilde{\bgamma}_{-l}=t%
\bgamma_{-l0}+ ( 1-t ) \widehat{\bgamma}_{n,-l}$ for $%
t\in( 0,1 ) $. Let $\widehat{\Delta}_{i}=\Delta_{i} (
\widehat{\bbeta}_{n},\widetilde{\btheta}_{-l},%
\widehat{\bgamma}_{n,l}^{\mathrm{OR}} ) $, $\widehat{\ryskus{V}}%
_{i}=\ryskus{V}_{i} ( \widehat{\bbeta}_{n},\widetilde{%
\btheta}_{-l},\widehat{\bgamma}_{n,l}^{\mathrm{OR}} )$,
$\widetilde{\underline{\bvarepsilon}_{i}}=\underline{%
\bvarepsilon}_{i}-\widehat{\Pi}_{n,\ryskus{X}} ( \underline{
\bvarepsilon}_{i} ) $,\break $\widetilde{\Delta\mu}(\underline{
\bldeta}_{i})=\Delta\mu(\underline{\bldeta}_{i})-\widehat{\Pi}_{n,%
\ryskus{X}} \{ \Delta\mu(\underline{\bldeta}_{i}) \} $, $%
\ryskus{B}_{ij,-l}= \{ ( \ryskus{B}_{ijl^{\prime
}}^{\statT},l^{\prime}\neq l )^{\statT} \}_{ (
d_{2}-1 ) J_{n}\times1}$,\break $\ryskus{B}_{i,-l}= \{ (
\ryskus{B}_{i1,-l},\ldots,\ryskus{B}_{im_{i},-l} )^{\statT} \}
_{m_{i}\times(
d_{2}-1 ) J_{n}}$. Thus, by (\ref{EQgnlstar}) and the proofs for Lemma
\ref{LEMgammahatnm}, with probability approaching $1$, there are
constants $%
0<C_{1,}C_{2}<\infty$ such that%
\begin{eqnarray*}
&&
\bigl\llVert\ryskus{g}_{n,l}^{\calS} \bigl( \widehat{\bgamma
}%
_{n,l}^{\mathrm{OR}},\widehat{\bbeta}_{n},
\widehat{%
\btheta}_{n,-l} \bigr) -\ryskus{g}_{n,l}^{\calS}
\bigl( \widehat{%
\bgamma}_{n,l}^{\mathrm{OR}},\widehat{
\bbeta}_{n},%
\widetilde{\btheta}_{-l0} \bigr)
\bigr\rrVert_{\infty}\\
&&\qquad\leq C_{1} \bigl( \lambda_{n}^{\min}
\bigr)^{-1}n_{\statT}^{-1}
\\
&&\qquad\quad{}\times\Biggl\llVert\Biggl( \sum_{i=1}^{n}
\bigl( \ryskus{B}_{i\cdot l}^{%
\calS} \bigr)^{\statT}\widehat{
\Delta}_{i}\widehat{\ryskus{V}}%
_{i}^{-1}
\ryskus{B}_{i,-l} \Biggr) \Biggl\{ \sum_{i=1}^{n}
\ryskus{B}%
_{i,-l}^{\statT}\Delta_{i0}
\ryskus{V}_{i0}^{-1} \bigl( \widetilde{\underline{
\bvarepsilon}}_{i}+\widetilde{\Delta\mu}(\underline{%
\bldeta}_{i}) \bigr) \Biggr\} \Biggr\rrVert_{\infty}
\\
&&\qquad\leq C_{2} \bigl( \lambda_{n}^{\min}
\bigr)^{-1} \bigl( \llVert\zeta_{1}\rrVert_{\infty}+
\llVert\zeta_{2}\rrVert_{\infty
} \bigr),
\end{eqnarray*}
where $\zeta_{1}=n_{\statT}^{-1} \{
\sum_{i=1}^{n} ( \ryskus{B}_{i\cdot l}^{\calS} )^{%
\statT}\Delta_{i0}\ryskus{V}_{i0}^{-1}\ryskus{B}_{i,-l} \} \{
\sum_{i=1}^{n}\ryskus{B}_{i,-l}^{\statT}\Delta_{i0}\ryskus{V}%
_{i0}^{-1} ( \widetilde{\Delta\mu}(\underline{\bldeta}%
_{i}) ) \} $,\break
$\zeta_{2}=n_{\statT}^{-1} \{ \sum_{i=1}^{n} (
\ryskus{B}_{i\cdot l}^{\calS} )^{\statT}\Delta_{i0}\ryskus{V%
}_{i0}^{-1}\ryskus{B}_{i,-l} \} ( \sum_{i=1}^{n}\ryskus{B}%
_{i,-l}^{\statT}\Delta_{i0}\ryskus{V}_{i0}^{-1}\widetilde{
\underline{\bvarepsilon}}_{i} ) $, and then\break $\llVert\zeta_{1}\rrVert
_{\infty}
\leq( \lambda_{n}^{\max} )^{2}\llVert\zeta_{3}\rrVert_{\infty}O (
J_{n}^{-p} ) $, where $\zeta_{3}=\Delta_{1}+\Delta_{2}+\Delta_{3}$,
$\Delta_{1}= ( \delta_{1s} )_{s=1}^{J_{n}^{%
\calS}}$, $\Delta_{2}= ( \delta_{2s} )_{s=1}^{J_{n}^{%
\calS}}$ and $\Delta_{3}= ( \delta_{3s} )_{s=1}^{J_{n}^{%
\calS}}$ with $\delta_{1s}=n_{\statT}^{-1}\sum%
_{i=1}^{n}\delta_{1s,i}$, $\delta_{2s}=n_{\statT%
}^{-1}\sum_{i=1}^{n}\delta_{2s,i}$ and $\delta_{3s}=n_{\statT}^{-1}\sum
_{i=1}^{n}\delta_{3s,i}$,
\begin{eqnarray*}
\delta_{1s,i} &=&\sum_{j=1}^{m_{i}}
\sum_{l^{\prime
}=1,l^{\prime}\neq l}^{d_{2}}\sum
_{s^{\prime
}=1}^{J_{n}}\bigl\llvert B_{s,l}^{\calS}
( Z_{ijl} ) \bigr\rrvert\bigl\llvert B_{s^{\prime},l^{\prime}} (
Z_{ijl^{\prime
}} ) \bigr\rrvert^{2},
\\
\delta_{2s,i} &=&\sum_{j=1}^{m_{i}}
\sum_{j\prime\neq
j}\sum_{l^{\prime}\neq l}\sum
_{s^{\prime
}=1}^{J_{n}}\bigl\llvert B_{s,l}^{\calS}
( Z_{ijl} ) \bigr\rrvert\bigl\llvert B_{s^{\prime},l^{\prime}} (
Z_{ijl^{\prime
}} ) \bigr\rrvert\bigl\llvert B_{s^{\prime},l^{\prime}} (
Z_{ij^{\prime}l^{\prime}} ) \bigr\rrvert,
\\
\delta_{3s,i} &=&\sum_{j=1}^{m_{i}}
\sum_{i^{\prime}\neq
i}\sum_{j^{\prime}}\sum
_{l^{\prime}\neq
l}\sum_{s^{\prime}=1}^{J_{n}}
\bigl\llvert B_{s,l}^{\calS} ( Z_{ijl} ) \bigr\rrvert
\bigl\llvert B_{s^{\prime},l^{\prime}} ( Z_{ijl^{\prime}} ) \bigr
\rrvert\bigl\llvert
B_{s^{\prime},l^{\prime
}} ( Z_{i^{\prime}j^{\prime}l^{\prime}} ) \bigr\rrvert.
\end{eqnarray*}
Let $\delta_{1s,i}^{\ast}=\delta_{1s,i}-E ( \delta_{1s,i} ) $.
It can be proved by B-spline properties that $E ( \delta_{1s,i} )
\asymp m_{i}J_{n}/\sqrt{J_{n}^{\calS}}$, $E ( \delta_{1s,i}^{\ast
} ) =0$, $E ( \delta_{1s,i}^{\ast} )^{2}\asymp
m_{i}J_{n}^{2}+m_{i}^{2}J_{n}^{2} ( J_{n}^{\calS} )^{-1}$,
and $E ( \llvert\delta_{1s,i}^{\ast}\rrvert^{k} ) \leq
C \{ m_{i}J_{n}^{k} ( J_{n}^{\calS} )^{k/2-1}+m_{i}^{2}J_{n}^{k} (
J_{n}^{\calS} )^{k/2-2} \} $ for $k\geq3$ and some constant $C>0$.
Thus, $E (
\llvert\delta_{1s,i}^{\ast}\rrvert^{k} ) \leq(
C^{\prime} ( J_{n}^{\calS} )^{1/2}J_{n} )^{k-2}k!E ( \delta
_{1s,ijl^{\prime}s^{\prime}}^{2} ) $ with $%
C^{\prime}=C^{1/ ( k-2 ) }$. By\vspace*{1pt} Bernstein's inequality in~\cite{B98},
\[
P \Biggl( \Biggl\llvert\sum_{i=1}^{n}
\delta_{1s,i}\Biggr\rrvert\geq t \Biggr) \leq2\exp\biggl\{ -
\frac{t^{2}}{4\sum_{i=1}^{n}E (
\delta_{1s,i}^{\ast} )^{2}+2C^{\prime} ( J_{n}^{\calS%
} )^{1/2}J_{n}t} \biggr\}.
\]
Let $t=c \{ \{ n_{\statT}J_{n}^{2}+ (
\sum_{i=1}^{n}m_{i}^{2} ) J_{n}^{2} ( J_{n}^{\calS%
} )^{-1} \} \log n_{\statT} \}^{1/2}$ for a large
constant $0<c<\infty$. There is a constant $0<c^{\prime}<\infty$
such that $E ( \delta_{1s,i}^{\ast} )^{2}\leq c^{\prime
} \{ m_{i}J_{n}^{2}+m_{i}^{2}J_{n}^{2} (
J_{n}^{\calS} )^{-1} \} $. For
$J_{n}^{\calS}=O ( ( \log n_{\statT} )^{-1}n_{\statT}^{1/2}m_{ ( n )
}^{1/2} ) $, one has
$P ( \llvert\sum_{i=1}^{n}\delta_{1s,i}\rrvert\geq t ) \leq
2n_{\statT}^{-c^{2}/ (
4c^{\prime} ) }$. By the Borel--Cantelli lemma,
\[
\max_{1\leq s\leq J_{n}^{\calS}}\bigl\llvert\delta_{1s}-E (
\delta_{1s} ) \bigr\rrvert=O_{\mathrm{a.s.}} \bigl\{
n_{\statT}^{-1/2}J_{n} \bigl( 1+m_{ ( n ) }/J_{n}^{\calS}
\bigr)^{1/2} ( \log n_{\statT} )^{1/2} \bigr\}.
\]
Since $E ( \delta_{1s} ) \asymp J_{n}/\sqrt{J_{n}^{\calS}}$,
one has $\llVert\Delta_{1}\rrVert_{\infty}=O_{\mathrm{a.s.}} ( J_{n}/%
\sqrt{J_{n}^{\calS}} ) $. Since $E ( \delta_{2s} )
\asymp n_{\statT}^{-1} (
\sum_{i=1}^{n}m_{i}^{2} ) /\sqrt{J_{n}^{\calS}}$ and $%
E ( \delta_{3s} ) \asymp n_{\statT}/\sqrt{J_{n}^{%
\calS}}$, similarly it can be proved that $\llVert\Delta_{2}\rrVert
_{\infty}=O_{\mathrm{a.s.}} ( m_{ ( n ) }/\sqrt{J_{n}^{%
\calS}} ) $ and $\llVert\Delta_{3}\rrVert_{\infty
}=O_{\mathrm{a.s.}} ( n_{\statT}/\sqrt{J_{n}^{\calS}} ) $.
Therefore, $\llVert\zeta_{1}\rrVert_{\infty}=O_{\mathrm{a.s.}} \{
( \lambda_{n}^{\max} )^{2}n_{\statT} ( J_{n}^{%
\calS} )^{-1/2}J_{n}^{-p} \} $. Following similar
reasoning, by Bernstein's inequality one can prove $\llVert\zeta
_{2}\rrVert_{\infty}=O_{\mathrm{a.s.}} ( ( \lambda_{n}^{\max} )^{2}\sqrt
{n_{\statT}\log n_{\statT}/J_{n}^{\calS}%
} ) $. Thus,
\[
\bigl\llVert\ryskus{g}_{n,l}^{\calS} \bigl( \widehat{\bgamma
}%
_{n,l}^{\mathrm{OR}},\widehat{\bbeta}_{n},
\widehat{%
\btheta}_{n,-l} \bigr) -\ryskus{g}_{n,l}^{\calS}
\bigl( \widehat{%
\bgamma}_{n,l}^{\mathrm{OR}},\widehat{
\bbeta}_{n},%
\widetilde{\btheta}_{-l0} \bigr)
\bigr\rrVert_{\infty}=O_{p} ( a_{n}+b_{n}),
\]
where $a_{n}=c_{n} ( n_{\statT}\log n_{\statT}/J_{n}^{\calS} )^{1/2}$
and\vspace*{1pt} $b_{n}=c_{n}n_{\statT} ( J_{n}^{\calS} )^{-1/2}J_{n}^{-p}$ with
$c_{n}=\break( \lambda_{n}^{\min} )^{-1} ( \lambda_{n}^{\max} )^{2}$.
Following similar reasoning, one can prove that
$\llVert\ryskus{g}_{n,l}^{\calS} ( \widehat{\bgamma}%
_{n,l}^{\mathrm{OR}},\break\widehat{\bbeta}_{n},\widetilde{%
\btheta}_{-l0} ) -\ryskus{g}_{n,l}^{\calS} ( \widehat{%
\bgamma}_{n,l}^{\mathrm{OR}},\bbeta,\widetilde{%
\btheta}_{-l0} ) \rrVert_{\infty}=O_{p} ( a_{n}+d_{n} )
$, where $d_{n}=c_{n}n_{\statT} ( J_{n}^{\calS} )^{-1/2}J_{n}^{-2p}$,
$\llVert\ryskus{g}_{n,l}^{\calS} (
\widehat{\bgamma}_{n,l}^{\mathrm{OR}},\bbeta,\widetilde{%
\btheta}_{-l0} ) -\ryskus{g}_{n,l}^{\calS} (
\widehat{\bgamma}_{n,l}^{\mathrm{OR}},\bbeta,\btheta%
_{-l} ) \rrVert_{\infty}=O_{p} ( b_{n} ) $, where\vspace*{1pt} $%
\ryskus{g}_{n,l}^{\calS} ( \widehat{\bgamma}_{n,l}^
\mathrm{OR},\bbeta,\btheta_{-l} ) =\ryskus{0}$. Thus, $\llVert
\ryskus{g}_{n,l}^{\calS} ( \widehat{\bgamma}_{n,l}^{\mathrm
{OR}},\widehat{\bbeta}_{n},\widehat{\btheta}%
_{n,-l} ) \rrVert_{\infty}=O_{p} ( a_{n}+b_{n} ) $.
By the Taylor expansion, there is $t\in( 0,1 ) $ such that $%
\widetilde{\bgamma}_{n,l}=t\widehat{\bgamma}_{n,l}^{\mathrm{OR}}+ (
1-t ) \widehat{\bgamma}_{n,l}^{\calS\calS}$,
\[
\widehat{\bgamma}_{n,l}^{\calS\calS}-\widehat{\bgamma
}%
_{n,l}^{\mathrm{OR}}=- \bigl\{ \partial
\ryskus{g}_{n,l}^{\calS} ( \widetilde{\bgamma}_{n,l},
\widehat{\bbeta}_{n},%
\widehat{\btheta}_{n,-l} )
/\partial\widetilde{\bgamma}_{n,l}^{\statT} \bigr
\}^{-1}\ryskus{g}_{n,l}^{\calS} \bigl( \widehat{\bgamma
}_{n,l}^{\mathrm{OR}},\widehat{\bbeta}_{n}%
,
\widehat{\btheta}_{n,-l} \bigr).
\]
$\partial\ryskus{g}_{n,l}^{\calS} ( \widetilde{\bgamma}%
_{n,l},\widehat{\bbeta}_{n},\widehat{\btheta}%
_{n,-l} ) /\partial\widetilde{\bgamma}_{n,l}^{\statT%
}=\Lambda_{n} ( 1+o_{p} ( 1 ) ) $, with $\Lambda_{n}=\sum_{i=1}^{n} (
\ryskus{B}_{i\cdot l}^{\calS} )^{\statT}\widetilde{\Delta
}_{i}\widetilde{\ryskus{V}}_{i}^{-1}\times
\widetilde{\Delta}_{i}\ryskus{B}_{i\cdot l}^{\calS}$, $\widetilde{%
\Delta}_{i}=\Delta_{i} ( \widehat{\bbeta}_{n},%
\widehat{\btheta}_{n,-l},\widetilde{\bgamma}_{n,l} ) $
and $\widetilde{\ryskus{V}}_{i}=\ryskus{V}_{i} ( \widehat{\bbeta%
}_{n},\widehat{\btheta}_{n,-l},\widetilde{\bgamma}%
_{n,l} ) $. There exist constants $0<c_{3}<C_{3}<\infty$, such that
with probability $1$, for $n_{\statT}$ sufficiently large, $%
c_{3}\lambda_{n}^{\min}n_{\statT}\leq\lambda_{\mathrm{min}} ( \Lambda
_{n} ) \leq\lambda_{\mathrm{max}} ( \Lambda_{n} ) \leq C_{3}\lambda
_{n}^{\max}n_{\statT}$. By Theorem\vspace*{1pt} 13.4.3 of~\cite{DL93},
one has $\llVert\Lambda_{n}^{-1}\rrVert_{\infty}=O_{\mathrm{a.s.}} \{ (
\lambda_{n}^{\min}n_{\statT} )^{-1} \} $. Therefore,
\begin{eqnarray*}
\bigl\llVert
\widehat{\bgamma}_{n,l}^{\calS\calS}-\widehat{\bgamma}%
_{n,l}^{\mathrm{OR}}\bigr\rrVert_{\infty}
&\leq&\bigl\llVert\bigl\{ \partial\ryskus{g}_{n,l}^{\calS} (
\widetilde{\bgamma}_{n,l},\widehat{\bbeta}_{n},%
\widehat{\btheta}_{n,-l} ) /\partial\widetilde{\bgamma}_{n,l}^{\statT
} \bigr\}^{-1}\bigr\rrVert_{\infty}\bigl\llVert\ryskus{g}%
_{n,l}^{\calS} \bigl( \widehat{\bgamma}_{n,l}^{\mathrm{OR}},%
\widehat{\bbeta}_{n},\widehat{\btheta}%
_{n,-l} \bigr) \bigr\rrVert_{\infty}\\
&=&O_{p} \bigl\{ \bigl( \lambda_{n}^{\max}/\lambda_{n}^{\min} \bigr)^{2} \bigl( \sqrt{\log
n_{\statT}/ ( J_{n}^{\calS}n_{%
\statT} ) }+ \bigl( J_{n}^{\calS} \bigr)^{-1/2}J_{n}^{-p} \bigr) \bigr\}.\hspace*{15pt}\qed
\end{eqnarray*}
\noqed\end{pf}
\begin{pf*}{Proof of Theorem \protect\ref{THMthetahatconvergence}}
By Lemma~\ref{LEMgammahats},
\begin{eqnarray*}
&&\sup_{z_{l}\in[ 0,1 ] }\bigl\llvert\widehat{\theta}%
_{n,l}^{\calS}
( z_{l},\widehat{\bbeta}_{n}%
,\widehat{\btheta
}_{n,-l} ) -\widehat{\theta}_{n,l}^{%
\calS} (
z_{l},\bbeta_{0},\btheta_{-l0} ) \bigr\rrvert
\\
&&\qquad\leq\sum_{s=1}^{J_{n}^{\calS}}\bigl\llvert
B_{s,l} ( z_{l} ) \bigr\rrvert\bigl\llVert\widehat{
\bgamma}_{n,l}^{%
\calS\calS}-\widehat{\bgamma}_{n,l}^{\mathrm{OR}}
\bigr\rrVert_{\infty}
\\
&&\qquad=O_{p} \bigl\{ \bigl( \lambda_{n}^{\max}/
\lambda_{n}^{\min} \bigr)^{2} \bigl( \sqrt{\log
n_{\statT}/n_{\statT}}%
+J_{n}^{-p}
\bigr) \bigr\}.
\end{eqnarray*}
By the above result and (\ref{EQBGB}),%
\[
\sup_{z_{l}\in[ 0,1 ] }\bigl\llvert\bigl( \ryskus{B}_{l}^{\calS}
( z_{l} )^{\statT}\Xi_{n,l}^{\ast}
\ryskus{B}_{l}^{%
\calS} ( z_{l} )
\bigr)^{-1/2} \bigl\{ \widehat{\theta}%
_{n,l}^{\calS}
( z_{l},\widehat{\bbeta}_{n}%
,\widehat{\btheta
}_{n,-l} ) -\widehat{\theta}_{n,l}^{%
\calS} (
z_{l},\bbeta_{0},\btheta_{-l0} ) \bigr\} \bigr
\rrvert=o_{p} ( 1 ).
\]
Thus, the asymptotic normality of $\widehat{\theta}%
_{n,l}^{\calS} ( z_{l},\widehat{\bbeta}_{n}%
,\widehat{\btheta}_{n,-l} ) $ follows from Theorem~\ref{THMthetahatnormality},
the above result and Slutsky's theorem.
\end{pf*}
\end{appendix}

\section*{Acknowledgments}

The author is grateful for the insightful comments from the Editor, an
Associate Editor and anonymous referees.



\printaddresses


\begin{thebibliography}{36}

\bibitem{B98}
\begin{bbook}[mr]
\bauthor{\bsnm{Bosq},~\bfnm{D.}\binits{D.}}
(\byear{1998}).
\btitle{Nonparametric Statistics for Stochastic Processes: Estimation and Prediction},
\bedition{2nd} ed.
\bseries{Lecture Notes in Statistics}
\bvolume{110}.
\bpublisher{Springer}, \blocation{New York}.
\bid{mr={1640691}}
\bptok{imsref}%
\end{bbook}
\endbibitem

\bibitem{BC05}
\begin{barticle}[mr]
\bauthor{\bsnm{Bun},~\bfnm{Maurice J.~G.}\binits{M.~J.~G.}} \AND
  \bauthor{\bsnm{Carree},~\bfnm{Martin~A.}\binits{M.~A.}}
(\byear{2005}).
\btitle{Bias-corrected estimation in dynamic panel data models}.
\bjournal{J. Bus. Econom. Statist.}
\bvolume{23}
\bpages{200--210}.
\bid{doi={10.1198/073500104000000532}, issn={0735-0015}, mr={2157271}}
\bptok{imsref}%
\end{barticle}
\endbibitem

\bibitem{dB01}
\begin{bbook}[mr]
\bauthor{\bparticle{de} \bsnm{Boor},~\bfnm{Carl}\binits{C.}}
(\byear{2001}).
\btitle{A Practical Guide to Splines},
\bedition{revised} ed.
\bseries{Applied Mathematical Sciences}
\bvolume{27}.
\bpublisher{Springer}, \blocation{New York}.
\bid{mr={1900298}}
\bptok{imsref}%
\end{bbook}
\endbibitem

\bibitem{DL93}
\begin{bbook}[mr]
\bauthor{\bsnm{DeVore},~\bfnm{Ronald~A.}\binits{R.~A.}} \AND
  \bauthor{\bsnm{Lorentz},~\bfnm{George~G.}\binits{G.~G.}}
(\byear{1993}).
\btitle{Constructive Approximation}.
\bseries{Grundlehren der Mathematischen Wissenschaften [Fundamental Principles
  of Mathematical Sciences]}
\bvolume{303}.
\bpublisher{Springer}, \blocation{Berlin}.
\bid{mr={1261635}}
\bptok{imsref}%
\end{bbook}
\endbibitem

\bibitem{HT90}
\begin{bbook}[mr]
\bauthor{\bsnm{Hastie},~\bfnm{T.~J.}\binits{T.~J.}} \AND
  \bauthor{\bsnm{Tibshirani},~\bfnm{R.~J.}\binits{R.~J.}}
(\byear{1990}).
\btitle{Generalized Additive Models}.
\bseries{Monographs on Statistics and Applied Probability}
\bvolume{43}.
\bpublisher{Chapman \& Hall}, \blocation{London}.
\bid{mr={1082147}}
\bptok{imsref}%
\end{bbook}
\endbibitem

\bibitem{HFZ05}
\begin{barticle}[mr]
\bauthor{\bsnm{He},~\bfnm{Xuming}\binits{X.}},
  \bauthor{\bsnm{Fung},~\bfnm{Wing~K.}\binits{W.~K.}} \AND
  \bauthor{\bsnm{Zhu},~\bfnm{Zhongyi}\binits{Z.}}
(\byear{2005}).
\btitle{Robust estimation in generalized partial linear models for clustered
  data}.
\bjournal{J. Amer. Statist. Assoc.}
\bvolume{100}
\bpages{1176--1184}.
\bid{doi={10.1198/016214505000000277}, issn={0162-1459}, mr={2236433}}
\bptok{imsref}%
\end{barticle}
\endbibitem

\bibitem{HS96}
\begin{barticle}[mr]
\bauthor{\bsnm{He},~\bfnm{Xuming}\binits{X.}} \AND
  \bauthor{\bsnm{Shi},~\bfnm{Peide}\binits{P.}}
(\byear{1996}).
\btitle{Bivariate tensor-product {$B$}-splines in a partly linear model}.
\bjournal{J.~Multivariate Anal.}
\bvolume{58}
\bpages{162--181}.
\bid{doi={10.1006/jmva.1996.0045}, issn={0047-259X}, mr={1405586}}
\bptok{imsref}%
\end{barticle}
\endbibitem

\bibitem{H86}
\begin{barticle}[mr]
\bauthor{\bsnm{Heckman},~\bfnm{Nancy~E.}\binits{N.~E.}}
(\byear{1986}).
\btitle{Spline smoothing in a partly linear model}.
\bjournal{J. R. Stat. Soc. Ser. B Stat. Methodol.}
\bvolume{48}
\bpages{244--248}.
\bid{issn={0035-9246}, mr={0868002}}
\bptok{imsref}%
\end{barticle}
\endbibitem

\bibitem{HRWY98}
\begin{barticle}[mr]
\bauthor{\bsnm{Hoover},~\bfnm{Donald~R.}\binits{D.~R.}},
  \bauthor{\bsnm{Rice},~\bfnm{John~A.}\binits{J.~A.}},
  \bauthor{\bsnm{Wu},~\bfnm{Colin~O.}\binits{C.~O.}} \AND
  \bauthor{\bsnm{Yang},~\bfnm{Li-Ping}\binits{L.-P.}}
(\byear{1998}).
\btitle{Nonparametric smoothing estimates of time-varying coefficient models
  with longitudinal data}.
\bjournal{Biometrika}
\bvolume{85}
\bpages{809--822}.
\bid{doi={10.1093/biomet/85.4.809}, issn={0006-3444}, mr={1666699}}
\bptok{imsref}%
\end{barticle}
\endbibitem

\bibitem{HKM06}
\begin{barticle}[mr]
\bauthor{\bsnm{Horowitz},~\bfnm{Joel}\binits{J.}},
  \bauthor{\bsnm{Klemel{\"a}},~\bfnm{Jussi}\binits{J.}} \AND
  \bauthor{\bsnm{Mammen},~\bfnm{Enno}\binits{E.}}
(\byear{2006}).
\btitle{Optimal estimation in additive regression models}.
\bjournal{Bernoulli}
\bvolume{12}
\bpages{271--298}.
\bid{doi={10.3150/bj/1145993975}, issn={1350-7265}, mr={2218556}}
\bptok{imsref}%
\end{barticle}
\endbibitem

\bibitem{HL05}
\begin{barticle}[mr]
\bauthor{\bsnm{Horowitz},~\bfnm{Joel~L.}\binits{J.~L.}} \AND
  \bauthor{\bsnm{Lee},~\bfnm{Sokbae}\binits{S.}}
(\byear{2005}).
\btitle{Nonparametric estimation of an additive quantile regression model}.
\bjournal{J. Amer. Statist. Assoc.}
\bvolume{100}
\bpages{1238--1249}.
\bid{doi={10.1198/016214505000000583}, issn={0162-1459}, mr={2236438}}
\bptok{imsref}%
\end{barticle}
\endbibitem

\bibitem{HM05}
\begin{barticle}[mr]
\bauthor{\bsnm{Horowitz},~\bfnm{Joel~L.}\binits{J.~L.}} \AND
  \bauthor{\bsnm{Mammen},~\bfnm{Enno}\binits{E.}}
(\byear{2004}).
\btitle{Nonparametric estimation of an additive model with a link function}.
\bjournal{Ann. Statist.}
\bvolume{32}
\bpages{2412--2443}.
\bid{doi={10.1214/009053604000000814}, issn={0090-5364}, mr={2153990}}
\bptok{imsref}%
\end{barticle}
\endbibitem

\bibitem{H03}
\begin{barticle}[mr]
\bauthor{\bsnm{Huang},~\bfnm{Jianhua~Z.}\binits{J.~Z.}}
(\byear{2003}).
\btitle{Local asymptotics for polynomial spline regression}.
\bjournal{Ann. Statist.}
\bvolume{31}
\bpages{1600--1635}.
\bid{doi={10.1214/aos/1065705120}, issn={0090-5364}, mr={2012827}}
\bptok{imsref}%
\end{barticle}
\endbibitem

\bibitem{HZZ06}
\begin{barticle}[mr]
\bauthor{\bsnm{Huang},~\bfnm{Jianhua~Z.}\binits{J.~Z.}},
  \bauthor{\bsnm{Zhang},~\bfnm{Liangyue}\binits{L.}} \AND
  \bauthor{\bsnm{Zhou},~\bfnm{Lan}\binits{L.}}
(\byear{2007}).
\btitle{Efficient estimation in marginal partially linear models for
  longitudinal/clustered data using splines}.
\bjournal{Scand. J. Stat.}
\bvolume{34}
\bpages{451--477}.
\bid{doi={10.1111/j.1467-9469.2006.00550.x}, issn={0303-6898}, mr={2368793}}
\bptok{imsref}%
\end{barticle}
\endbibitem

\bibitem{LZ86}
\begin{barticle}[mr]
\bauthor{\bsnm{Liang},~\bfnm{Kung~Yee}\binits{K.~Y.}} \AND
  \bauthor{\bsnm{Zeger},~\bfnm{Scott~L.}\binits{S.~L.}}
(\byear{1986}).
\btitle{Longitudinal data analysis using generalized linear models}.
\bjournal{Biometrika}
\bvolume{73}
\bpages{13--22}.
\bid{doi={10.1093/biomet/73.1.13}, issn={0006-3444}, mr={0836430}}
\bptok{imsref}%
\end{barticle}
\endbibitem

\bibitem{LC00}
\begin{barticle}[mr]
\bauthor{\bsnm{Lin},~\bfnm{Xihong}\binits{X.}} \AND
  \bauthor{\bsnm{Carroll},~\bfnm{Raymond~J.}\binits{R.~J.}}
(\byear{2000}).
\btitle{Nonparametric function estimation for clustered data when the predictor
  is measured without/with error}.
\bjournal{J. Amer. Statist. Assoc.}
\bvolume{95}
\bpages{520--534}.
\bid{doi={10.2307/2669396}, issn={0162-1459}, mr={1803170}}
\bptok{imsref}%
\end{barticle}
\endbibitem

\bibitem{LC01}
\begin{barticle}[mr]
\bauthor{\bsnm{Lin},~\bfnm{Xihong}\binits{X.}} \AND
  \bauthor{\bsnm{Carroll},~\bfnm{Raymond~J.}\binits{R.~J.}}
(\byear{2001}).
\btitle{Semiparametric regression for clustered data}.
\bjournal{Biometrika}
\bvolume{88}
\bpages{1179--1185}.
\bid{doi={10.1093/biomet/88.4.1179}, issn={0006-3444}, mr={1872228}}
\bptok{imsref}%
\end{barticle}
\endbibitem

\bibitem{LWWC04}
\begin{barticle}[mr]
\bauthor{\bsnm{Lin},~\bfnm{Xihong}\binits{X.}},
  \bauthor{\bsnm{Wang},~\bfnm{Naisyin}\binits{N.}},
  \bauthor{\bsnm{Welsh},~\bfnm{Alan~H.}\binits{A.~H.}} \AND
  \bauthor{\bsnm{Carroll},~\bfnm{Raymond~J.}\binits{R.~J.}}
(\byear{2004}).
\btitle{Equivalent kernels of smoothing splines in nonparametric regression for
  clustered/longitudinal data}.
\bjournal{Biometrika}
\bvolume{91}
\bpages{177--193}.
\bid{doi={10.1093/biomet/91.1.177}, issn={0006-3444}, mr={2050468}}
\bptok{imsref}%
\end{barticle}
\endbibitem

\bibitem{L00}
\begin{barticle}[mr]
\bauthor{\bsnm{Linton},~\bfnm{Oliver~B.}\binits{O.~B.}}
(\byear{2000}).
\btitle{Efficient estimation of generalized additive nonparametric regression
  models}.
\bjournal{Econometric Theory}
\bvolume{16}
\bpages{502--523}.
\bid{doi={10.1017/S0266466600164023}, issn={0266-4666}, mr={1790289}}
\bptok{imsref}%
\end{barticle}
\endbibitem

\bibitem{LY10}
\begin{barticle}[mr]
\bauthor{\bsnm{Liu},~\bfnm{Rong}\binits{R.}} \AND
  \bauthor{\bsnm{Yang},~\bfnm{Lijian}\binits{L.}}
(\byear{2010}).
\btitle{Spline-backfitted kernel smoothing of additive coefficient model}.
\bjournal{Econometric Theory}
\bvolume{26}
\bpages{29--59}.
\bid{doi={10.1017/S0266466609090604}, issn={0266-4666}, mr={2587102}}
\bptok{imsref}%
\end{barticle}
\endbibitem

\bibitem{MSW12}
\begin{bmisc}[auto:STB|2012/11/23|13:23:43]
\bauthor{\bsnm{Ma},~\bfnm{S.}\binits{S.}},
  \bauthor{\bsnm{Song},~\bfnm{Q.}\binits{Q.}} \AND
  \bauthor{\bsnm{Wang},~\bfnm{L.}\binits{L.}}
(\byear{2013}).
\bhowpublished{Simultaneous variable selection and estimation in semiparametric
  modeling of longitudinal/clustered data. \textit{Bernoulli}
\textbf{19} 252--274}.
\bptok{imsref}%
\end{bmisc}
\endbibitem

\bibitem{MY11}
\begin{barticle}[mr]
\bauthor{\bsnm{Ma},~\bfnm{Shujie}\binits{S.}} \AND
  \bauthor{\bsnm{Yang},~\bfnm{Lijian}\binits{L.}}
(\byear{2011}).
\btitle{Spline-backfitted kernel smoothing of partially linear additive model}.
\bjournal{J. Statist. Plann. Inference}
\bvolume{141}
\bpages{204--219}.
\bid{doi={10.1016/j.jspi.2010.05.028}, issn={0378-3758}, mr={2719488}}
\bptok{imsref}%
\end{barticle}
\endbibitem

\bibitem{MLN99}
\begin{barticle}[mr]
\bauthor{\bsnm{Mammen},~\bfnm{E.}\binits{E.}},
  \bauthor{\bsnm{Linton},~\bfnm{O.}\binits{O.}} \AND
  \bauthor{\bsnm{Nielsen},~\bfnm{J.}\binits{J.}}
(\byear{1999}).
\btitle{The existence and asymptotic properties of a backfitting projection
  algorithm under weak conditions}.
\bjournal{Ann. Statist.}
\bvolume{27}
\bpages{1443--1490}.
\bid{doi={10.1214/aos/1017939137}, issn={0090-5364}, mr={1742496}}
\bptok{imsref}%
\end{barticle}
\endbibitem

\bibitem{M90}
\begin{bmisc}[auto:STB|2012/11/23|13:23:43]
\bauthor{\bsnm{Munnell},~\bfnm{A.~H.}\binits{A.~H.}}
(\byear{1990}).
\bhowpublished{How does public infrastructure affect regional economic performance.
\textit{New England Econ. Rev.}
Sep. 11--33.}
\bptok{imsref}%
\end{bmisc}
\endbibitem

\bibitem{OR97}
\begin{barticle}[mr]
\bauthor{\bsnm{Opsomer},~\bfnm{Jean~D.}\binits{J.~D.}} \AND
  \bauthor{\bsnm{Ruppert},~\bfnm{David}\binits{D.}}
(\byear{1997}).
\btitle{Fitting a bivariate additive model by local polynomial regression}.
\bjournal{Ann. Statist.}
\bvolume{25}
\bpages{186--211}.
\bid{doi={10.1214/aos/1034276626}, issn={0090-5364}, mr={1429922}}
\bptok{imsref}%
\end{barticle}
\endbibitem

\bibitem{SY10}
\begin{barticle}[mr]
\bauthor{\bsnm{Song},~\bfnm{Qiongxia}\binits{Q.}} \AND
  \bauthor{\bsnm{Yang},~\bfnm{Lijian}\binits{L.}}
(\byear{2010}).
\btitle{Oracally efficient spline smoothing of nonlinear additive
  autoregression models with simultaneous confidence band}.
\bjournal{J. Multivariate Anal.}
\bvolume{101}
\bpages{2008--2025}.
\bid{doi={10.1016/j.jmva.2010.04.004}, issn={0047-259X}, mr={2671198}}
\bptok{imsref}%
\end{barticle}
\endbibitem

\bibitem{W99}
\begin{bmisc}[auto:STB|2012/11/23|13:23:43]
\bauthor{\bsnm{Walterskirchen},~\bfnm{E.}\binits{E.}}
(\byear{1999}).
\bhowpublished{The relationship between growth, employment and unemployment in
  the EU. European economists for an alternative economic policy, Workshop in
  Barcelona}.
\bptok{imsref}%
\end{bmisc}
\endbibitem

\bibitem{WY09}
\begin{barticle}[mr]
\bauthor{\bsnm{Wang},~\bfnm{Jing}\binits{J.}} \AND
  \bauthor{\bsnm{Yang},~\bfnm{Lijian}\binits{L.}}
(\byear{2009}).
\btitle{Polynomial spline confidence bands for regression curves}.
\bjournal{Statist. Sinica}
\bvolume{19}
\bpages{325--342}.
\bid{issn={1017-0405}, mr={2487893}}
\bptok{imsref}%
\end{barticle}
\endbibitem

\bibitem{WY07}
\begin{barticle}[mr]
\bauthor{\bsnm{Wang},~\bfnm{Li}\binits{L.}} \AND
  \bauthor{\bsnm{Yang},~\bfnm{Lijian}\binits{L.}}
(\byear{2007}).
\btitle{Spline-backfitted kernel smoothing of nonlinear additive autoregression
  model}.
\bjournal{Ann. Statist.}
\bvolume{35}
\bpages{2474--2503}.
\bid{doi={10.1214/009053607000000488}, issn={0090-5364}, mr={2382655}}
\bptok{imsref}%
\end{barticle}
\endbibitem

\bibitem{WCL05}
\begin{barticle}[mr]
\bauthor{\bsnm{Wang},~\bfnm{Naisyin}\binits{N.}},
  \bauthor{\bsnm{Carroll},~\bfnm{Raymond~J.}\binits{R.~J.}} \AND
  \bauthor{\bsnm{Lin},~\bfnm{Xihong}\binits{X.}}
(\byear{2005}).
\btitle{Efficient semiparametric marginal estimation for longitudinal/clustered
  data}.
\bjournal{J. Amer. Statist. Assoc.}
\bvolume{100}
\bpages{147--157}.
\bid{doi={10.1198/016214504000000629}, issn={0162-1459}, mr={2156825}}
\bptok{imsref}%
\end{barticle}
\endbibitem

\bibitem{WLC02}
\begin{barticle}[mr]
\bauthor{\bsnm{Welsh},~\bfnm{Alan~H.}\binits{A.~H.}},
  \bauthor{\bsnm{Lin},~\bfnm{Xihong}\binits{X.}} \AND
  \bauthor{\bsnm{Carroll},~\bfnm{Raymond~J.}\binits{R.~J.}}
(\byear{2002}).
\btitle{Marginal longitudinal nonparametric regression: Locality and efficiency
  of spline and kernel methods}.
\bjournal{J. Amer. Statist. Assoc.}
\bvolume{97}
\bpages{482--493}.
\bid{doi={10.1198/016214502760047014}, issn={0162-1459}, mr={1941465}}
\bptok{imsref}%
\end{barticle}
\endbibitem

\bibitem{WY96}
\begin{barticle}[mr]
\bauthor{\bsnm{Wild},~\bfnm{C.~J.}\binits{C.~J.}} \AND
  \bauthor{\bsnm{Yee},~\bfnm{T.~W.}\binits{T.~W.}}
(\byear{1996}).
\btitle{Additive extensions to generalized estimating equation methods}.
\bjournal{J. R. Stat. Soc. Ser. B Stat. Methodol.}
\bvolume{58}
\bpages{711--725}.
\bid{issn={0035-9246}, mr={1410186}}
\bptok{imsref}%
\end{barticle}
\endbibitem

\bibitem{XY03}
\begin{barticle}[mr]
\bauthor{\bsnm{Xie},~\bfnm{Minge}\binits{M.}} \AND
  \bauthor{\bsnm{Yang},~\bfnm{Yaning}\binits{Y.}}
(\byear{2003}).
\btitle{Asymptotics for generalized estimating equations with large cluster
  sizes}.
\bjournal{Ann. Statist.}
\bvolume{31}
\bpages{310--347}.
\bid{doi={10.1214/aos/1046294467}, issn={0090-5364}, mr={1962509}}
\bptok{imsref}%
\end{barticle}
\endbibitem

\bibitem{XY06}
\begin{barticle}[mr]
\bauthor{\bsnm{Xue},~\bfnm{Lan}\binits{L.}} \AND
  \bauthor{\bsnm{Yang},~\bfnm{Lijian}\binits{L.}}
(\byear{2006}).
\btitle{Additive coefficient modeling via polynomial spline}.
\bjournal{Statist. Sinica}
\bvolume{16}
\bpages{1423--1446}.
\bid{issn={1017-0405}, mr={2327498}}
\bptok{imsref}%
\end{barticle}
\endbibitem

\bibitem{ZSW98}
\begin{barticle}[mr]
\bauthor{\bsnm{Zhou},~\bfnm{S.}\binits{S.}},
  \bauthor{\bsnm{Shen},~\bfnm{X.}\binits{X.}} \AND
  \bauthor{\bsnm{Wolfe},~\bfnm{D.~A.}\binits{D.~A.}}
(\byear{1998}).
\btitle{Local asymptotics for regression splines and confidence regions}.
\bjournal{Ann. Statist.}
\bvolume{26}
\bpages{1760--1782}.
\bid{doi={10.1214/aos/1024691356}, issn={0090-5364}, mr={1673277}}
\bptok{imsref}%
\end{barticle}
\endbibitem

\end{thebibliography}
\end{document}